\documentclass[11pt,a4paper]{article}

\usepackage[margin=2.45cm]{geometry}
\usepackage[T1]{fontenc}
\usepackage[utf8]{inputenc}
\usepackage{lmodern}
\usepackage{amsmath,amssymb,amsthm,mathtools}
\usepackage{enumitem}
\usepackage{graphicx}
\usepackage{booktabs}
\usepackage{array}
\usepackage{placeins}
\usepackage{float}

\usepackage[colorlinks=true,linkcolor=blue,citecolor=blue,urlcolor=blue]{hyperref}
\hypersetup{pdftitle={A posteriori error estimates for fully coupled McKean-Vlasov FBSDEs with jumps},
            pdfauthor={Chunrong Feng, Christian Garry},
            pdfkeywords={McKean-Vlasov FBSDEs, jump diffusions, a posteriori error estimates, deep BSDE methods, stochastic numerical analysis}}

\graphicspath{{figs/}}
\newif\ifshowdoi
\showdoitrue
\newcommand{\doi}[1]{\ifshowdoi;\ \href{https://doi.org/#1}{doi:#1}\fi}
\newcommand{\geninput}[1]{\IfFileExists{#1}{\input{#1}}{%
  \errmessage{missing generated input #1 -- run the producing script}}}

\newcommand{\twosidedBandVal}{[0.31,\,1.00]}

\newcommand{\gridUniformLoVal}{0.318}
\newcommand{\gridUniformHiVal}{0.319}
\newcommand{\gridMaxEffVal}{1.000}
\newcommand{\gridMinVal}{20}
\newcommand{\gridMaxVal}{320}

\newcommand{\npathsVal}{200{,}000}

\newcommand{\strongBandVal}{[0.19,\,1.00]}

\newcommand{\weakBandVal}{[0.41,\,1.00]}

\newcommand{\seedSpreadPctVal}{1.2}
\newcommand{\oracleLawLoVal}{0.318}
\newcommand{\particleLawLoVal}{0.289}
\newcommand{\lawCheckGridLoVal}{20}
\newcommand{\lawCheckGridHiVal}{80}
\newcommand{\tauFineVal}{960}
\newcommand{\tauPathsVal}{120{,}000}
\newcommand{\tauGridListVal}{15,30,60,120,240}
\newcommand{\tauSlopeVal}{1.12}
\newcommand{\tauSlopeRVal}{1.12}
\newcommand{\dcrBandVal}{[1.001,\,1.018]}

\newcommand{\regimeBalLoVal}{0.837}
\newcommand{\regimeBalHiVal}{0.907}
\newcommand{\regimeBalAchievedBandVal}{[1.000,\,1.003]}

\newcommand{\regimeEstLoVal}{0.790}
\newcommand{\regimeEstHiVal}{0.824}
\newcommand{\regimeEstAchievedBandVal}{[10.008,\,10.024]}
\newcommand{\regimeEstTargetVal}{10}

\newcommand{\lqinvBandVal}{[0.47,\,1.18]}

\newcommand{\lqinvDpiStartVal}{0.24}

\newcommand{\lqinvDpiEndVal}{2.4{\times}10^{-5}}
\newcommand{\lqinvNseedsVal}{10}

\newcommand{\nonlinBaseBandVal}{[0.30,\,0.73]}
\newcommand{\nonlinActiveBandVal}{[0.31,\,0.70]}

\newcommand{\flBandVal}{[0.37,\,1.00]}

\newcommand{\flLamVal}{(0.15,0.10,0.08,0.06)}

\newcommand{\flGridBody}{%
  $\delta X$ & 0.371 & 0.373 \\
  $\delta Y$ & 0.529 & 0.534 \\
  $\delta Z$ & 0.495 & 0.501 \\
  $\delta U$ & 0.496 & 0.509 \\
  $\delta M$ & 1.000 & 1.000 \\
  all components & 0.540 & 0.546 \\
}

\newcommand{\monoBaseVal}{0.95}
\newcommand{\monoStrongVal}{1.95}
\newcommand{\monoWeakVal}{0.25}
\newcommand{\monoFullLawVal}{0.84}

\newcommand{\sineBandVal}{[0.48,\,1.08]}

\newcommand{\sineMmassVal}{1.7\times10^{-4}}
\newcommand{\sineAlphaVal}{0.4}
\newcommand{\sineBetaTwoVal}{0.5}

\newcommand{\contOffLoVal}{0.494}
\newcommand{\contOffHiVal}{0.509}

\newcommand{\contLambdaVal}{1}

\newcommand{\infOffLoVal}{0.498}
\newcommand{\infOffHiVal}{0.502}

\newcommand{\infLambdaVal}{1}
\newcommand{\infEpsVal}{10^{-3}}

\newcommand{\infTrapEpsVal}{0.05}
\newcommand{\infTrapFacVal}{1.11}
\newcommand{\infTrapMatchedVal}{0.498}
\newcommand{\infTrapMismatchVal}{0.552}
\newcommand{\infTrapCompletedVal}{0.500}

\newcommand{\deepBandVal}{[0.29,\,0.68]}

\newcommand{\deepEpiStartVal}{2.66}
\newcommand{\deepEpiEndVal}{1.8{\times}10^{-4}}
\newcommand{\deepDpiEndVal}{8.5{\times}10^{-5}}
\newcommand{\deepBudgetMaxVal}{3200}
\newcommand{\deepNseedsVal}{10}
\newcommand{\deepHiddenVal}{64}

\newcommand{\lsmcBandVal}{[0.42,\,0.67]}

\newcommand{\lsmcBudgetMaxVal}{30}
\newcommand{\lsmcNseedsVal}{10}
\newcommand{\lsmcNtrainVal}{32{,}768}
\newcommand{\lsmcCYshareVal}{99.98\%}

\newcommand{\stratValDetrendCorrVal}{0.96}
\newcommand{\stratValSpearMedVal}{0.96}

\newcommand{\stratLqinvDetrendCorrVal}{0.94}
\newcommand{\stratLqinvSpearMedVal}{0.95}
\newcommand{\stratLqinvSpearLoVal}{-0.28}
\newcommand{\stratLsmcDetrendCorrVal}{0.98}
\newcommand{\stratLsmcSpearMedVal}{0.94}

\newcommand{\rankGeometryTableBody}{$m<n$ & $X$ & $\operatorname{ran}(G^\top)$ & 0.303 & 0.306 \\
$m<n$ & $X$ & $\ker G$ & 1.015 & 1.020 \\
$m>n$ & $Y$ & $\operatorname{ran}(G)$ & 0.492 & 0.499 \\
$m>n$ & $Y$ & $\ker(G^\top)$ & 0.836 & 0.837 \\}

\newcommand{\rankKerMltnVal}{1.02}
\newcommand{\rankRanMltnVal}{0.30}

\newcommand{\nigOffBandVal}{[0.498,\,0.505]}

\newcommand{\nigTypesVal}{4}

\newcommand{\nigTypeEnergyVal}{1}

\newcommand{\candStatsBCVal}{800{,}000}
\newcommand{\candStatsTableBody}{Neural solver & 90 & 0.021\% & 2.6\% & 4.9\%\\
Nonlinear mean, $\lambda=0$ & 30 & 0.17\% & 2.5\% & 4.3\%\\
Nonlinear mean, $\lambda=1$ & 30 & 0.72\% & 8.8\% & 9.9\%\\
Picard/LSMC & 100 & 0.0033\% & 0.055\% & 0.065\%\\
Mean-field inventory control & 90 & 0.036\% & 4\% & 8.7\%\\}

\newcommand{\ictlNpartVal}{100{,}000}
\newcommand{\ictlEnsVal}{16}
\newcommand{\ictlCellsVal}{90}
\newcommand{\ictlSeedsVal}{10}
\newcommand{\ictlBVal}{1.5{\times}10^{-5}}
\newcommand{\ictlCinvVal}{10.51}
\newcommand{\ictlBandVal}{[0.47,\,1.20]}
\newcommand{\ictlDetrendEDVal}{0.93}

\newcommand{\ictlLateFromVal}{800}
\newcommand{\ictlLateEFacVal}{14}
\newcommand{\ictlLateRegLoVal}{2.8{\times}10^{-5}}
\newcommand{\ictlLateRegHiVal}{4.4{\times}10^{-5}}
\newcommand{\ictlBelowBVal}{9}
\newcommand{\ictlSlackVal}{10}

\newcommand{\ictlSelAgreeVal}{8}
\newcommand{\ictlSelEHitsVal}{5}

\newcommand{\ictlSelMaxVal}{4.1}
\newcommand{\ictlSelDiffVal}{-8.4{\times}10^{-7}}
\newcommand{\ictlSelDiffSeVal}{1.4{\times}10^{-6}}
\newcommand{\ictlRefCountVal}{3}
\newcommand{\ictlRefNpartVal}{300{,}000}
\newcommand{\ictlRefShiftVal}{1.2\%}

\newcommand{\hoLqinvObjLoVal}{0.95}
\newcommand{\hoLqinvEstMedVal}{1.00}

\newcommand{\hoNeurObjLoVal}{0.95}
\newcommand{\hoNeurEstMedVal}{1.00}

\newcommand{\refineWorstDevVal}{0.0012}

\newcommand{\selBTVal}{800{,}000}
\newcommand{\selLqinvFinalWorstVal}{1.00}
\newcommand{\selLqinvFinalHitsVal}{10}
\newcommand{\selLqinvFinalSeedsVal}{10}

\newcommand{\selLqinvEWorstVal}{1.06}
\newcommand{\selLqinvEHitsVal}{9}
\newcommand{\selLqinvESeedsVal}{10}

\newcommand{\selLqinvLWorstVal}{1.80}
\newcommand{\selLqinvLHitsVal}{9}
\newcommand{\selLqinvLSeedsVal}{10}
\newcommand{\selLsmcFinalWorstVal}{1.09}
\newcommand{\selLsmcFinalHitsVal}{0}
\newcommand{\selLsmcFinalSeedsVal}{10}

\newcommand{\selLsmcEWorstVal}{1.01}
\newcommand{\selLsmcEHitsVal}{8}
\newcommand{\selLsmcESeedsVal}{10}
\newcommand{\selValFinalWorstVal}{2.37}
\newcommand{\selValFinalHitsVal}{9}
\newcommand{\selValFinalSeedsVal}{10}

\newcommand{\selValEWorstVal}{1.00}
\newcommand{\selValEHitsVal}{10}
\newcommand{\selValESeedsVal}{10}

\newcommand{\selValLWorstVal}{9.78}
\newcommand{\selValLHitsVal}{9}
\newcommand{\selValLSeedsVal}{10}

\newcommand{\snlSeedsVal}{5}
\newcommand{\snlCellsVal}{30}
\newcommand{\snlBudgetMaxVal}{3200}
\newcommand{\snlBudgetPenultVal}{1600}
\newcommand{\snlEStartVal}{4.3}
\newcommand{\snlEEndVal}{1.5{\times}10^{-3}}
\newcommand{\snlDStartVal}{3.9}
\newcommand{\snlDEndVal}{1.2{\times}10^{-3}}
\newcommand{\snlDPenultVal}{2.2{\times}10^{-3}}
\newcommand{\snlBandVal}{[0.57,\,1.35]}

\newcommand{\snlDetrendVal}{0.880}
\newcommand{\snlFloorVal}{1.9{\times}10^{-4}}
\newcommand{\snlFloorShareVal}{17\%}
\newcommand{\snlRefDxVal}{0.005}
\newcommand{\snlRefLVal}{12}
\newcommand{\snlRefQVal}{64}
\newcommand{\snlRefDirectVal}{3{\times}10^{-9}}
\newcommand{\snlFinerDxVal}{0.0025}
\newcommand{\snlFinerLVal}{16}
\newcommand{\snlFinerQVal}{128}
\newcommand{\snlFinerTailVal}{10^{-16}}
\newcommand{\snlFinerRootVal}{10^{-13}}

\newcommand{\closedFormSettingsRows}{%
base scalar sweep (\S\ref{B:subsec:benchmark-protocol}) & $20,40,80,160,320$ & $200{,}000$ & 1 & $\{0.001,\allowbreak 0.003,\allowbreak 0.01,\allowbreak 0.03,\allowbreak 0.1,\allowbreak 0.3\}$\\
feedback $\kappa_{\mathrm F}=\kappa_{\mathrm B}\in\{0.3,\allowbreak 2\}$ & $20,40,80,160,320$ & $200{,}000$ & 1 each & $\{0.001,\allowbreak 0.003,\allowbreak 0.01,\allowbreak 0.03,\allowbreak 0.1,\allowbreak 0.3\}$\\
particle-mean law input (\S\ref{B:ssec:twosided}) & $20,40,80$ & $200{,}000$ & 1 & $\{0.001,\allowbreak 0.003,\allowbreak 0.01,\allowbreak 0.03,\allowbreak 0.1,\allowbreak 0.3\}$\\
vector rank cases (\S\ref{B:ssec:rank}) & $20,40,80,160$ & $120{,}000$ & 1 & $\{0.001,\allowbreak 0.01,\allowbreak 0.1,\allowbreak 0.3\}$\\
full-tuple mean coupling (\S\ref{B:ssec:fulllaw-setup}) & $20,40,80,160$ & $120{,}000$ & 1 & $\{0.001,\allowbreak 0.01,\allowbreak 0.1,\allowbreak 0.3\}$\\
continuous-time simulation (\S\ref{B:ssec:tau}) & $15,30,60,120,240$; $N_{\mathrm{fine}}=960$ & $120{,}000$ & 1 & unperturbed; targets $r\in\{1,10\}$\\
sine terminal (\S\ref{B:ssec:sine}) & $20,40,80$ & $200{,}000$ & 5 & $\{0.001,\allowbreak 0.01,\allowbreak 0.1\}$\\
continuous marks (\S\ref{B:ssec:contmark}) & $20,40,80,160$ & $100{,}000$ & 1 & $\{0.01,\allowbreak 0.1,\allowbreak 0.3\}$\\
infinite activity, NIG (\S\ref{B:ssec:nigmark}) & $20,40,80,160$ & $120{,}000$ & 1 & $\{0.001,\allowbreak 0.01,\allowbreak 0.1,\allowbreak 0.3\}$\\
infinite activity, cutoff (\S\ref{B:ssec:infmark}) & $20,40,80$ & $100{,}000$ & 1 & $\{0.01,\allowbreak 0.1,\allowbreak 0.3\}$\\
integrand effectivity & $160$ & $100{,}000$ & 3 & target $d=10^{-3}$\\
}

\IfFileExists{results-macros.tex}{\providecommand{\rv}[1]{\csname res@#1\endcsname}
\expandafter\gdef\csname res@apocert.coupling\endcsname{0.2}
\expandafter\gdef\csname res@apocert.G\endcsname{1}
\expandafter\gdef\csname res@apocert.betaOne\endcsname{1}
\expandafter\gdef\csname res@apocert.betaStrong\endcsname{2}
\expandafter\gdef\csname res@apocert.betaWeak\endcsname{0.3}
\expandafter\gdef\csname res@apocert.T\endcsname{1}
\expandafter\gdef\csname res@apocert.sigmaZero\endcsname{0.5}
\expandafter\gdef\csname res@apocert.nuE\endcsname{1}
\expandafter\gdef\csname res@apocert.xZeroLaw\endcsname{\mathcal N(0,1)}
\expandafter\gdef\csname res@vector.coupling\endcsname{0.15}
\expandafter\gdef\csname res@vector.betaOne\endcsname{1}
\expandafter\gdef\csname res@vector.G.mltn.11\endcsname{1}
\expandafter\gdef\csname res@vector.G.mltn.12\endcsname{0.5}
\expandafter\gdef\csname res@vector.sminSq.mltn\endcsname{1.25}
\expandafter\gdef\csname res@vector.cG.mltn\endcsname{0.8}
\expandafter\gdef\csname res@vector.G.meqn.11\endcsname{1}
\expandafter\gdef\csname res@vector.G.meqn.12\endcsname{0.3}
\expandafter\gdef\csname res@vector.G.meqn.21\endcsname{0.2}
\expandafter\gdef\csname res@vector.G.meqn.22\endcsname{1}
\expandafter\gdef\csname res@vector.smin.meqn\endcsname{0.75125}
\expandafter\gdef\csname res@vector.cG.meqn\endcsname{1.7719}
\expandafter\gdef\csname res@vector.G.mgtn.11\endcsname{1}
\expandafter\gdef\csname res@vector.G.mgtn.21\endcsname{0.5}
\expandafter\gdef\csname res@vector.sminSq.mgtn\endcsname{1.25}
\expandafter\gdef\csname res@vector.cG.mgtn\endcsname{0.8}
\expandafter\gdef\csname res@vector.noiseEntry\endcsname{0.5}
\expandafter\gdef\csname res@vector.xZeroLaw\endcsname{\mathcal N(0,I_n)}
\expandafter\gdef\csname res@fulllaw.xZeroLaw\endcsname{\mathcal N(0,1)}
\expandafter\gdef\csname res@fulllaw.T\endcsname{1}
\expandafter\gdef\csname res@fulllaw.betaOne\endcsname{1}
\expandafter\gdef\csname res@fulllaw.sigmaZero\endcsname{0.5}
\expandafter\gdef\csname res@fulllaw.nuE\endcsname{1}
\expandafter\gdef\csname res@sine.a\endcsname{1}
\expandafter\gdef\csname res@sine.csin\endcsname{0.3}
\expandafter\gdef\csname res@sine.k\endcsname{2}
\expandafter\gdef\csname res@sine.sigmaZero\endcsname{0.4}
\expandafter\gdef\csname res@sine.hZero\endcsname{0.3}
\expandafter\gdef\csname res@sine.nuE\endcsname{1}
\expandafter\gdef\csname res@sine.T\endcsname{1}
\expandafter\gdef\csname res@sine.xZeroLaw\endcsname{\mathcal N(0.3,0.5^2)}
\expandafter\gdef\csname res@select.eWorstPct\endcsname{6}
}{%
  \errmessage{missing results-macros.tex -- run scripts/gen_macros.py}}
\newcommand{\incfig}[2][\linewidth]{%
  \IfFileExists{figs/#2}{\includegraphics[width=#1]{#2}}{%
    \errmessage{missing figure figs/#2 -- run the plotting script}}}

\newcommand{\R}{\mathbb{R}}
\newcommand{\E}{\mathbb{E}}
\newcommand{\Prob}{\mathbb{P}}
\newcommand{\F}{\mathcal{F}}
\newcommand{\Ei}{\E_i}
\newcommand{\diff}{\mathrm{d}}
\newcommand{\Nt}{\widetilde{N}}                 %
\newcommand{\inner}[2]{\langle #1,\,#2\rangle}
\newcommand{\Lnu}{L^2(\nu)}
\newcommand{\Gstar}{G^{\!\top}}                 %
\newcommand{\Law}{\mathrm{Law}}
\newcommand{\Sc}[1]{\mathcal{S}^2(0,T;#1)}       %
\newcommand{\HWc}[1]{\mathcal{H}^2_{W}(0,T;#1)}  %
\newcommand{\HNc}[1]{\mathcal{H}^2_{N}(0,T;#1)}  %
\newcommand{\Sspi}[1]{\mathcal{S}^2(\pi;#1)}     %
\newcommand{\HWpi}[1]{\mathcal{H}^2_{W}(\pi;#1)} %
\newcommand{\HNpi}[1]{\mathcal{H}^2_{N}(\pi;#1)} %
\newcommand{\Mpi}[1]{\mathcal{M}^2(\pi;#1)}      %
\newcommand{\ap}[1]{\widehat{#1}}               %
\newcommand{\ex}[1]{#1^{\pi}}                    %
\newcommand{\au}[1]{\bar{#1}}                    %
\newcommand{\sm}[1]{\check{#1}}                  %
\newcommand{\Bc}{\mathcal{B}}                    %
\newcommand{\EG}{\mathcal{E}_G}                  %
\newcommand{\Epi}{\mathcal{E}_\pi}              %
\newcommand{\Dpi}{\mathcal{D}_\pi}              %
\newcommand{\Epip}{\mathcal{E}_{\pi,p}}          %
\newcommand{\Dpip}{\mathcal{D}_{\pi,p}}          %
\newcommand{\Dc}{\mathcal{D}_{\mathrm c}}        %
\newcommand{\Sd}{\mathcal{S}}                    %
\newcommand{\Pp}{\mathcal{P}}                    %
\newcommand{\Reg}{\mathfrak{R}}                  %
\DeclareMathOperator{\tr}{tr}

\numberwithin{equation}{section}
\theoremstyle{plain}
\newtheorem{theorem}{Theorem}[section]
\newtheorem{proposition}[theorem]{Proposition}
\newtheorem{lemma}[theorem]{Lemma}
\newtheorem{corollary}[theorem]{Corollary}
\theoremstyle{definition}
\newtheorem{assumption}{Assumption}

\newtheorem{definition}[theorem]{Definition}

\newtheorem{convention}[theorem]{Convention}

\theoremstyle{remark}
\newtheorem{remark}[theorem]{Remark}

\DeclareRobustCommand{\samevaluetext}[5]{%
  \begingroup
  \edef\cmpvalone{#1}\edef\cmpvaltwo{#2}\edef\cmpvalthree{#3}%
  \ifx\cmpvalone\cmpvaltwo
    \ifx\cmpvalone\cmpvalthree
      \gdef\cmpvalout{#4}%
    \else
      \gdef\cmpvalout{#5}%
    \fi
  \else
    \gdef\cmpvalout{#5}%
  \fi
  \endgroup\cmpvalout}
\DeclareRobustCommand{\slopephrase}[3]{%
  \samevaluetext{#2}{#3}{#2}%
    {have #1 log--log slope #2}%
    {have #1 log--log slopes #2 and #3, respectively}}

\DeclareRobustCommand{\tauSlopeFitted}{\slopephrase{fitted}{\tauSlopeVal}{\tauSlopeRVal}}

\DeclareRobustCommand{\hoEstCorrPhrase}{%
  \samevaluetext{\hoNeurEstMedVal}{\hoLqinvEstMedVal}{\hoNeurEstMedVal}%
    {is $\hoNeurEstMedVal$ for the estimator in both neural studies}%
    {is $\hoNeurEstMedVal$ and $\hoLqinvEstMedVal$ for the estimator in the linear and inventory
     neural studies, respectively}}

\DeclareRobustCommand{\hoObjLoPhrase}{%
  \samevaluetext{\hoNeurObjLoVal}{\hoLqinvObjLoVal}{\hoNeurObjLoVal}%
    {rounds to $\hoNeurObjLoVal$ in each study}%
    {rounds to $\hoNeurObjLoVal$ and $\hoLqinvObjLoVal$ in the linear and inventory studies,
     respectively}}

\title{\emph{A posteriori} error estimates for fully coupled\\
McKean--Vlasov FBSDEs with jumps}
\author{Chunrong Feng and Christian Garry$^{*}$\\[4pt]
\small Department of Mathematical Sciences, Durham University, Durham DH1 3LE, UK\\[2pt]
\small $^{*}$Corresponding author: \texttt{christian.t.garry@durham.ac.uk}}
\newcommand{\builddate}{November 2026}
\date{\builddate}

\begin{document}
\maketitle

\begin{abstract}
We prove a two-sided \emph{a posteriori} error estimate for numerical approximations of
fully coupled McKean--Vlasov forward--backward SDEs with jumps under full-tuple law
dependence. The coefficients may depend on the joint law of $(X,Y,Z,U)$, and we allow
infinite jump activity. We construct a computable estimator $\Epi$ for the squared full-tuple
error $\Dpi$ relative to the exact Euler solution. The estimator uses
the boundary and cumulative residuals of any admissible grid approximation. Under monotonicity and on
sufficiently fine grids, $c\,\Epi \le \Dpi \le C\,\Epi$. The efficiency bound holds without
monotonicity or a mesh restriction.
Analogous two-sided bounds hold for the continuous-time error up to a discretisation
remainder that vanishes under coefficient time regularity.
Closed-form and state-nonlinear benchmarks cover full-tuple law dependence,
mark-dependent jump integrands, nonzero orthogonal martingales and infinite jump activity.
In our experiments, the estimator tracks the discrete error during neural training and across
Picard iterations. In a mean-field inventory model, we use the estimator and a computable
discretisation term to bound the excess cost of the computed control over the discrete optimum.
\end{abstract}

\smallskip
\noindent\textbf{Keywords:} McKean--Vlasov FBSDEs; jump diffusions; \emph{a posteriori} error estimates; deep BSDE methods; stochastic numerical analysis.

\noindent\textbf{MSC 2020:} 65C30; 60H35; 60H10.

\section{Introduction}
\label{B:sec:intro}

For fully coupled McKean--Vlasov forward--backward stochastic differential equations with jumps
(MV-FBSDEJs), the exact Euler solution is generally unavailable. We use the boundary and
cumulative residuals of a computed approximation to bound its error from above and below.
MV-FBSDEJs arise in mean-field control and games and in probabilistic representations on
Wasserstein space. Jumps can model execution shocks, defaults and impulse-type
interactions. We consider jump systems driven by a
Brownian motion $W$ and an independent compensated Poisson random measure $\Nt$ with
$\sigma$-finite intensity $\diff t\,\nu(\diff e)$ on a countably generated mark space:
\begin{equation}\label{B:eq:mvfbsdej}
\left\{\;\begin{aligned}
X_t&=\xi_0+\!\int_0^t\! b(s,\Theta_s,\mu_s)\diff s+\!\int_0^t\!\sigma(s,\Theta_s,\mu_s)\diff W_s
   +\!\int_0^t\!\!\!\int_E\! h(s,\Theta_s^-,e,\mu_s^-)\Nt(\diff s,\diff e),\\
Y_t&=g(X_T,\Law(X_T))+\!\int_t^T\! f(s,\Theta_s,\mu_s)\diff s-\!\int_t^T\! Z_s\diff W_s
   -\!\int_t^T\!\!\!\int_E\! U_s(e)\Nt(\diff s,\diff e).
\end{aligned}\right.
\end{equation}
Here $\Theta_t:=(X_t,Y_t,Z_t,U_t)$ and $\mu_t:=\Law(\Theta_t)$, with $U_t$
taking values in $\Lnu$. The jump coefficient uses the predictable tuple
$\Theta_s^-:=(X_{s-},Y_{s-},Z_s,U_s)$ and its law $\mu_s^-:=\Law(\Theta_s^-)$.

The monotonicity hypothesis of Li and Min~\cite[(H3.2)]{LiMin} requires finite jump activity.
Our earlier work~\cite{CWPpaper} establishes well-posedness for~\eqref{B:eq:mvfbsdej},
allowing infinite jump activity. Here we develop the corresponding discrete stability and
\emph{a posteriori} error theory.

To approximate \eqref{B:eq:mvfbsdej}, fix a uniform grid $t_i=i\tau$, $0\le i\le N$, with
$\tau=T/N$, and consider the discrete tuples
\[
\begin{aligned}
\ex\Theta&=(\ex X_i,\ex Y_i,\ex Z_i,\ex U_i,\ex M_i),\\
\ap\Theta&=(\ap X_i,\ap Y_i,\ap Z_i,\ap U_i,\ap M_i).
\end{aligned}
\]
Here $\ex\Theta$ is the exact solution of the forward--backward Euler scheme, and
$\ap\Theta$ is an admissible approximation on the same grid.
The component $\ex M$ is the orthogonal martingale in the discrete Kunita--Watanabe decomposition.

The true discrete error of $\ap\Theta$ is its squared distance from the exact Euler solution,
\begin{equation}\label{B:eq:Dpi}
\begin{aligned}
\Dpi:={}&\max_{0\le i\le N}\E|\ap X_i-\ex X_i|^2+\max_{0\le i\le N}\E|\ap Y_i-\ex Y_i|^2\\
&+\sum_{i=0}^{N-1}\tau\Big(\E|\ap Z_i-\ex Z_i|_F^2+\E\|\ap U_i-\ex U_i\|_{\Lnu}^2\Big)
+\E|\ap M_N-\ex M_N|^2 .
\end{aligned}
\end{equation}
To estimate $\Dpi$ without computing the exact Euler solution $\ex\Theta$, we use
\begin{equation}\label{B:eq:intro-epi}
\Epi:=\E\big|\ap X_0-\xi_0\big|^2+\E\big|\ap Y_N-g(\ap X_N,\Law(\ap X_N))\big|^2
+\max_{i<N}\E\big|C_i^X\big|^2+\max_{i<N}\E\big|C_i^Y\big|^2,
\end{equation}
where $C_i^X$ and $C_i^Y$ are the cumulative forward and backward residuals defined in
\eqref{B:eq:CX} and \eqref{B:eq:CY}, respectively. Theorem~\ref{B:thm:discrete} shows that, on
sufficiently fine grids,
\begin{equation}\label{B:eq:intro-equiv}
c\,\Epi \le \Dpi \le C\,\Epi .
\end{equation}

\subsection{Related work and positioning}

Bender and Steiner~\cite{BenderSteiner} established computable upper and lower
\emph{a posteriori} bounds for the squared $L^2$-error of approximations to decoupled Brownian
BSDEs, both against the implicit Euler scheme and, up to discretisation terms, against the
continuous-time solution. Their discrete estimate also applies to the decoupled jump setting
of Bouchard and Elie~\cite{BouchardElie}, as observed in~\cite[Remark~2.3]{BenderSteiner}.
Reisinger, Stockinger and Zhang~\cite{RSZ} extended their estimator to fully coupled Brownian
McKean--Vlasov forward--backward SDEs, using cumulative residuals, discrete perturbation
stability and an orthogonal martingale completion. On sufficiently fine grids they obtain
two-sided bounds against the Euler scheme~\cite[Theorems~3.1--3.2]{RSZ}. For approximations
with $\ap M=0$, they also bound the continuous-time error from both sides up to a
discretisation remainder~\cite[Theorem~4.3]{RSZ}.

For Brownian coupled FBSDEs, the deep-BSDE convergence
analyses~\cite{HanLong,HuangNO,NegyesiHO} bound the error in terms of the terminal loss and a
time-discretisation term. For the deep-BSDE trajectories in~\cite[Corollary~4.5]{RSZ}, the
interior equations hold by construction, so the residual estimator reduces to the terminal
loss. Zhang~\cite{ZhangControlMismatch} includes a control-mismatch term when the forward
equation uses an auxiliary approximation of the backward variable.

Wang et al.~\cite{Wang} consider coupled FBSDEs with jumps whose forward coefficients may
depend on $Y$ and whose coefficients have no law dependence. Under their coupling conditions,
they bound the error of a deep-learning approximation by its terminal loss and a
time-discretisation term. The jump component of the bound measures the error in a weighted
integral of the jump integrand. Gnoatto, Oberpriller and Picarelli~\cite{GnoattoJumps} analyse
a variant of the deep BSDEJ solver of~\cite{AnderssonJumps} for decoupled FBSDEs with finite
and infinite jump activity. Their one-sided \emph{a posteriori} bounds control the
$\nu$-integrated jump term.

Meherrem, Guluzada and \c{C}etin~\cite{MeherremGuluzadaCetin} treat a shooting formulation of
fully coupled McKean--Vlasov FBSDEs with finite-activity jumps. Their continuous-time
bound~\cite[Theorem~5.4]{MeherremGuluzadaCetin} controls the error, including the full $\Lnu$
error of the jump integrand, by the terminal loss of the shooting problem. They observe that
the same argument applies to the full-tuple, infinite-activity setting of~\cite{CWPpaper}
under their monotonicity condition. In their shooting formulation, the initial backward value
and the two integrands determine $X$ and $Y$ through the interior equations. Our estimator
measures the cumulative residuals of both equations for a supplied grid approximation. To the
best of our knowledge, no previous work gives a two-sided \emph{a posteriori} estimate for
arbitrary admissible grid approximations of fully coupled McKean--Vlasov FBSDEs with jumps.

\subsection{Contributions and structure}

\begin{enumerate}[label=\textup{(\Alph*)},leftmargin=2.4em,itemsep=3pt]
\item \emph{Discrete theory.} Theorem~\ref{B:thm:discrete} gives the two-sided discrete
estimate for systems with full-tuple law dependence and an $\Lnu$-valued jump integrand.
Efficiency holds under Lipschitz continuity without a mesh restriction; reliability
additionally requires integrability at the origin, $G$-monotonicity and a sufficiently
fine grid. We also establish discrete perturbation
stability and well-posedness of the Euler scheme.
\item \emph{Continuous theory.} Theorem~\ref{B:thm:cont} is the jump analogue of
\cite[Theorem~4.3]{RSZ}: it extends the two-sided estimate to the continuous-time error,
including the full $\Lnu$ error of the jump integrand. The discretisation term vanishes
under coefficient time
regularity. The squared-error discretisation term is $O(\tau)$ if the coefficient time
modulus is $O(\tau^{1/2})$ and the path-regularity modulus is $O(\tau)$.
\item \emph{Numerical validation.} Sections~\ref{B:sec:validation} and~\ref{B:sec:results}
test the estimator on closed-form benchmarks, a state-nonlinear problem with an
independently computed Euler reference, and a mean-field inventory model. Together these
cover full-tuple law dependence, mark-dependent jump integrands, nonzero orthogonal
martingales and infinite jump activity. In the neural and Picard/LSMC experiments, the
estimator tracks the discrete error. For the inventory model, Proposition~\ref{B:prop:regret}
bounds control regret in terms of the estimator and a computable discretisation baseline.
\end{enumerate}

Section~\ref{B:sec:model} introduces the model and states the standing assumptions.
Section~\ref{B:sec:scheme} constructs the one-step Brownian--Poisson projection.
Following the discrete-stability approach of~\cite{RSZ}, Section~\ref{B:sec:main}
proves reliability by subtracting the cumulative residuals from the candidate's
$X$- and $Y$-components and applying Proposition~\ref{B:prop:pert} to the corrected
tuple. Section~\ref{B:sec:contres} compares the Euler solution with the continuous
solution sampled on the grid, using conditional cell projections for $Z$ and $U$.

\section{Model and assumptions}
\label{B:sec:model}
\subsection{Stochastic basis and solution spaces}\label{B:ssec:basis}
On $[0,T]$ we work on a filtered probability space $(\Omega,\F,(\F_t)_{t\ge0},\Prob)$
supporting a $d$-dimensional Brownian motion $W$, $d\ge1$, and a Poisson random measure
$N(\diff t,\diff e)$ on $[0,T]\times E$ independent of $W$.
Its intensity measure is $\diff t\,\nu(\diff e)$ over a $\sigma$-finite, countably generated mark
space $(E,\mathcal E,\nu)$, so that $L^2(\nu;\R^m)$ is separable; the compensated measure is
$\Nt(\diff t,\diff e):=N(\diff t,\diff e)-\diff t\,\nu(\diff e)$,
see~\cite[\S~1, p.~58]{BBP}.
We take $(\F_t)_{t\ge0}$ to be the usual augmentation of the filtration generated by $W$, $N$
and an independent initial $\sigma$-algebra $\F_0$, which carries any auxiliary randomness
used to construct the approximation. The resulting filtration retains the Brownian--Poisson predictable
representation property~\cite[\S5.3]{Applebaum}, \cite[Theorem~13.49]{HeWangYan}.
Assumptions~\ref{B:ass:integ}--\ref{B:ass:mono} below are imposed on this enlarged basis.

In \eqref{B:eq:mvfbsdej}, we take
\[
X_t\in\R^n,\qquad Y_t\in\R^m,\qquad Z_t\in\R^{m\times d},\qquad U_t\in L^2(\nu;\R^m).
\]

Set
\[
\mathcal H:=\R^n\times\R^m\times\R^{m\times d}\times L^2(\nu;\R^m),
\]
where $|\cdot|$ is the Euclidean norm (on any $\R^j$) and the Frobenius and
$L^2(\nu;\R^m)$ norms are
\[
|z|_F^2:=\sum_{i,j}z_{ij}^2,
\qquad
\|u\|_{\Lnu}^2:=\int_E|u(e)|^2\,\nu(\diff e).
\]
We give $\mathcal H$ the Hilbert product norm
\[
\|\theta\|_{\mathcal H}^2:=|x|^2+|y|^2+|z|_F^2+\|u\|_{\Lnu}^2,
\qquad\theta=(x,y,z,u)\in\mathcal H.
\]
Throughout, $W_2$ denotes the 2-Wasserstein distance induced by this norm on
$\mathcal P_2(\mathcal H)$, and by the Euclidean norm on $\mathcal P_2(\R^n)$.

\paragraph{Continuous solution class.} Write
\[
\begin{aligned}
\Sc{\R^j}&:=\bigl\{V\ \text{c\`adl\`ag adapted}:\ \E\sup_{t\le T}|V_t|^2<\infty\bigr\},\\[4pt]
\HWc{\R^{m\times d}}&:=\bigl\{Z\ \text{predictable}:\ \E\!\int_0^T\!|Z_t|_F^2\diff t<\infty\bigr\},\\[4pt]
\HNc{L^2(\nu;\R^m)}&:=\bigl\{U\ \text{predictable}:\ \E\!\int_0^T\!\|U_t\|_{\Lnu}^2\diff t<\infty\bigr\}.
\end{aligned}
\]
Set
\[
\mathcal S:=\Sc{\R^n}\times\Sc{\R^m}\times\HWc{\R^{m\times d}}
\times\HNc{L^2(\nu;\R^m)} .
\]
A \emph{solution} is an element $\Theta\in\mathcal S$ satisfying \eqref{B:eq:mvfbsdej}.

We use predictable representatives of $Z$ and $U$ and Borel representatives of the law flows
$\mu$ and $\mu^-$, as in~\cite[Lemma~A.1]{CWPpaper}.

\subsection{Coefficients and standing assumptions}\label{B:subsec:standing}

The coefficients are defined on the following domains:
\[
\begin{alignedat}{2}
b&:\ [0,T]\times\mathcal H\times\mathcal P_2(\mathcal H)\to\R^n,
&\qquad f&:\ [0,T]\times\mathcal H\times\mathcal P_2(\mathcal H)\to\R^m,\\
\sigma&:\ [0,T]\times\mathcal H\times\mathcal P_2(\mathcal H)\to\R^{n\times d},
& h&:\ [0,T]\times\mathcal H\times\mathcal P_2(\mathcal H)\to L^2(\nu;\R^n),\\
g&:\ \R^n\times\mathcal P_2(\R^n)\to\R^m.
\end{alignedat}
\]
We write $h(t,\theta,e,\mu):=h(t,\theta,\mu)(e)$.
The coefficients are deterministic and jointly Borel measurable.

\begin{definition}[Diagonal input]\label{B:def:diag}
Fix $t\in[0,T]$. Define
\[
\mathcal D_t
:=
\Bigl\{
(\Theta,\mu)\in
L^2(\Omega,\F_t;\mathcal H)\times\mathcal P_2(\mathcal H)
:
\mu=\Law(\Theta)
\Bigr\}.
\]
An element of $\mathcal D_t$ is called a \emph{diagonal input at time $t$}. Interior
coefficient differences are evaluated between pairs in $\mathcal D_t$.
\end{definition}

\smallskip\noindent\textbf{Difference notation.}
For $(\Theta,\mu),(\Theta',\mu')\in\mathcal D_t$, $\delta$ denotes the difference
\[
\begin{aligned}
\delta\Theta&:=\Theta-\Theta',\quad\text{componentwise }\delta X,\delta Y,\delta Z,\delta U,\\[2pt]
\delta\psi&:=\psi(t,\Theta,\mu)-\psi(t,\Theta',\mu'),\quad \psi\in\{b,\sigma,f\},\\[2pt]
\delta h(e)&:=h(t,\Theta,e,\mu)-h(t,\Theta',e,\mu'),
\end{aligned}
\]
and, for $X_T,X'_T\in L^2(\Omega,\F_T;\R^n)$,
\[
\delta X_T:=X_T-X'_T,\qquad
\delta g_T:=g(X_T,\Law(X_T))-g(X'_T,\Law(X'_T)).
\]
Whenever an interior assumption is stated $\diff t$-a.e., the exceptional null set is
understood to be independent of the pair in $\mathcal D_t$.

\begin{assumption}[Integrability at the origin]\label{B:ass:integ}
$\xi_0\in L^2(\Omega,\F_0;\R^n)$ and, with $\delta_0$ the Dirac mass at the relevant origin,
the origin values satisfy
\[
\begin{aligned}
b(\cdot,0,\delta_0)\ &\in\ L^2(0,T;\R^n),
&\qquad\quad \sigma(\cdot,0,\delta_0)\ &\in\ L^2(0,T;\R^{n\times d}),\\[2.2ex]
h(\cdot,0,\cdot,\delta_0)\ &\in\ L^2\bigl(0,T;L^2(\nu;\R^n)\bigr),
&\qquad\quad f(\cdot,0,\delta_0)\ &\in\ L^2(0,T;\R^m).
\end{aligned}
\]
\end{assumption}

\begin{assumption}[Lipschitz coefficients]\label{B:ass:lip}
There is $L\ge0$ such that, for $\diff t$-a.e.\ $t\in[0,T]$ and all
$(\Theta,\mu),(\Theta',\mu')\in\mathcal D_t$, the following holds $\Prob$-a.s.:
\[
|\delta b|+|\delta\sigma|_F+\|\delta h\|_{\Lnu}+|\delta f|
\le L\big(|\delta X|+|\delta Y|+|\delta Z|_F+\|\delta U\|_{\Lnu}+W_2(\mu,\mu')\big).
\]
For all $X_T,X'_T\in L^2(\Omega,\F_T;\R^n)$, the following holds $\Prob$-a.s.:
\[
|\delta g_T|\le L\big(|\delta X_T|+W_2(\Law(X_T),\Law(X'_T))\big).
\]
\end{assumption}

\begin{assumption}[Coupling matrix]\label{B:ass:G}
A fixed $G\in\R^{m\times n}$ of full rank couples the forward and backward errors, and
$\Gstar\in\R^{n\times m}$ denotes its Euclidean adjoint. Write $\|G\|$ for the operator norm, so that
$|Gv|\le\|G\|\,|v|$ and $|\Gstar v|\le\|G\|\,|v|$, and $s_{\min}(G)>0$ for its smallest
singular value; set $c_G:=s_{\min}(G)^{-2}$. Then
\[
|x|^2\le c_G|Gx|^2\ \ (m\ge n),\qquad |y|^2\le c_G|\Gstar y|^2\ \ (n\ge m).
\]
If $\Gstar$ is injective, the latter also holds for matrices and for elements of $\Lnu$:
\[
|Z|_F^2\le c_G|\Gstar Z|_F^2,
\qquad
\|U\|_{\Lnu}^2\le c_G\|\Gstar U\|_{\Lnu}^2.
\]
\end{assumption}

\begin{assumption}[Terminal and interior jump-extended $G$-monotonicity]\label{B:ass:mono}
For all $X_T,X'_T\in L^2(\Omega,\F_T;\R^n)$,
\[
\E\inner{\delta g_T}{G\delta X_T}\ge\alpha\,\E|G\delta X_T|^2,
\]
and, for $\diff t$-a.e.\ $t\in[0,T]$ and all $(\Theta,\mu),(\Theta',\mu')\in\mathcal D_t$,
\[
\begin{aligned}
\E\Big[&\inner{\delta b}{\Gstar\delta Y}+\inner{\delta\sigma}{\Gstar\delta Z}_F
+\inner{-\delta f}{G\delta X}+\!\int_E\!\inner{\delta h(e)}{\Gstar\delta U(e)}\nu(\diff e)\Big]\\
&\le-\beta_1\E\big[|\Gstar\delta Y|^2+|\Gstar\delta Z|_F^2+\|\Gstar\delta U\|_{\Lnu}^2\big]
-\beta_2\E|G\delta X|^2.
\end{aligned}
\]
The constants satisfy
\[
\alpha,\beta_1,\beta_2\ge0,\qquad \alpha+\beta_1>0,\qquad \beta_1+\beta_2>0,
\]
with $\beta_1>0$ if $m<n$ and $\alpha,\beta_2>0$ if $m>n$.
\end{assumption}

\begin{lemma}[Synchronous coupling]\label{B:lem:law}
For $(\Theta,\mu),(\Theta',\mu')\in\mathcal D_t$,
\[
W_2(\mu,\mu')^2\le\E\big[|\delta X|^2+|\delta Y|^2+|\delta Z|_F^2+\|\delta U\|_{\Lnu}^2\big],
\]
and likewise $W_2(\Law(X_T),\Law(X'_T))^2\le\E|\delta X_T|^2$.
\end{lemma}
\begin{proof}
The joint law of $(\Theta,\Theta')$ is a coupling of $\mu$ and $\mu'$. Since $W_2(\mu,\mu')^2$
is the infimum over all coupling costs, it is bounded above by
$\E\|\Theta-\Theta'\|_{\mathcal H}^2$. Applying the same argument to $(X_T,X_T')$ gives the
terminal bound.
\end{proof}

\section{Discrete scheme and estimator}
\label{B:sec:scheme}

Fix the uniform grid
\[
\pi:=\{t_i:=i\tau\}_{i=0}^N,\qquad \tau:=T/N,
\]
and write
\[
\Delta W_i:=W_{t_{i+1}}-W_{t_i},\qquad
\Delta\Nt_i(\diff e):=\Nt((t_i,t_{i+1}],\diff e),\qquad
\Ei[\,\cdot\,]:=\E[\,\cdot\mid\F_{t_i}].
\]
The \emph{exact Euler solution} satisfies, for $i=0,\dots,N-1$,
\begin{align}
\Delta \ex X_i&=b_i\tau+\sigma_i\Delta W_i+\int_E h_i\,\Delta\Nt_i, & \ex X_0&=\xi_0,\label{B:eq:dfwd}\\
\Delta \ex Y_i&=-f_i\tau+\ex Z_i\Delta W_i+\int_E \ex U_i\,\Delta\Nt_i+\Delta \ex M_i, & \ex Y_N&=g(\ex X_N,\Law(\ex X_N)),\label{B:eq:dbwd}
\end{align}
where $\Delta A_i:=A_{i+1}-A_i$.
The Euler coefficients are evaluated at the left endpoints $t_i$:
\[
\begin{aligned}
\ex\mu_i
&:=\Law(\ex X_i,\ex Y_i,\ex Z_i,\ex U_i),\\[2pt]
\psi_i
&:=\psi\bigl(t_i,\ex X_i,\ex Y_i,\ex Z_i,\ex U_i,\ex\mu_i\bigr),
&&\psi\in\{b,\sigma,f\},\\[2pt]
h_i(e)
&:=h\bigl(t_i,\ex X_i,\ex Y_i,\ex Z_i,\ex U_i,e,\ex\mu_i\bigr).
\end{aligned}
\]

\begin{convention}[Coefficient representatives at grid nodes]\label{B:conv:nodewise}
We fix Borel representatives of $b,\sigma,h,f$ satisfying the Lipschitz inequality in
\ref{B:ass:lip} at every grid node. When \ref{B:ass:mono} is assumed, these
representatives also satisfy its interior monotonicity inequality at every grid node. The
constants are independent of the grid.
\end{convention}

\noindent Such representatives can be obtained by redefining the coefficients on the exceptional
$\diff t$-null set. If only \ref{B:ass:lip} is required, set $b=\sigma=h=f=0$ there. If
\ref{B:ass:mono} is also required, set
\[
(b,\sigma,h,f)
:=
\bigl(
-\beta_1\Gstar y,\,
-\beta_1\Gstar z,\,
-\beta_1\Gstar u,\,
\beta_2 Gx
\bigr),
\qquad
\theta=(x,y,z,u)\in\mathcal H.
\]
The latter choice satisfies the interior monotonicity inequality in \ref{B:ass:mono}
with equality and the Lipschitz inequality in \ref{B:ass:lip} after enlarging $L$ to
\[
\max\{L,\beta_1\|G\|,\beta_2\|G\|\}.
\]
Both choices leave the continuous-time equations unchanged, since the coefficients are modified
only on a $\diff t$-null set.

\begin{definition}[Admissible orthogonal martingales]\label{B:def:incr-adm}
A square-integrable $\F_{t_{i+1}}$-measurable increment $\Delta M_i$ is \emph{admissible} if
\[
\Ei[\Delta M_i]=0,\qquad \Ei[\Delta M_i\Delta W_i^{\!\top}]=0,
\]
and
\[
\Ei\Big[\Delta M_i\Big(\!\int_E a\,\Delta\Nt_i\Big)^{\!\top}\Big]=0
\]
for every square-integrable $\F_{t_i}$-measurable $L^2(\nu)$-integrand $a$ of compatible
dimension.

An adapted square-integrable process $M=(M_i)_{i=0}^N$ with $M_0=0$ is \emph{admissible} if
each increment $\Delta M_i$ is admissible.
\end{definition}

\paragraph{Discrete solution class.}
Write
\[
\begin{aligned}
\Sspi{\R^j}&:=\bigl\{(V_i)_{i\le N}\ \F_{t_i}\text{-adapted}:\ \max_{0\le i\le N}\E|V_i|^2<\infty\bigr\},\\[4pt]
\HWpi{\R^{m\times d}}&:=\bigl\{(Z_i)_{i<N}\ \F_{t_i}\text{-adapted}:\ \textstyle\sum_{i=0}^{N-1}\tau\,\E|Z_i|_F^2<\infty\bigr\},\\[4pt]
\HNpi{L^2(\nu;\R^m)}&:=\bigl\{(U_i)_{i<N}\ \F_{t_i}\text{-adapted}:\ \textstyle\sum_{i=0}^{N-1}\tau\,\E\|U_i\|_{\Lnu}^2<\infty\bigr\},\\[4pt]
\Mpi{\R^m}&:=\bigl\{M\in\Sspi{\R^m}:\ M\ \text{admissible}\bigr\}.
\end{aligned}
\]
Set
\[
\mathcal S_\pi:=\Sspi{\R^n}\times\Sspi{\R^m}\times\HWpi{\R^{m\times d}}
\times\HNpi{L^2(\nu;\R^m)}\times\Mpi{\R^m},
\]
with squared norm
\begin{equation}\label{B:eq:Spi-norm}
\|\Theta\|_{\mathcal S_\pi}^2:=\max_{0\le i\le N}\E|X_i|^2+\max_{0\le i\le N}\E|Y_i|^2
+\sum_{i=0}^{N-1}\tau\big(\E|Z_i|_F^2+\E\|U_i\|_{\Lnu}^2\big)+\E|M_N|^2.
\end{equation}
We call elements of $\mathcal S_\pi$ admissible grid approximations.

\begin{lemma}[One-step Kunita--Watanabe projection]\label{B:lem:kw}
Let
\[
R\in L^2(\Omega,\F_{t_{i+1}};\R^m),\qquad \Ei[R]=0.
\]
There exist unique
\[
Z\in L^2(\Omega,\F_{t_i};\R^{m\times d}),\qquad
U\in L^2(\Omega,\F_{t_i};L^2(\nu;\R^m)),
\]
and an admissible increment
\[
\Delta M\in L^2(\Omega,\F_{t_{i+1}};\R^m)
\]
such that
\begin{equation}\label{B:eq:kw-decomp}
R=Z\,\Delta W_i+\int_E U(e)\,\Delta\Nt_i(\diff e)+\Delta M.
\end{equation}
Moreover,
\[
Z=\tau^{-1}\Ei[R\,\Delta W_i^{\!\top}],
\]
while $U$ is uniquely characterised by
\[
\Ei\Big[\inner{R}{\int_E a\,\Delta\Nt_i}\Big]=\tau\int_E\inner{U(e)}{a(e)}\,\nu(\diff e)
\]
for every square-integrable $\F_{t_i}$-measurable test integrand $a$.
\end{lemma}
\noindent The proof is given in Appendix~\ref{B:app:kw}.

For $d\ge1$, a nonzero admissible increment is
\[
\Delta M_i=\bigl((\Delta W^{(1)}_i)^2-\tau\bigr)e_1,
\]
where $e_1$ is the first coordinate vector of $\R^m$.

\subsection{Discrete stability and well-posedness}\label{B:sec:stab}
Let
\[
\chi_i\in\R^n,\qquad \gamma_i\in\R^{n\times d},\qquad \eta_i\in L^2(\nu;\R^n),\qquad
\phi_i\in\R^m,\qquad \zeta\in L^2(\Omega,\F_T;\R^m),
\]
where $\chi_i,\gamma_i,\eta_i,\phi_i$ are square-integrable $\F_{t_i}$-measurable perturbations
and $\zeta$ is the terminal perturbation. A \emph{perturbed discrete solution}
$\Theta'\in\mathcal S_\pi$ is evaluated at
\[
\begin{aligned}
\mu'_i&:=\Law(X'_i,Y'_i,Z'_i,U'_i),\\[2pt]
\psi'_i&:=\psi(t_i,X'_i,Y'_i,Z'_i,U'_i,\mu'_i),\qquad \psi\in\{b,\sigma,f\},\\[2pt]
h'_i(e)&:=h(t_i,X'_i,Y'_i,Z'_i,U'_i,e,\mu'_i),
\end{aligned}
\]
and satisfies, at each grid step,
\begin{align}
\Delta X'_i&=b'_i\tau+\sigma'_i\Delta W_i+\int_E h'_i\,\Delta\Nt_i+\chi_i\tau+\gamma_i\Delta W_i+\int_E\eta_i\,\Delta\Nt_i,\label{B:eq:pfwd}\\
\Delta Y'_i&=-f'_i\tau+Z'_i\Delta W_i+\int_E U'_i\,\Delta\Nt_i+\Delta M'_i+\phi_i\tau,\qquad Y'_N=g(X'_N,\Law(X'_N))+\zeta.\label{B:eq:pbwd}
\end{align}
Write
\[
\delta A_i:=\ex A_i-A'_i,\qquad A\in\{X,Y,Z,U,M\},
\]
with the same convention for coefficient differences. For the perturbation source
$s=(\chi,\gamma,\eta,\phi,\zeta)$, the squared norm of its interior terms is
\begin{equation}\label{B:eq:Pnorm}
\Pp :=
\sum_{i=0}^{N-1}\tau\,\E\!\left[
|\chi_i|^2+|\gamma_i|_F^2+\|\eta_i\|_{\Lnu}^2+|\phi_i|^2
\right].
\end{equation}
Including the terminal perturbation $\zeta$ gives the squared source norm
\begin{equation}\label{B:eq:srcnorm}
\|s\|_{\mathrm{src}}^2 := \Pp+\E|\zeta|^2.
\end{equation}

\begin{proposition}[Coefficient-perturbation stability]\label{B:prop:pert}
Under Assumptions~\ref{B:ass:integ}--\ref{B:ass:mono} and
Convention~\ref{B:conv:nodewise}, let $\ex\Theta$ be the exact discrete solution and
$\Theta'$ the perturbed discrete solution of \eqref{B:eq:pfwd}--\eqref{B:eq:pbwd}. There are $\tau_{\mathrm{stab}}>0$ and
$C_{\mathrm{stab}}<\infty$, depending only on
$L,T,\|G\|,c_G,\alpha,\beta_1,\beta_2$, such that for all $\tau<\tau_{\mathrm{stab}}$,
\begin{equation}
\label{B:eq:pert}
\begin{split}
&\max_{0\le i\le N}\E|\delta X_i|^2+\max_{0\le i\le N}\E|\delta Y_i|^2+\sum_{i=0}^{N-1}\tau\,\E|\delta Z_i|_F^2\\
&\qquad+\sum_{i=0}^{N-1}\tau\,\E\|\delta U_i\|_{\Lnu}^2+\E|\delta M_N|^2
\le C_{\mathrm{stab}}\big(\E|\delta X_0|^2+\|s\|_{\mathrm{src}}^2\big).
\end{split}
\end{equation}
The same estimate holds for two perturbed solutions, with $s$ replaced by $s^1-s^2$.
\end{proposition}

\noindent The proof is given in Appendix~\ref{B:app:stab}.

\begin{theorem}[Discrete well-posedness]\label{B:thm:wp}
Under Assumptions~\ref{B:ass:integ}--\ref{B:ass:mono} and
Convention~\ref{B:conv:nodewise}, there is $\tau_0>0$ such that for every
$\tau<\tau_0$ the discrete Euler scheme \eqref{B:eq:dfwd}--\eqref{B:eq:dbwd} admits a unique solution
$\ex\Theta\in\mathcal S_\pi$.
\end{theorem}
\begin{proof}
Existence is proved by continuation in Appendix~\ref{B:app:wp}. Uniqueness follows from
Proposition~\ref{B:prop:pert} with zero perturbation data.
\end{proof}

\subsection{The \emph{a posteriori} estimator}
The estimator is defined for any $\ap\Theta\in\mathcal S_\pi$. For $0\le i\le N-1$, set
\[
\ap\mu_i:=\Law(\ap X_i,\ap Y_i,\ap Z_i,\ap U_i),
\]
and
\[
\begin{aligned}
\ap\psi_i
&:=\psi\bigl(t_i,\ap X_i,\ap Y_i,\ap Z_i,\ap U_i,\ap\mu_i\bigr),
&&\psi\in\{b,\sigma,f\},\\[2pt]
\ap h_i(e)
&:=h\bigl(t_i,\ap X_i,\ap Y_i,\ap Z_i,\ap U_i,e,\ap\mu_i\bigr).
\end{aligned}
\]

\begin{definition}[\emph{A posteriori} estimator]\label{B:def:Epi}
Define the boundary residuals
\[
r_0^X:=\ap X_0-\xi_0,\qquad r_T^Y:=\ap Y_N-g(\ap X_N,\Law(\ap X_N)),
\]
and, for $0\le i\le N-1$, the cumulative residuals
\begin{equation}
\label{B:eq:CX}
C_i^X:=\ap X_{i+1}-\ap X_0-\sum_{j=0}^{i}\Big(\ap b_j\tau+\ap\sigma_j\Delta W_j+\int_E\ap h_j\,\Delta\Nt_j\Big),
\end{equation}
\begin{equation}
\label{B:eq:CY}
C_i^Y:=\ap Y_{i+1}-\ap Y_0+\sum_{j=0}^{i}\Big(\ap f_j\tau-\ap Z_j\Delta W_j-\int_E\ap U_j\,\Delta\Nt_j\Big)-\ap M_{i+1}.
\end{equation}
The \emph{a posteriori estimator} is
\begin{equation}
\label{B:eq:Epi}
\Epi:=\E|r_0^X|^2+\E|r_T^Y|^2+\max_{i<N}\E|C_i^X|^2+\max_{i<N}\E|C_i^Y|^2.
\end{equation}
\end{definition}

The discrete error \eqref{B:eq:Dpi} is
\[
\Dpi=\|\ap\Theta-\ex\Theta\|_{\mathcal S_\pi}^2 .
\]

\section{Discrete \emph{a posteriori} estimates}
\label{B:sec:main}
\begin{theorem}[Discrete \emph{a posteriori} estimates]\label{B:thm:discrete}
Let $\ap\Theta\in\mathcal S_\pi$ be arbitrary. Use the coefficient representatives fixed in
Convention~\ref{B:conv:nodewise}.
\begin{enumerate}[label=\textup{(\roman*)},ref=\roman*,leftmargin=2.4em,itemsep=3pt]
\item\label{B:item:eff} \emph{Efficiency.} Let Assumption~\ref{B:ass:lip} hold. For any solution
$\ex\Theta$ of \eqref{B:eq:dfwd}--\eqref{B:eq:dbwd},
\begin{equation}
\label{B:eq:eff}
\Epi\;\le\;C_{\mathrm{eff}}\,\Dpi
\end{equation}
on \emph{any} grid. Here
\[
C_{\mathrm{eff}}\lesssim(1+T)^2(1+L)^2.
\]

\item\label{B:item:rel} \emph{Reliability.} Let Assumptions~\ref{B:ass:integ}--\ref{B:ass:mono} hold. There are $\tau_0>0$ and
$C<\infty$, depending only on $L,T,\|G\|,c_G,\alpha,\beta_1,\beta_2$, such that, for every
$0<\tau<\tau_0$,
\begin{equation}
\label{B:eq:rel}
\Dpi\;\le\;C\,\Epi.
\end{equation}
\end{enumerate}

Consequently, for every $0<\tau<\tau_0$,
\begin{equation}
\label{B:eq:discrete}
c\,\Epi\;\le\;\Dpi\;\le\;C\,\Epi,\qquad c:=C_{\mathrm{eff}}^{-1}.
\end{equation}
\end{theorem}
\begin{proof}
Write $\delta\Theta:=\ex\Theta-\ap\Theta$, with the same convention for coefficient
differences. Throughout, $C$ denotes a generic constant independent of the grid.

\medskip\noindent\emph{Efficiency.}
Since $\ex\Theta$ satisfies \eqref{B:eq:dfwd}--\eqref{B:eq:dbwd}, the residual definitions
\eqref{B:eq:CX}--\eqref{B:eq:CY} give, for $0\le i\le N-1$,
\begin{equation}
\label{B:eq:Ceff}
\begin{aligned}
-C_i^X&=\delta X_{i+1}-\delta X_0-\sum_{j=0}^{i}\Big(\delta b_j\tau+\delta\sigma_j\Delta W_j+\int_E\delta h_j\,\Delta\Nt_j\Big),\\[2pt]
-C_i^Y&=\delta Y_{i+1}-\delta Y_0+\sum_{j=0}^{i}\Big(\delta f_j\tau-\delta Z_j\Delta W_j-\int_E\delta U_j\,\Delta\Nt_j\Big)-\delta M_{i+1}.
\end{aligned}
\end{equation}
Since $\ex M$ and $\ap M$ are admissible, so is $\delta M$, and
$\delta M_{i+1}=\sum_{j=0}^{i}\Delta\delta M_j$. Martingale differences at distinct steps are
orthogonal, and the conditional orthogonality of $\Delta W_j$ and $\Delta\Nt_j$, together with
admissibility, eliminates the same-step cross terms. Hence
\[
\begin{aligned}
\E\Big|\sum_{j=0}^{i}\Big(\delta\sigma_j\Delta W_j+\int_E\delta h_j\,\Delta\Nt_j\Big)\Big|^2
&=\sum_{j=0}^{i}\tau\,\E\big(|\delta\sigma_j|_F^2+\|\delta h_j\|_{\Lnu}^2\big),\\
\E\Big|\sum_{j=0}^{i}\Big(\delta Z_j\Delta W_j+\int_E\delta U_j\,\Delta\Nt_j\Big)+\delta M_{i+1}\Big|^2
&=\sum_{j=0}^{i}\tau\,\E\big(|\delta Z_j|_F^2+\|\delta U_j\|_{\Lnu}^2\big)+\E|\delta M_{i+1}|^2,
\end{aligned}
\]
with $\E|\delta M_{i+1}|^2\le\E|\delta M_N|^2$, while
\[
\begin{aligned}
\E\Big|\sum_{j=0}^{i}\delta o_j\tau\Big|^2
&\le\Big(\sum_{j=0}^{i}\tau\Big)\Big(\sum_{j=0}^{i}\tau\,\E|\delta o_j|^2\Big)\\
&\le T\sum_{j=0}^{i}\tau\,\E|\delta o_j|^2,\qquad o\in\{b,f\}.
\end{aligned}
\]
By Assumption~\ref{B:ass:lip} and Lemma~\ref{B:lem:law},
\[
\begin{aligned}
&\sum_{j=0}^{N-1}\tau\,\E\big(|\delta b_j|^2+|\delta\sigma_j|_F^2+\|\delta h_j\|_{\Lnu}^2+|\delta f_j|^2\big)\\
&\qquad\le C\sum_{j=0}^{N-1}\tau\,\E\big(|\delta X_j|^2+|\delta Y_j|^2+|\delta Z_j|_F^2+\|\delta U_j\|_{\Lnu}^2
+W_2(\ap\mu_j,\ex\mu_j)^2\big)\\
&\qquad\le C\,\Dpi,
\end{aligned}
\]
so
\[
\max_{i<N}\E|C_i^X|^2+\max_{i<N}\E|C_i^Y|^2\le C\,\Dpi.
\]
For the boundary residuals,
\[
\E|r_0^X|^2=\E|\delta X_0|^2\le\Dpi.
\]
Assumption~\ref{B:ass:lip} and Lemma~\ref{B:lem:law}, applied to $g$, with
\[
\delta g_N:=g(\ex X_N,\Law(\ex X_N))-g(\ap X_N,\Law(\ap X_N)),
\]
give
\[
\begin{aligned}
\E|r_T^Y|^2&=\E|\delta Y_N-\delta g_N|^2\\
&\le2\,\E|\delta Y_N|^2+2\,\E|\delta g_N|^2\\
&\le C\big(\E|\delta Y_N|^2+\E|\delta X_N|^2\big)\\
&\le C\,\Dpi.
\end{aligned}
\]
Hence
\[
\Epi\le C_{\mathrm{eff}}\,\Dpi,\qquad C_{\mathrm{eff}}\lesssim(1+T)^2(1+L)^2.
\]

\medskip\noindent\emph{Reliability.}
Define the auxiliary tuple $\au\Theta=(\au X,\au Y,\au Z,\au U,\au M)$ by
\[
\au Z\equiv\ap Z,\qquad \au U\equiv\ap U,\qquad \au M\equiv\ap M,\qquad
\au X_0=\ap X_0,\qquad \au Y_0=\ap Y_0,
\]
and let its $X$- and $Y$-components evolve according to
\begin{equation}
\label{B:eq:aux}
\Delta\au X_i:=\ap b_i\tau+\ap\sigma_i\Delta W_i+\int_E\ap h_i\,\Delta\Nt_i,\qquad
\Delta\au Y_i:=-\ap f_i\tau+\ap Z_i\Delta W_i+\int_E\ap U_i\,\Delta\Nt_i+\Delta\ap M_i .
\end{equation}
Set
\[
\begin{aligned}
\au\mu_i
&:=\Law(\au X_i,\au Y_i,\au Z_i,\au U_i),\\
\au\psi_i
&:=\psi(t_i,\au X_i,\au Y_i,\au Z_i,\au U_i,\au\mu_i),
\qquad \psi\in\{b,\sigma,h,f\}.
\end{aligned}
\]

Summing \eqref{B:eq:aux} gives, for $0\le i\le N-1$,
\begin{equation}
\label{B:eq:gap}
\ap X_{i+1}-\au X_{i+1}=C_i^X,\qquad \ap Y_{i+1}-\au Y_{i+1}=C_i^Y.
\end{equation}
\emph{Perturbed scheme.} The auxiliary tuple $\au\Theta$ satisfies
\eqref{B:eq:pfwd}--\eqref{B:eq:pbwd} with perturbations
\[
\begin{aligned}
\chi_i&=\ap b_i-\au b_i, &\qquad \gamma_i&=\ap\sigma_i-\au\sigma_i,\\[2pt]
\eta_i(e)&=\ap h_i(e)-\au h_i(e), &\qquad \phi_i&=\au f_i-\ap f_i,
\end{aligned}
\]
and terminal perturbation $\zeta=\au Y_N-g(\au X_N,\Law(\au X_N))$.

By Lipschitz continuity, law coupling, and $\au Z\equiv\ap Z$,
$\au U\equiv\ap U$,
\begin{equation}
\label{B:eq:pertbd}
\E\big[|\chi_i|^2+|\gamma_i|_F^2+\|\eta_i\|_{\Lnu}^2+|\phi_i|^2\big]
\le C\,\E\big[|\ap X_i-\au X_i|^2+|\ap Y_i-\au Y_i|^2\big].
\end{equation}
By \eqref{B:eq:Pnorm}, \eqref{B:eq:gap}, and \eqref{B:eq:pertbd},
\[
\begin{aligned}
\Pp
&\le C\sum_{i=0}^{N-1}\tau\,\E\big[|\ap X_i-\au X_i|^2+|\ap Y_i-\au Y_i|^2\big]
=C\sum_{i=1}^{N-1}\tau\,\E\big[|C^X_{i-1}|^2+|C^Y_{i-1}|^2\big]\\
&\le CT\big(\max_{i<N}\E|C^X_i|^2+\max_{i<N}\E|C^Y_i|^2\big)\\
&\le C\,\Epi .
\end{aligned}
\]

Moreover,
\[
\E|\ex X_0-\au X_0|^2=\E|r_0^X|^2\le\Epi
\]
and
\[
\zeta=(\au Y_N-\ap Y_N)+r_T^Y+\big[g(\ap X_N,\Law(\ap X_N))-g(\au X_N,\Law(\au X_N))\big],
\]
so, by Lipschitz $g$, Lemma~\ref{B:lem:law} and \eqref{B:eq:gap},
\[
\begin{aligned}
\E|\zeta|^2
&\le C\big(\E|\au Y_N-\ap Y_N|^2+\E|r_T^Y|^2+\E|\ap X_N-\au X_N|^2\big)\\
&=C\big(\E|C^Y_{N-1}|^2+\E|r_T^Y|^2+\E|C^X_{N-1}|^2\big)\\
&\le C\,\Epi .
\end{aligned}
\]

Thus by Proposition~\ref{B:prop:pert}, for $\tau<\tau_0$,
\begin{equation}
\label{B:eq:auxstabdisc}
\|\au\Theta-\ex\Theta\|_{\mathcal S_\pi}^2\le C\,\Epi .
\end{equation}
By \eqref{B:eq:gap} and the definition of $\au\Theta$,
\[
\|\ap\Theta-\au\Theta\|_{\mathcal S_\pi}^2
=\max_{i<N}\E|C_i^X|^2+\max_{i<N}\E|C_i^Y|^2\le\Epi.
\]
Hence, by \eqref{B:eq:auxstabdisc},
\[
\Dpi=\|\ap\Theta-\ex\Theta\|_{\mathcal S_\pi}^2
\le2\|\ap\Theta-\au\Theta\|_{\mathcal S_\pi}^2+2\|\au\Theta-\ex\Theta\|_{\mathcal S_\pi}^2
\le C\,\Epi.
\]
\end{proof}

\section{Continuous-time estimate}
\label{B:sec:contres}
\noindent Under Assumptions~\ref{B:ass:integ}--\ref{B:ass:mono}, the continuous
system~\eqref{B:eq:mvfbsdej} admits a unique solution
$\Theta\in\mathcal S$~\cite[Theorems~3.3 and~3.5]{CWPpaper}.

\begin{definition}[Continuous-time error]\label{B:def:Dc}
For $\ap\Theta\in\mathcal S_\pi$, define
\begin{equation}
\label{B:eq:Dc}
\begin{aligned}
\Dc:={}&\max_{0\le i\le N-1}\sup_{t\in[t_i,t_{i+1}]}\big(\E|\ap X_i-X_t|^2+\E|\ap Y_i-Y_t|^2\big)\\
&+\E|\ap X_N-X_T|^2+\E|\ap Y_N-Y_T|^2\\
&+\sum_{i=0}^{N-1}\E\!\int_{t_i}^{t_{i+1}}\!\big(|\ap Z_i-Z_t|_F^2+\|\ap U_i-U_t\|_{\Lnu}^2\big)\diff t+\E|\ap M_N|^2.
\end{aligned}
\end{equation}
The term $\E|\ap M_N|^2$ accounts for the discrete orthogonal-martingale component, which
has no continuous-time counterpart.
\end{definition}

\begin{definition}[Grid sampling]\label{B:def:grid-sampling}
Define the grid-sampled tuple $\sm\Theta:=(\sm X,\sm Y,\sm Z,\sm U,0)$ by
\[
\sm X_i:=X_{t_i},\qquad \sm Y_i:=Y_{t_i},
\]
\[
\sm Z_i:=\frac1\tau\Ei\Big[\int_{t_i}^{t_{i+1}}Z_s\diff s\Big],\qquad
\sm U_i:=\frac1\tau\Ei\Big[\int_{t_i}^{t_{i+1}}U_s\diff s\Big].
\]
Then $\sm Z_i$ and $\sm U_i$ are the $L^2$-projections of $Z$ and $U$ onto the
$\F_{t_i}$-measurable constants on $[t_i,t_{i+1}]$.
\end{definition}

\begin{definition}[Path-regularity modulus]\label{B:def:Reg}
The path-regularity modulus of $\Theta$ on $\pi$ is
\begin{equation}
\label{B:eq:Reg}
\begin{split}
\Reg(\tau)&:=\max_{0\le i\le N-1}\sup_{t\in[t_i,t_{i+1}]}\big(\E|X_t-\sm X_i|^2+\E|Y_t-\sm Y_i|^2\big)\\
&\quad+\sum_{i=0}^{N-1}\E\!\int_{t_i}^{t_{i+1}}\!\big(|Z_s-\sm Z_i|_F^2+\|U_s-\sm U_i\|_{\Lnu}^2\big)\diff s.
\end{split}
\end{equation}
\end{definition}

\begin{lemma}[The path-regularity modulus vanishes]\label{B:lem:regvanish}
Let Assumptions~\ref{B:ass:integ}--\ref{B:ass:lip} hold and let
$\Theta=(X,Y,Z,U)\in\mathcal S$ be any solution of \eqref{B:eq:mvfbsdej}. Then
\[
\Reg(\tau)\longrightarrow0\qquad\text{as }\tau\to0 .
\]
\end{lemma}

The proof is given in Appendix~\ref{B:app:regvanish}.

\begin{assumption}[Coefficient time-regularity]\label{B:ass:reg}
The common growth-weighted time-modulus of the coefficients
$b,\sigma,h,f$,
\begin{equation}
\label{B:eq:omega}
\omega(\tau):=\max_{\psi\in\{b,\sigma,h,f\}}\;\sup_{|t-s|\le\tau}\;
\sup_{\substack{\theta\in\mathcal H\\ \mu\in\mathcal P_2(\mathcal H)}}
\frac{\|\psi(t,\theta,\mu)-\psi(s,\theta,\mu)\|}{1+\|\theta\|_{\mathcal H}+W_2(\mu,\delta_0)},
\end{equation}
vanishes as $\tau\to0$.
\end{assumption}

By continuity in time, the coefficients in \eqref{B:eq:omega} satisfy the Lipschitz and
monotonicity inequalities at every grid node, as required by Convention~\ref{B:conv:nodewise}.

\begin{proposition}[Discretisation error]\label{B:prop:discerr}
Under Assumptions~\ref{B:ass:integ}--\ref{B:ass:reg}, there are
$\tau_0>0$ and $C<\infty$ such that, for every $0<\tau<\tau_0$,
\[
\begin{gathered}
\max_{0\le i\le N}\E|\sm X_i-\ex X_i|^2+\max_{0\le i\le N}\E|\sm Y_i-\ex Y_i|^2
+\sum_{i=0}^{N-1}\tau\big(\E|\sm Z_i-\ex Z_i|_F^2+\E\|\sm U_i-\ex U_i\|_{\Lnu}^2\big)
+\E|\ex M_N|^2\\[3pt]
\hphantom{{\le}\;C\big(\omega(\tau)^2+\Reg(\tau)\big).}%
{\le}\;C\big(\omega(\tau)^2+\Reg(\tau)\big).
\end{gathered}
\]
\end{proposition}
\begin{proof}
By \eqref{B:eq:Spi-norm}, Definition~\ref{B:def:grid-sampling}, and conditional Jensen and Cauchy--Schwarz,
\[
\|\sm\Theta\|_{\mathcal S_\pi}^2\le\E\sup_{t\le T}|X_t|^2+\E\sup_{t\le T}|Y_t|^2
+\E\int_0^T\big(|Z_s|_F^2+\|U_s\|_{\Lnu}^2\big)\diff s<\infty.
\]
Hence $\sm\Theta\in\mathcal S_\pi$. For $0\le i<N$, set
\[
\sm\Theta_i:=(\sm X_i,\sm Y_i,\sm Z_i,\sm U_i),
\qquad
\sm\mu_i:=\Law(\sm\Theta_i),
\]
and write
\[
\begin{aligned}
\sm b_i &:= b(t_i,\sm\Theta_i,\sm\mu_i),
&\qquad \sm\sigma_i &:= \sigma(t_i,\sm\Theta_i,\sm\mu_i),\\
\sm f_i &:= f(t_i,\sm\Theta_i,\sm\mu_i),
&\qquad \sm h_i(e) &:= h(t_i,\sm\Theta_i,e,\sm\mu_i).
\end{aligned}
\]
Define the auxiliary tuple $\au\Theta:=(\au X,\au Y,\sm Z,\sm U,0)$ by
$\au X_0=\sm X_0$, $\au Y_0=\sm Y_0$, and
\begin{equation}\label{B:eq:auxcont}
\begin{aligned}
\Delta\au X_i
&:=\sm b_i\tau+\sm\sigma_i\Delta W_i+\int_E\sm h_i\,\Delta\Nt_i,\\
\Delta\au Y_i
&:=-\sm f_i\tau+\sm Z_i\Delta W_i+\int_E\sm U_i\,\Delta\Nt_i.
\end{aligned}
\end{equation}
The boundary conditions in \eqref{B:eq:mvfbsdej} give $r_0^X=0$ and $r_T^Y=0$.
Applying \eqref{B:eq:gap} with $\ap\Theta=\sm\Theta$ gives, for $0\le i<N$,
\[
C_i^X=\sm X_{i+1}-\au X_{i+1},
\qquad
C_i^Y=\sm Y_{i+1}-\au Y_{i+1}.
\]
Thus Lemma~\ref{B:lem:freezing-error} gives
\[
\begin{gathered}
\Epi(\sm\Theta)
=\max_{0\le i\le N-1}\E|\sm X_{i+1}-\au X_{i+1}|^2
+\max_{0\le i\le N-1}\E|\sm Y_{i+1}-\au Y_{i+1}|^2\\[3pt]
\hphantom{{\le}\;C\bigl(\omega(\tau)^2+\Reg(\tau)\bigr).}%
{\le}\;C\bigl(\omega(\tau)^2+\Reg(\tau)\bigr).
\end{gathered}
\]
For $0<\tau<\tau_0$, Theorem~\ref{B:thm:discrete}(\ref{B:item:rel}) therefore gives
\[
\|\sm\Theta-\ex\Theta\|_{\mathcal S_\pi}^2\le C\,\Epi(\sm\Theta)
\le C\bigl(\omega(\tau)^2+\Reg(\tau)\bigr),
\]
which, by \eqref{B:eq:Spi-norm}, is the asserted estimate.
\end{proof}
In particular, by Assumption~\ref{B:ass:reg} and Lemma~\ref{B:lem:regvanish},
\[
\E|\ex M_N|^2\le C\bigl(\omega(\tau)^2+\Reg(\tau)\bigr)\longrightarrow0\qquad\text{as }\tau\to0.
\]

\begin{theorem}[Continuous-time two-sided \emph{a posteriori} estimate]\label{B:thm:cont}
Under Assumptions~\ref{B:ass:integ}--\ref{B:ass:reg}, there are
$0<c\le C<\infty$ and $\tau_0>0$ such that, for every $0<\tau<\tau_0$,
\begin{equation}
\label{B:eq:cont}
\big(c\,\Epi-C(\omega(\tau)^2+\Reg(\tau))\big)_+\;\le\;\Dc\;\le\;C\big(\Epi+\omega(\tau)^2+\Reg(\tau)\big).
\end{equation}
If, in addition,
\[
\omega(\tau)=O(\tau^{1/2}),\qquad \Reg(\tau)=O(\tau)\qquad\text{as }\tau\downarrow0,
\]
then, after increasing $C$ and decreasing $\tau_0$ if necessary,
\begin{equation}
\label{B:eq:cont-rate}
(c\,\Epi-C\tau)_+\;\le\;\Dc\;\le\;C(\Epi+\tau),\qquad 0<\tau<\tau_0.
\end{equation}
\end{theorem}
\begin{proof}
Set
\[
\sm\Dpi:=\|\ap\Theta-\sm\Theta\|_{\mathcal S_\pi}^2,\qquad
R_\tau:=\omega(\tau)^2+\Reg(\tau).
\]
For $t\in[t_i,t_{i+1}]$,
\[
\begin{aligned}
&\E|\ap X_i-X_t|^2+\E|\ap Y_i-Y_t|^2\\
&\qquad\le
2\E|\ap X_i-\sm X_i|^2
+2\E|\ap Y_i-\sm Y_i|^2
+2\E|X_t-\sm X_i|^2
+2\E|Y_t-\sm Y_i|^2.
\end{aligned}
\]
Moreover, by the cellwise projection property,
\[
\int_{t_i}^{t_{i+1}}\E|\ap Z_i-Z_t|_F^2\,\diff t
=
\tau\E|\ap Z_i-\sm Z_i|_F^2
+
\int_{t_i}^{t_{i+1}}\E|Z_t-\sm Z_i|_F^2\,\diff t,
\]
with the analogous identity for $U$. Since $\sm X_N=X_T$ and $\sm Y_N=Y_T$, these estimates give
\begin{equation}\label{B:eq:projDc}
\sm\Dpi\le2\Dc,\qquad
\Dc\le3\,\sm\Dpi+2\Reg(\tau).
\end{equation}
By Proposition~\ref{B:prop:discerr} and the squared triangle inequality in $\mathcal S_\pi$,
\begin{equation}\label{B:eq:projDpi}
\sm\Dpi\le2\Dpi+C R_\tau,\qquad
\Dpi\le2\,\sm\Dpi+C R_\tau.
\end{equation}
By Theorem~\ref{B:thm:discrete} and \eqref{B:eq:projDc}--\eqref{B:eq:projDpi},
\[
\begin{aligned}
\Dc&\le C\bigl(\sm\Dpi+\Reg(\tau)\bigr)\le C\bigl(\Dpi+R_\tau\bigr)\le C\bigl(\Epi+R_\tau\bigr),\\[4pt]
\Epi&\le C\,\Dpi\le C\bigl(\sm\Dpi+R_\tau\bigr)\le C\bigl(\Dc+R_\tau\bigr).
\end{aligned}
\]
Thus, for suitable $0<c\le C$,
\[
\bigl(c\,\Epi-C R_\tau\bigr)_+\le\Dc\le C\bigl(\Epi+R_\tau\bigr),
\]
which is \eqref{B:eq:cont}. Under the additional rate assumptions, $R_\tau=O(\tau)$, giving
\eqref{B:eq:cont-rate}.
\end{proof}

\section{Numerical framework}\label{B:sec:validation}

\subsection{Closed-form benchmarks}\label{B:sec:numerics}

\subsubsection{Linear benchmark and protocol}\label{B:subsec:benchmark-protocol}
The benchmark is the monotone linear MV-FBSDEJ with coefficients
\[
\begin{gathered}
b=\bar a\,G^+G\,\E X-\kappa_{\mathrm F}\Gstar Y,\qquad \sigma=\sigma_0-\kappa_{\mathrm F}\Gstar Z,\qquad h(e)=h_0-\kappa_{\mathrm F}\Gstar U(e),\\[4pt]
f=\kappa_{\mathrm B} GX+\bar c\,G\,\E X,\qquad g=GX+\bar g\,G\,\E X,
\end{gathered}
\]
with a single jump mark of intensity mass $\nu(E)$.

Here $G^+$ is the Moore--Penrose pseudoinverse, so $G^+G$ is the orthogonal projector onto
$(\ker G)^\perp$. The mean-field terms depend on the law only through $G\,\E X=\E[GX]$.

The base scalar benchmark takes
\begin{equation}\label{B:eq:base-scalar}
\begin{gathered}
n=m=d=1,\qquad G=\rv{apocert.G},\qquad X_0\sim\rv{apocert.xZeroLaw},\qquad T=\rv{apocert.T},\\[4pt]
\kappa_{\mathrm F}=\kappa_{\mathrm B}=\rv{apocert.betaOne},\qquad \bar a=\bar c=\bar g=\rv{apocert.coupling},\qquad
\sigma_0=h_0=\rv{apocert.sigmaZero},\qquad \nu(E)=\rv{apocert.nuE}.
\end{gathered}
\end{equation}
The vector rank cases have the same $T$ and $\nu(E)$, and take
\[
X_0\sim\rv{vector.xZeroLaw},\qquad d=1,\qquad \bar a=\bar c=\bar g=\rv{vector.coupling},\qquad
\kappa_{\mathrm F}=\kappa_{\mathrm B}=\rv{vector.betaOne},
\]
with every entry of $\sigma_0\in\R^{n\times1}$ and $h_0\in\R^n$ equal to $\rv{vector.noiseEntry}$.
The coupling matrices and their smallest singular values are
\[
\begin{array}{cccc}
\toprule
(n,m) & G & s_{\min}(G) & c_G=s_{\min}(G)^{-2}\\
\midrule
(2,1) & (\rv{vector.G.mltn.11},\;\rv{vector.G.mltn.12}) & \sqrt{\rv{vector.sminSq.mltn}} & \rv{vector.cG.mltn}\\[4pt]
(2,2) & \begin{pmatrix}\rv{vector.G.meqn.11}&\rv{vector.G.meqn.12}\\ \rv{vector.G.meqn.21}&\rv{vector.G.meqn.22}\end{pmatrix}
      & \rv{vector.smin.meqn} & \rv{vector.cG.meqn}\\[10pt]
(1,2) & (\rv{vector.G.mgtn.11},\;\rv{vector.G.mgtn.21})^{\!\top} & \sqrt{\rv{vector.sminSq.mgtn}} & \rv{vector.cG.mgtn}\\
\bottomrule
\end{array}
\]

The scalar experiments also test $\kappa_{\mathrm F}=\kappa_{\mathrm B}\in\{\rv{apocert.betaWeak},\rv{apocert.betaStrong}\}$.
The vector cases cover $m<n$, $m=n$ and $m>n$.
Remark~\ref{B:rem:direct-benchmark-mono} verifies Assumption~\ref{B:ass:mono} for these choices, with
positive interior coercivity constants and $\alpha=1$.

\paragraph{Exact discrete reference.} The exact discrete solution has the affine form
\[
\ex Y_i=P_i\ex X_i+p_i,\qquad P_N=G.
\]
By centring and uniqueness,
\[
\E\ex X_i=\E\ex Y_i=0,
\]
and hence $p_i=0$.
Taking conditional expectations in the forward and backward steps gives
\[
\begin{aligned}
\Ei[\ex X_{i+1}]&=\ex X_i-\kappa_{\mathrm F}\tau\Gstar\ex Y_i,\\
\ex Y_i&=\Ei[\ex Y_{i+1}]+\kappa_{\mathrm B}\tau G\ex X_i.
\end{aligned}
\]
Substituting the centred affine relation yields
\[
\ex Y_i=P_{i+1}\bigl(\ex X_i-\kappa_{\mathrm F}\tau\Gstar\ex Y_i\bigr)+\kappa_{\mathrm B}\tau G\ex X_i.
\]
Using $\ex Y_i=P_i\ex X_i$ gives the Riccati recursion
\begin{equation}\label{B:eq:riccati}
P_i=(I+\kappa_{\mathrm F}\tau P_{i+1}\Gstar)^{-1}(P_{i+1}+\kappa_{\mathrm B}\tau G).
\end{equation}
The forward and backward steps give two representations of the centred increment:
\[
\begin{aligned}
\ex Y_{i+1}-\Ei[\ex Y_{i+1}]
&=P_{i+1}\sigma_i\Delta W_i+\int_E P_{i+1}h_i\,\Delta\Nt_i\\
&=\ex Z_i\Delta W_i+\int_E\ex U_i\,\Delta\Nt_i+\Delta\ex M_i.
\end{aligned}
\]
By Lemma~\ref{B:lem:kw},
\[
\ex Z_i=P_{i+1}\sigma_i,\qquad \ex U_i=P_{i+1}h_i,\qquad \Delta\ex M_i=0.
\]
For this benchmark,
\[
\sigma_i=\sigma_0-\kappa_{\mathrm F}\Gstar\ex Z_i,\qquad h_i=h_0-\kappa_{\mathrm F}\Gstar\ex U_i.
\]
Solving for the integrands gives
\[
\begin{aligned}
\ex Z_i&=(I+\kappa_{\mathrm F}P_{i+1}\Gstar)^{-1}P_{i+1}\sigma_0,\\
\ex U_i&=(I+\kappa_{\mathrm F}P_{i+1}\Gstar)^{-1}P_{i+1}h_0.
\end{aligned}
\]
Both matrices inverted above are invertible along the positive Riccati branch
(Appendix~\ref{B:app:wp}).
Since every increment vanishes, $\ex M\equiv0$. Perturbing only $M$ gives
\[
\Dpi=\Epi=\E|\ap M_N|^2,
\]
so the effectivity is one.
For nonzero perturbations confined to $Z$ and $U$ in the scalar benchmark,
Appendix~\ref{B:app:deferred} gives
\[
\frac{\Dpi}{\Epi}=\frac{1}{1+\kappa_{\mathrm F}^{\,2}}.
\]

\paragraph{Perturbation protocol.} To test the estimator independently of a particular numerical
solver, we perturb the exact Euler tuple $\ex\Theta$. For amplitude $\varepsilon>0$, set
\[
\ap\Theta=\ex\Theta+\delta\Theta.
\]
In the scalar benchmark, we perturb each component separately and all five together, using the
following perturbations.
\[
\begin{aligned}
\delta X_i&=\varepsilon\bigl(\sin\ex X_i+t_i/T\bigr),
&
\delta Y_i&=\varepsilon\bigl(\cos\ex Y_i+t_i/T\bigr),\\
\delta Z_i&=\varepsilon\bigl(\cos\ex X_i+\tfrac12\bigr),
&
\delta U_i&=\varepsilon\bigl(\sin\ex X_i+\tfrac12\bigr),\\
\Delta\delta M_i
&=
\varepsilon\frac{(\Delta W_i)^2-\tau}{\sqrt{2\tau}}.
\end{aligned}
\]
For the vector sweeps, where $d=1$, the corresponding perturbations are
\[
\begin{aligned}
\delta X_i
&=\varepsilon\bigl(\sin\ex X_i+(t_i/T)\mathbf 1_n\bigr),
&
\delta Y_i
&=\varepsilon\bigl(\cos\ex Y_i+(t_i/T)\mathbf 1_m\bigr),
\\[3pt]
\delta Z_i
&=\varepsilon\bigl(\cos(\ex X_i)_1+\tfrac12\bigr)\mathbf 1_m,
&
\delta U_i
&=\varepsilon\bigl(\sin(\ex X_i)_1+\tfrac12\bigr)\mathbf 1_m,
\\[3pt]
\Delta\delta M_i
&=
\varepsilon\frac{(\Delta W_i)^2-\tau}{\sqrt{2\tau}}\,\mathbf 1_m.
\end{aligned}
\]
Here $\mathbf 1_k$ denotes the vector of ones in $\R^k$.
Using $\ex X_i$ as a control variate for $\ap X_i$, we evaluate the mean-field coefficients at
\[
\widehat m_i
=\E\ex X_i
+\frac{1}{n_{\mathrm{MC}}}
\sum_{\ell=1}^{n_{\mathrm{MC}}}\delta X_i^{(\ell)},
\]
using the same $n_{\mathrm{MC}}$ paths as for $\Epi$ and $\Dpi$. We use this estimator in the
vector and mark-space variants. In the full-tuple-law study, we apply it to the corresponding
$X,Y,Z,U$ moments. In the rectangular cases, we test the two projections in separate runs:
\[
\begin{aligned}
\delta X_i &\quad\text{onto }\operatorname{ran}(\Gstar)\text{ and onto }\ker(G)
&&\text{when }m<n,\\
\delta Y_i &\quad\text{onto }\operatorname{ran}(G)\text{ and onto }\ker(\Gstar)
&&\text{when }m>n.
\end{aligned}
\]
We rescale each perturbation on an independent calibration batch to give a target
squared error of $\varepsilon^2$, then use the same scale on fresh scoring paths. Study-specific mark-space
directions and normalisations are stated where they are used. Table~\ref{B:tab:studies} gives the time
steps, paths per run, seeds and perturbation scales of every closed-form study; the complete per-study
settings are stated in the corresponding experiment scripts in the repository.

\subsubsection{Full-tuple-law variant}\label{B:ssec:fulllaw-setup}
In the scalar case $n=m=d=1$, $G=1$,
\[
\begin{gathered}
\mathfrak m_t=\lambda_X\E X_t+\lambda_Y\E Y_t+\lambda_Z\E Z_t+\lambda_U\E\!\int_E U_t\,\diff\nu,
\qquad(\lambda_X,\lambda_Y,\lambda_Z,\lambda_U)=\flLamVal,\\[3pt]
b=\mathfrak m_t-\kappa_{\mathrm F}Y,\qquad \sigma=\sigma_0-\kappa_{\mathrm F}Z,\qquad h(e)=h_0-\kappa_{\mathrm F}U(e),\qquad
f=\kappa_{\mathrm B}X,\qquad g(x,\rho)=x,\\[3pt]
X_0\sim\rv{fulllaw.xZeroLaw},\qquad T=\rv{fulllaw.T},\qquad \kappa_{\mathrm F}=\kappa_{\mathrm B}=\rv{fulllaw.betaOne},\qquad
\sigma_0=h_0=\rv{fulllaw.sigmaZero},\qquad \nu(E)=\rv{fulllaw.nuE}.
\end{gathered}
\]
The exact discrete solution has the affine form
\[
\ex Y_i=A_i\ex X_i+C_i,
\]
where $A_i$ satisfies \eqref{B:eq:riccati}. Writing $\bar x_i:=\E\ex X_i$, the
mean-field statistic, state mean and affine offset satisfy
\[
\begin{aligned}
\mathfrak m_i&=\lambda_X\bar x_i+\lambda_Y(A_i\bar x_i+C_i)+\lambda_Z\ex Z_i+\lambda_U\nu(E)\ex U_i,\\
\bar x_{i+1}&=\bar x_i+\tau\bigl(\mathfrak m_i-\kappa_{\mathrm F}(A_i\bar x_i+C_i)\bigr),\\
(1+\kappa_{\mathrm F}\tau A_{i+1})C_i&=C_{i+1}+\tau A_{i+1}\mathfrak m_i,
\end{aligned}
\]
with $\bar x_0=\E X_0$ and $C_N=0$.
The deterministic system $(\bar x_i,C_i)$ is solved by Picard iteration.
Remark~\ref{B:rem:fulllaw-dissipation} verifies Assumption~\ref{B:ass:mono} for this benchmark, with
$\lambda_{\min}\ge\monoFullLawVal$.

\subsubsection{Sine-terminal benchmark}
\label{B:ssec:sine}
Take
\[
X_{i+1}=X_i+\sigma_0\Delta W_i+h_0\Delta\Nt_i,\qquad
f(x)=\kappa_{\mathrm B}x,\qquad
g(x)=ax+c_{\sin}\sin(kx),
\]
with
\[
a=\rv{sine.a},\qquad c_{\sin}=\rv{sine.csin},\qquad k=\rv{sine.k},\qquad \kappa_{\mathrm F}=0,\qquad
\kappa_{\mathrm B}=\sineBetaTwoVal,
\]
and
\[
\sigma_0=\rv{sine.sigmaZero},\qquad h_0=\rv{sine.hZero},\qquad \nu(E)=\rv{sine.nuE},\qquad
T=\rv{sine.T},\qquad X_0\sim\rv{sine.xZeroLaw}.
\]
Since
\[
g'(x)=a+c_{\sin}k\cos(kx)\ge a-c_{\sin}k=\sineAlphaVal,
\]
Assumption~\ref{B:ass:mono} holds with $\alpha=\sineAlphaVal$, $\beta_1=0$ and $\beta_2=\sineBetaTwoVal$.

The nonlinear terminal datum yields $\ex M\ne0$. Set
\[
\varphi:=\E\bigl[e^{\mathrm i k(\sigma_0\Delta W_i+h_0\Delta\Nt_i)}\bigr],\qquad
\Psi_i:=\varphi^{\,N-i},\qquad
A_i:=a+\kappa_{\mathrm B}(T-t_i).
\]
Then
\[
\ex Y_i=A_iX_i+c_{\sin}\,\mathrm{Im}\bigl(\Psi_i\,e^{\mathrm i kX_i}\bigr),
\]
with $\ex Z_i$, $\ex U_i$ and $\Delta\ex M_i$ given by the one-step decomposition of
Lemma~\ref{B:lem:kw}.

\subsubsection{Continuous marks}
\label{B:ssec:contmark}
Take
\[
E=[0,1],\qquad \nu(\diff e)=\Lambda\,\diff e,\qquad \Lambda=\contLambdaVal,
\]
and
\[
h(e)=h_0\chi(e)-\kappa_{\mathrm F}U(e),\qquad \chi(e)=1+\tfrac{\sqrt2}{2}\cos(\pi e).
\]
The remaining data are given by \eqref{B:eq:base-scalar}.
The affine coefficient $A_i$ satisfies \eqref{B:eq:riccati}, and
\[
\ex U_i(e)=\frac{A_{i+1}h_0}{1+\kappa_{\mathrm F}A_{i+1}}\,\chi(e).
\]
We use the orthonormal cosine basis
\[
\varphi_0(e):=1,\qquad
\varphi_j(e):=\sqrt2\cos(j\pi e),\quad j\ge1,
\]
of $L^2([0,1],\diff e)$. Since
\[
\langle\chi,\varphi_2\rangle_{\Lnu}=0,
\]
the two $U$-perturbations are
\[
\delta U_i^{\parallel}(e)=\varepsilon\bigl(\sin\ex X_i+\tfrac12\bigr)\chi(e),\qquad
\delta U_i^{\perp}(e)=\varepsilon\bigl(\sin\ex X_i+\tfrac12\bigr)\varphi_2(e).
\]
The second is an off-profile mark-shape perturbation, absent in the single-mark benchmark.
When we perturb all five components together, we take
$\delta U_i=\delta U_i^\parallel+\delta U_i^\perp$.

\newcommand{\matchedCellsVal}{12}
\newcommand{\matchedBetaLowBand}{[0.901,\,0.935]}
\newcommand{\matchedBetaHighBand}{[0.197,\,0.204]}
\newcommand{\matchedBetaBase}{1}
\newcommand{\matchedBetaLow}{0.3}
\newcommand{\matchedBetaHigh}{2}
\newcommand{\matchedPredLow}{0.917}

For the integrand-effectivity comparison, we perturb $Z$, $U$ along $\chi$, or $U$ along
$\varphi_2$, with $\kappa_{\mathrm B}=1$ and
$\kappa_{\mathrm F}\in\{\matchedBetaLow,\matchedBetaBase,\matchedBetaHigh\}$. On an independent
calibration batch $A$, set
\[
\varepsilon
=\sqrt{\frac{d}{\widehat{\mathcal D}_{\pi,A}(1)}},
\qquad d=10^{-3},
\]
where $\widehat{\mathcal D}_{\pi,A}(1)$ is the measured error at unit amplitude. Freeze
$\varepsilon$ and evaluate the approximation on an independent scoring batch $B$.

\subsubsection{Normal-inverse-Gaussian benchmark}
\label{B:ssec:nigmark}
This benchmark tests the estimator under infinite jump activity with exactly sampled jump
increments on each time cell. Take $J=\nigTypesVal$ mark types and
\[
E=\{1,\dots,J\}\times(\R\setminus\{0\}),\qquad \nu=\textstyle\sum_{j=1}^{J}\delta_j\otimes\nu_0,
\]
with
\[
\nu_0(\diff e)=\frac{\delta\alpha}{\pi}\,\frac{K_1(\alpha|e|)}{|e|}\,\diff e,\qquad \alpha=\delta=1.
\]
Here $K_1$ is the modified Bessel function of the second kind of order one.
The small-jump asymptotic $K_1(z)\sim z^{-1}$ as $z\downarrow0$ gives
$\int|e|\,\nu_0(\diff e)=\infty$ and $\nu_0(\R\setminus\{0\})=\infty$, hence $\nu(E)=\infty$.
The second moments satisfy
\[
\int e^2\,\nu_0(\diff e)=\frac{\delta}{\alpha}=\nigTypeEnergyVal,
\qquad
\int_E e^2\,\nu(\diff e)=J.
\]
For each mark type, the jump integral is the $L^2$ limit of its compensated
truncations as the cutoff tends to zero.
For each mark type and time cell, draw an inverse-Gaussian variable $V$ with
\[
\E V=\frac{\delta\tau}{\alpha},
\qquad \operatorname{shape}(V)=(\delta\tau)^2.
\]
Conditional on $V$, draw the increment $S$ from $\mathcal N(0,V)$. The pairs $(V,S)$ are
independent across cells and mark types and independent of $W$. Then
\[
\E e^{\mathrm iuS}=\E e^{-u^2V/2}=\exp\bigl\{\delta\tau\bigl(\alpha-\sqrt{\alpha^2+u^2}\bigr)\bigr\}
=\exp\Bigl\{\tau\!\int\bigl(\cos(ue)-1\bigr)\,\nu_0(\diff e)\Bigr\},
\]
so $S$ has exactly the law of a typewise compensated jump increment over one cell.
Define the feature functions and driving profile by
\[
\begin{gathered}
q_j(k,e):=\mathbf 1_{\{j=k\}}e,\qquad \langle q_j,q_k\rangle_{\Lnu}=\delta_{jk},\\
\chi_{\mathrm{NIG}}(j,e)=a_je,\qquad
(a_1,a_2,a_3,a_4)=\frac{(1,\,0.7,\,0.4,\,0)}{\sqrt{1.65}}.
\end{gathered}
\]
Along the exact reference, the coefficients of Section~\ref{B:ssec:contmark},
with $\chi$ replaced by $\chi_{\mathrm{NIG}}$, give
\[
h_i(j,e)=\gamma_ia_je,\qquad \ex U_i(j,e)=u_ia_je
\]
for deterministic $\gamma_i,u_i$. The profile and off-profile perturbations are
\[
\delta U_i^\parallel(j,e)=\varepsilon\,\chi_{\mathrm{NIG}}(j,e),\qquad
\delta U_i^\perp(j,e)=\varepsilon\,\mathbf 1_{\{j=4\}}e.
\]
When we perturb all five components together, we take
$\delta U_i=\delta U_i^\parallel$. We test the off-profile direction separately.

\subsubsection{Infinite activity with a cutoff}
\label{B:ssec:infmark}
Take
\[
E=(0,1],\qquad \nu(\diff e)=\Lambda e^{-2}\,\diff e,\qquad \Lambda=\infLambdaVal,
\]
so that
\[
\nu(E)=\infty,\qquad \int_Ee\,\nu(\diff e)=\infty,\qquad \int_Ee^2\,\nu(\diff e)=\Lambda.
\]
Using the cosine basis $(\varphi_j)_{j\ge0}$ above, define
\[
q_j(e):=e\,\varphi_j(e).
\]
Then
\[
\langle q_j,q_k\rangle_{\Lnu}
=\Lambda\int_0^1\varphi_j(e)\varphi_k(e)\,\diff e
=\Lambda\delta_{jk}.
\]
With
\[
\chi_\infty(e):=e\,\chi(e),\qquad V_J:=\operatorname{span}\{q_0,\dots,q_{J-1}\},\qquad J=4,
\]
and the coefficients of Section~\ref{B:ssec:contmark} with $\chi$ replaced by $\chi_\infty$, the exact
jump integrand is
\[
\ex U_i(e)=\frac{A_{i+1}h_0}{1+\kappa_{\mathrm F}A_{i+1}}\,\chi_\infty(e).
\]
The profile and off-profile perturbations are
\[
\delta U_i^\parallel(e)=\varepsilon\bigl(\sin\ex X_i+\tfrac12\bigr)\chi_\infty(e),\qquad
\delta U_i^\perp(e)=\varepsilon\bigl(\sin\ex X_i+\tfrac12\bigr)q_2(e).
\]
When we perturb all five components together, we take
$\delta U_i=\delta U_i^\parallel+\delta U_i^\perp$.
We simulate marks $e>\eta$ explicitly, using $\eta=\infEpsVal$ in the main experiment. The
covariance of the omitted small-jump contribution on $V_J$ is
\[
\Sigma(\eta)_{jk}:=\int_0^\eta q_j(e)q_k(e)\,\nu(\diff e)
=\Lambda\int_0^\eta\varphi_j(e)\varphi_k(e)\,\diff e=O(\eta).
\]
On each time cell, we replace the omitted small-jump contribution by a centred Gaussian
increment with covariance $\tau\Sigma(\eta)$. This increment is independent of the
$\sigma$-field generated by $\F_{t_i}$, $\Delta W_i$ and the simulated jumps with
$e>\eta$ on that cell. Together with the
explicitly simulated marks $e>\eta$, the completed increments preserve the conditional second moments
of every integrand in $V_J$ with respect to the full L\'evy measure; higher moments, and hence
expectations of nonlinear functions of the state, can differ. The reference and perturbed
tuples are propagated on the same completed driver.

The off-profile direction $q_2$ has squared norm $\Lambda$ under the full measure and
$\Lambda-\Sigma(\eta)_{22}$ under the truncated measure. On the same truncated paths, the squared
amplitude $\varepsilon^2(\sin\ex X_i+\tfrac12)^2$ is common to both evaluations of the perturbation
energy. Using the full norm in $\Dpi$ therefore multiplies the effectivity $\Dpi/\Epi$ by
\[
\frac{\Lambda}{\Lambda-\Sigma(\eta)_{22}}=1+O(\eta).
\]
Covariance completion restores the full-measure variance in the $q_2$ direction.

For the continuous-time comparison, first-order path regularity is known in Brownian and
decoupled jump settings~\cite{RSZPathReg,Zhang,GeissLabart,CarmonaDelarue,BouchardElie}, but
these results do not directly yield the required $\Lnu$ time regularity of $U$ under
$U$-dependent coefficients and full-tuple law dependence. We therefore verify the required rate
for the following two benchmarks.

\begin{proposition}[Path regularity of the affine benchmarks]
\label{B:prop:benchreg}
For the linear benchmark and its weighted infinite-activity variant of
Section~\ref{B:ssec:infmark},
\[
\Reg(\tau)=O(\tau).
\]
\end{proposition}

The proof is given in Appendix~\ref{B:app:deferred}. Since the coefficients are
time-independent, $\omega\equiv0$, so \eqref{B:eq:cont-rate} applies to both benchmarks.

\paragraph{Continuous-time reference.} We use the coarse grids
\[
N\in\{\tauGridListVal\}.
\]
The remainder $\Reg(\tau)$ and the unperturbed population error are evaluated from exact second
moments. In the simulated comparison, $\Dc$ is approximated using a nested fine Euler reference at
$N_{\mathrm{fine}}=\tauFineVal$, with cell suprema taken over fine-grid nodes and $Z,U$ integrals
evaluated by fine-subcell sums.

\begin{table}[tb]\centering\footnotesize
\caption{Numerical settings for the closed-form benchmark studies.}
\label{B:tab:studies}
\begin{tabular}{@{}>{\raggedright\arraybackslash}p{0.25\linewidth}>{\raggedright\arraybackslash}p{0.17\linewidth}rc>{\raggedright\arraybackslash}p{0.29\linewidth}@{}}
\toprule
Study & Time steps $N$ & Paths per run & Seeds & Perturbation scale\\
\midrule
\closedFormSettingsRows
\bottomrule
\end{tabular}
\end{table}

\subsection{Solver-produced approximations}\label{B:sec:deepsolver}

We next evaluate the estimator on approximations produced by neural and Picard/LSMC solvers.

\paragraph{Reference-free split-sample scoring.}
For the population error-tracking studies, we use the split-sample protocol of
Table~\ref{B:tab:protocol}, with an independent pilot batch $A$.
At each checkpoint, the pilot batch estimates $\E[X]$ in the neural studies of the linear and
nonlinear-mean benchmarks and in the Picard/LSMC study. In the inventory study, batch $A$ estimates
both $\E[X]$ and $\E[Y]$. We freeze these
estimates as the statistic path $s^A$. A fresh scoring batch $B$ is then generated from
initial data and driving noise independent of the training and pilot data. The same frozen
path $s^A$ is used both to propagate these trajectories and to evaluate their residuals. Thus
the batch estimates $\Epi^{\mathrm{plug}}(s^A)$, and the reported effectivities are
$\Dpi/\Epi^{\mathrm{plug}}$.

Conditional on the trained parameters and the pilot statistics $s^A$, all laws and
expectations refer to the fresh evaluation trajectories. Corollary~\ref{B:cor:statbridge} relates the
plug-in score $\Epi^{\mathrm{plug}}(s^A)$ to the own-law estimator $\Epi$ of
Theorem~\ref{B:thm:discrete}.

Because the same frozen pilot statistics are used to propagate $\ap X$ and to evaluate the forward
coefficients, $C_i^X\equiv0$ analytically under this protocol and appears numerically only at
round-off level. The imposed initial condition also gives $r_0^X=0$. Hence
\[
\Epi^{\mathrm{plug}}
=\max_{i<N}\E|C_i^Y|^2+\E|r_T^Y|^2.
\]
The neural and Picard/LSMC solvers fit only $(\ap Y,\ap Z,\ap U)$.
The forward component $\ap X$ is generated by the forward recursion, and for scoring we take
$\ap M\equiv0$.

To test sensitivity to the pilot statistics, batch $C$ of $\candStatsBCVal$ trajectories is
generated for each checkpoint using the same trained candidate and frozen statistic path $s^A$ as
batch $B$. The resulting statistics $s^C$ are then used to rescore batch $B$, and this is compared
with scoring the same batch using $s^A$ (Corollary~\ref{B:cor:statbridge}). Since the batch-$B$
trajectories are unchanged, $\Dpi$ is unchanged; we report the corresponding change in the
plug-in effectivity $\Dpi/\Epi^{\mathrm{plug}}$.

\begin{table}[tb]\centering\small
\caption{Evaluation settings for the linear, nonlinear-mean and inventory population error-tracking studies.}
\label{B:tab:protocol}
\begin{tabular}{@{}ll@{}}
\toprule
Quantity & Setting\\
\midrule
grid & $N=40$\\
pilot batch & $B_A=800{,}000$\\
scoring batch & $B_B=100{,}000$\\
scoring convention & split-sample, $\ap M\equiv0$\\
\bottomrule
\end{tabular}
\end{table}

\subsubsection{Neural solver}\label{B:ssec:neural}

On the linear benchmark of Section~\ref{B:subsec:benchmark-protocol}, we train an FBSNN-type
residual solver adapted from~\cite{Raissi2018} (see also~\cite{EHanJentzen}).
The network takes $(t,x)$ as input and has three hidden layers of $\deepHiddenVal$ $\tanh$ units each.
We use Adam with learning rate $10^{-3}$ and batches of $2048$ paths.
For each of $\deepNseedsVal$ seeds, we evaluate the network after
\[
K\in\{0,25,50,100,200,400,800,1600,3200\}
\]
updates. Define the one-step backward residual by
\[
\rho_j^Y
:=\ap Y_{j+1}-\ap Y_j+\ap f_j^{\mathrm{tr}}\tau
-\ap Z_j\Delta W_j-\int_E\ap U_j\,\Delta\Nt_j.
\]
The training objective is
\[
\mathcal L
=\E\bigl|\ap Y_N-g(\ap X_N,\rho_N^{\mathrm{tr}})\bigr|^2
+\sum_{j=0}^{N-1}\tau^{-1}\E|\rho_j^Y|^2.
\]
For a training batch of $n$ paths, set
\[
\mu_i^{\mathrm{tr}}
:=\frac1n\sum_{\ell=1}^{n}
\delta_{(\ap X_i^{(\ell)},\,\ap Y_i^{(\ell)},\,\ap Z_i^{(\ell)},\,\ap U_i^{(\ell)})},
\qquad
\rho_N^{\mathrm{tr}}
:=\frac1n\sum_{\ell=1}^{n}\delta_{\ap X_N^{(\ell)}}.
\]
Evaluating \eqref{B:eq:Epi} with these laws gives $\Epi^{\mathrm{train}}$.

Since $C_i^Y=\sum_{j=0}^i\rho_j^Y$, Cauchy--Schwarz gives
\[
\max_{i<N}\E|C_i^Y|^2
\le T\sum_{j=0}^{N-1}\tau^{-1}\E|\rho_j^Y|^2.
\]
Since $r_0^X=0$ and $C_i^X=0$, including the terminal residual yields
\[
\Epi^{\mathrm{train}}\le\max(1,T)\,\mathcal L.
\]

\subsubsection{Nonlinear mean dependence}\label{B:ssec:nonlinear-mean}

To test nonlinear dependence on the mean field, let
\[
m_i:=\E[X_i],\qquad T=1,\qquad X_0\sim\mathcal N(0.75,1),\qquad \nu(E)=1,
\]
and take
\[
\begin{gathered}
b=0.2m-y,\qquad \sigma=0.5-z,\qquad h=0.5-u,\\
f=x+0.2m+\lambda\tanh(m/0.5),\qquad g=x+0.2m,\qquad \lambda\in\{0,1\}.
\end{gathered}
\]
Since $\tanh$ is increasing,
\[
-\lambda(m_1-m_2)\bigl[\tanh(m_1/0.5)-\tanh(m_2/0.5)\bigr]\le0,
\]
so the nonlinear mean term preserves the monotonicity bound of the linear benchmark.

The exact discrete solution is
\[
\ex Y_i=A_i\ex X_i+C_i,\qquad A_i\equiv1,\qquad \ex Z_i=\ex U_i=\tfrac14,
\]
with
\[
\begin{aligned}
C_N&=0.2m_N,\\
C_i&=\frac{(0.2A_{i+1}+0.2)m_i\tau+\lambda\tau\tanh(m_i/0.5)+C_{i+1}}{1+A_{i+1}\tau},\\
m_0&=0.75,\\
m_{i+1}&=m_i+\tau\bigl[0.2m_i-(A_im_i+C_i)\bigr].
\end{aligned}
\]
The coupled $(m_i,C_i)$ system is solved to equation residual below $10^{-11}$. For $\lambda=1$,
the nonlinear term remains non-negligible throughout the grid, with
\[
\lambda\tanh(m_i/0.5)>0.3.
\]
We train the FBSNN of Section~\ref{B:ssec:neural} for five matched seeds at $\lambda=0$ and $1$,
saving checkpoints at $K\in\{0,100,400,800,1600,3200\}$.

\subsubsection{Picard/LSMC}

On the linear benchmark of Section~\ref{B:subsec:benchmark-protocol}, we use Picard iteration with
least-squares Monte Carlo regression. For each seed, we reuse the same $\lsmcNtrainVal$ noise paths across Picard iterations and
initialise $Y^{(0)}=Z^{(0)}=U^{(0)}=0$. At each iteration, we simulate $X$ forward using the current
fields, then fit new $Y,Z,U$ fields by backward regression.

Since the exact decoupling fields are affine, we fit $(\ap Y_i,\ap Z_i,\ap U_i)$ at each node by
ordinary least squares using a constant and the standardised state. In the scalar single-mark benchmark, the backward step and
Lemma~\ref{B:lem:kw} give
\[
\begin{aligned}
Y_i&=\Ei[Y_{i+1}]+f_i\tau,\\
Z_i&=\tau^{-1}\Ei[Y_{i+1}\Delta W_i],\\
U_i&=\bigl(\tau\nu(E)\bigr)^{-1}\Ei[Y_{i+1}\Delta\Nt_i].
\end{aligned}
\]

Let $\widetilde\Phi^{(k+1)}$ denote the new regression fit for $\Phi=(Y,Z,U)$.
We use the damped update
\[
\Phi^{(k+1)}
=(1-\theta)\Phi^{(k)}+\theta\,\widetilde\Phi^{(k+1)},
\qquad \theta=\tfrac34.
\]
This was the largest tested value of $\theta$ that passed the stability check on all development seeds.
For each of $\lsmcNseedsVal$ separate seeds, we record the iterates $k\in\{0,1,2,3,4,6,8,12,20,30\}$.

\subsubsection{Mean-field inventory control}
\label{B:ssec:lqinv}
We consider a scalar additive-noise specialisation of the linear--quadratic mean-field jump-control
problem of Tang and Meng~\cite{TangMeng}. The controlled state and objective are
\[
dX_t=\bigl[-\kappa X_t-\bar a m_t+u_t\bigr]\,dt+\sigma_0\,dW_t+h_0\,d\Nt_t,\qquad m_t=\E[X_t],
\]
and
\[
J(u)=\E\Big[\int_0^T\bigl(Q(X_t-m_t)^2+\bar Q\,m_t^2+Ru_t^2\bigr)\,dt
+Q_T(X_T-m_T)^2+\bar Q_Tm_T^2\Big].
\]
Assume $R>0$ and $Q,\bar Q,Q_T,\bar Q_T\ge0$.
Here $X_t$ is the dealer inventory and $u_t$ the trading rate. With the adjoint normalisation
$p_t=2Y_t$, the mean-field Pontryagin principle gives
\[
u_t^\star=-R^{-1}Y_t
\]
and the optimality system
\[
\begin{aligned}
dX_t&=\bigl[-\kappa X_t-\bar a m_t-R^{-1}Y_t\bigr]\,dt+\sigma_0\,dW_t+h_0\,d\Nt_t,\\
-dY_t&=\bigl[-\kappa Y_t-\bar a\,\E[Y_t]+Q(X_t-m_t)+\bar Q\,m_t\bigr]\,dt-Z_t\,dW_t-U_t\,d\Nt_t,\\
Y_T&=Q_T(X_T-m_T)+\bar Q_Tm_T.
\end{aligned}
\]
The backward mean term $-\bar a\,\E[Y_t]$ is the mean-field adjoint of the forward crowding term
$-\bar a m_t$.

The exact continuous-time solution is affine,
\[
Y_t=P_t(X_t-m_t)+\Pi_tm_t,\qquad Z_t=P_t\sigma_0,\qquad U_t=P_th_0.
\]
For the Euler scheme, the exact discrete reference is determined by
\[
\begin{aligned}
P_N&=Q_T,&
P_i&=\frac{(1-\kappa\tau)P_{i+1}+\tau Q}{1+\kappa\tau+\tau P_{i+1}/R},\\
\Pi_N&=\bar Q_T,&
\Pi_i&=\frac{(1-(\kappa+\bar a)\tau)\Pi_{i+1}+\tau\bar Q}{1+(\kappa+\bar a)\tau+\tau\Pi_{i+1}/R}.
\end{aligned}
\]
The mean satisfies
\[
m_{i+1}=m_i\Bigl[1-\tau\Bigl(\kappa+\bar a+\frac{\Pi_i}{R}\Bigr)\Bigr],\qquad m_0=\E[X_0],
\]
and
\[
\ex Y_i=P_i(\ex X_i-m_i)+\Pi_im_i,\qquad \ex Z_i=P_{i+1}\sigma_0,\qquad \ex U_i=P_{i+1}h_0.
\]

The reported runs use
\[
\begin{gathered}
T=1,\qquad \kappa=0.5,\qquad \bar a=0.3,\qquad R=Q=1,\qquad \bar Q=0.8,\\
Q_T=0.9,\qquad \bar Q_T=0.7,\qquad \sigma_0=0.4,\qquad h_0=0.3,\qquad \nu(E)=1,
\end{gathered}
\]
with
\[
X_0\sim\mathcal N(0.5,0.2^2),\qquad N=40.
\]
Appendix~\ref{B:app:inv-admiss} verifies Assumptions~\ref{B:ass:integ}--\ref{B:ass:mono} with
$\beta_1=0$, $\beta_2=0.8$ and $\alpha=0.7$.

The FBSNN of Section~\ref{B:ssec:neural} is trained on this benchmark with the same budgets and
seeds. For split-sample evaluation, the pilot freezes both $\E[X]$ and $\E[Y]$, the latter being
required by the backward coefficient $-\bar a\,\E[Y_t]$.

\paragraph{Control regret.} For the control comparison, each checkpoint is evaluated as the
closed-loop feedback
\[
\ap u_i=-R^{-1}\ap Y_i
\]
on a system of $p$ interacting particles with independent initial data and driving noises. Write
\[
X_i=(X_i^1,\dots,X_i^p)^{\!\top},\qquad u_i=(u_i^1,\dots,u_i^p)^{\!\top},
\]
and, for $v\in\R^p$,
\[
\bar v:=\frac1p\mathbf 1^{\!\top}v,\qquad \tilde v:=v-\bar v\,\mathbf 1,\qquad
|v|_p^2:=\frac1p\sum_{a=1}^p|v^a|^2.
\]
Replacing $\E[X_i]$ and $\E[Y_i]$ by the empirical means $\bar X_i$ and $\bar Y_i$ gives the
corresponding particle FBSDE in $\R^p$.

Its Brownian and jump integrands are $p\times p$ matrices. For $A\in\R^{p\times p}$, we use the
normalized Frobenius norm
\[
|A|_{p,F}^2:=\frac1p|A|_F^2,
\]
with jump-integrand norm $\nu(E)|A|_{p,F}^2$. We write $\Epip$ and $\Dpip$ for the resulting
estimator and discrete error.

Its Euler solution is
\[
\ex Y_i=P_i\tilde X^\pi_i+\Pi_i\bar X^\pi_i\mathbf 1.
\]
Define
\[
\mathsf P_i:=P_iI+(\Pi_i-P_i)\frac{\mathbf 1\mathbf 1^{\!\top}}{p}.
\]
Then
\[
\ex Z_i=\sigma_0\mathsf P_{i+1},\qquad \ex U_i=h_0\mathsf P_{i+1}.
\]
For the neural candidate, $\ap Z_i$ and $\ap U_i$ are diagonal matrices formed from the
particles' $Z$- and $U$-outputs, respectively.

The left-endpoint discretisation of the control objective is
\[
J_{\pi,p}(u):=\E\Big[\sum_{i=0}^{N-1}\tau\big(Q|\tilde X_i|_p^2+\bar Q\,\bar X_i^2+R|u_i|_p^2\big)
+Q_T|\tilde X_N|_p^2+\bar Q_T\,\bar X_N^2\Big].
\]
Its minimiser over adapted, square-integrable controls is
\[
u^\star_i=-k_i\tilde X_i-\bar k_i\bar X_i\mathbf 1,
\]
where, with
\[
\ell:=1-\kappa\tau,\qquad \bar\ell:=1-(\kappa+\bar a)\tau,
\]
\[
\begin{aligned}
S_N&=Q_T,&
S_i&=\tau Q+\frac{\ell^2RS_{i+1}}{R+\tau S_{i+1}},&
k_i&=\frac{\ell S_{i+1}}{R+\tau S_{i+1}},\\
\bar S_N&=\bar Q_T,&
\bar S_i&=\tau\bar Q+\frac{\bar\ell^{\,2}R\bar S_{i+1}}{R+\tau\bar S_{i+1}},&
\bar k_i&=\frac{\bar\ell\,\bar S_{i+1}}{R+\tau\bar S_{i+1}}.
\end{aligned}
\]
Define the control regret and Euler-FBSDE baseline by
\[
\mathcal R_{\pi,p}(u):=J_{\pi,p}(u)-J_{\pi,p}(u^\star),
\qquad
B_{\pi,p}:=\mathcal R_{\pi,p}(u^{\mathrm{FB}}),
\]
where $u^{\mathrm{FB}}_i:=-R^{-1}\ex Y_i$ is the Euler-FBSDE feedback.
In this benchmark, the Euler-FBSDE feedback gains differ from those of the minimiser of
$J_{\pi,p}$, so even the exact Euler solution incurs positive regret $B_{\pi,p}$.

\begin{proposition}[Control regret]\label{B:prop:regret}
Assume that, for $0\le i\le N$,
\[
P_i,\ \Pi_i\ge0,\qquad \varphi_i:=1-\tau(\kappa+P_i/R)\in[0,1],\qquad
\bar\varphi_i:=1-\tau(\kappa+\bar a+\Pi_i/R)\in[0,1].
\]
\begin{enumerate}[label=\textup{(\roman*)},ref=\roman*,leftmargin=2.4em,itemsep=3pt]
\item\label{B:item:sq} For every adapted, square-integrable control $u$,
\[
\mathcal R_{\pi,p}(u)=\E\sum_{i=0}^{N-1}\tau\Big[(R+\tau S_{i+1})|\tilde u_i+k_i\tilde X_i|_p^2
+(R+\tau\bar S_{i+1})|\bar u_i+\bar k_i\bar X_i|^2\Big].
\]
\item\label{B:item:regbound} Let $\ap\Theta$ be a square-integrable tuple with $\ap X_0=X_0$ whose
forward component follows the particle Euler recursion under $\ap u=-R^{-1}\ap Y$, with adapted
$\ap Z,\ap U$ and admissible $\ap M$. Then
\[
\Big|\sqrt{\mathcal R_{\pi,p}(\ap u)}-\sqrt{B_{\pi,p}}\Big|\le\sqrt{C_{\mathrm{inv}}\,\Epip(\ap\Theta)},
\qquad C_{\mathrm{inv}}:=\frac{10\,T\,w_{\max}}{R^2}\,(1+T\Delta_k)^2,
\]
with
\[
w_{\max}:=\max_i\max\{R+\tau S_{i+1},R+\tau\bar S_{i+1}\},\qquad
\Delta_k:=\max_i\max\{|k_i-P_i/R|,|\bar k_i-\Pi_i/R|\}.
\]
\item\label{B:item:Btau} $\Delta_k=O(\tau)$, and $B_{\pi,p}=O(\tau^2)$ uniformly in $p$.
\end{enumerate}
\end{proposition}

\noindent The proof is given in Appendix~\ref{B:app:regret}. At
fixed $\tau$ and $p$, part~(\ref{B:item:regbound}) gives
\[
\Epip(\ap\Theta)\to0\quad\Longrightarrow\quad\mathcal R_{\pi,p}(\ap u)\to B_{\pi,p}.
\]
Thus an increasingly accurate solution of the Euler-FBSDE approaches the Euler-FBSDE regret
$B_{\pi,p}$, rather than necessarily zero regret. By parts~(\ref{B:item:regbound})
and~(\ref{B:item:Btau}),
\[
\mathcal R_{\pi,p}(\ap u)\le C\bigl(\Epip(\ap\Theta)+\tau^2\bigr),
\]
with $C$ independent of $p$.

\begin{remark}[Estimator accuracy and control regret]\label{B:rem:regret-order}
The estimator and control regret measure different aspects of the approximation. This is seen
explicitly by interpolating between the Euler-FBSDE feedback and the discrete-optimal feedback.

Let $\Theta^\star$ be the tuple generated by the discrete-optimal feedback $u^\star$ on the same
initial data and driving noise as $\ex\Theta$, with
\[
Y_i^\star:=-Ru_i^\star,\qquad i<N,
\qquad
Y_N^\star:=Q_T\tilde X_N^\star+\bar Q_T\bar X_N^\star\mathbf 1.
\]
By Lemma~\ref{B:lem:kw},
\[
Y_{i+1}^\star-\Ei[Y_{i+1}^\star]
=
Z_i^\star\Delta W_i
+\int_E U_i^\star(e)\,\Delta\Nt_i(\diff e)
+\Delta M_i^\star.
\]

For $0\le\theta\le1$, set
\[
\Theta^{(\theta)}
:=
(1-\theta)\ex\Theta+\theta\Theta^\star,
\qquad
u^{(\theta)}
:=
-R^{-1}Y^{(\theta)}
=
(1-\theta)u^{\mathrm{FB}}+\theta u^\star.
\]
By linearity, $\Theta^{(\theta)}$ satisfies the hypotheses of
Proposition~\ref{B:prop:regret}(\ref{B:item:regbound}). The estimator residuals vanish for
$\ex\Theta$, and both squared terms in Proposition~\ref{B:prop:regret}(\ref{B:item:sq}) vanish for
$u^\star$. Hence
\[
\Epip\bigl(\Theta^{(\theta)}\bigr)
=
\theta^2\Epip(\Theta^\star),
\qquad
\mathcal R_{\pi,p}\bigl(u^{(\theta)}\bigr)
=
(1-\theta)^2B_{\pi,p}.
\]
Thus, as $\theta$ decreases from $1$ to $0$, the estimator falls from $\Epip(\Theta^\star)$ to zero
while the regret rises from zero to $B_{\pi,p}$.
\end{remark}

For each training seed, we use two disjoint evaluation cohorts, each consisting of $\ictlEnsVal$
independent interacting systems of $p=\ictlNpartVal$ particles. The same cohorts are used across
that seed's checkpoints, so comparisons across training budgets are paired. Within each cohort,
$\Epip$, $\Dpip$ and $\mathcal R_{\pi,p}$ are computed from the same trajectories, with regret
evaluated using Proposition~\ref{B:prop:regret}(\ref{B:item:sq}).

\subsubsection{State-nonlinear benchmark}\label{B:ssec:snl}
To test nonlinear coupling in both directions, take
\[
b(y)=-\tanh y,\qquad \sigma=\sigma_0,\qquad h=h_0,\qquad f(x)=g(x)=x+\tfrac12\tanh x,
\]
with
\[
\sigma_0=h_0=\tfrac12,\qquad \nu(E)=1,\qquad T=1,\qquad X_0\sim\mathcal N(0,1),\qquad N=40.
\]
Since $\tanh$ is increasing,
\[
\delta b\,\delta Y-\delta f\,\delta X\le-|\delta X|^2,\qquad \delta g\,\delta X\ge|\delta X|^2,
\]
so Assumption~\ref{B:ass:mono} holds with
\[
G=1,\qquad \beta_1=0,\qquad \beta_2=1,\qquad \alpha=1.
\]

The exact discrete solution is Markov,
\[
\ex Y_i=v_i(\ex X_i),\qquad v_N=g.
\]
For independent
\[
\xi\sim\mathcal N(0,1),\qquad K\sim\mathrm{Poisson}(\nu(E)\tau),
\]
define
\[
H_i(s):=\E\Bigl[v_{i+1}\bigl(s+\sigma_0\sqrt\tau\,\xi+h_0(K-\nu(E)\tau)\bigr)\Bigr].
\]
Then $v_i(x)$ is determined by
\[
v_i(x)=H_i\bigl(x-\tau\tanh v_i(x)\bigr)+\tau f(x).
\]
Indeed, suppose inductively that $v_{i+1}$ is continuous, increasing and of at most linear
growth, as holds for $v_N=g$. Then $H_i$ is finite, continuous and increasing. For
\[
F_i(x,y):=H_i(x-\tau\tanh y)+\tau f(x)-y,
\]
the map $F_i(x,\cdot)$ is continuous and strictly decreasing, with
\[
\lim_{y\to-\infty}F_i(x,y)=+\infty,
\qquad
\lim_{y\to+\infty}F_i(x,y)=-\infty.
\]
Hence it has a unique zero $v_i(x)$, which depends continuously on $x$. Moreover, if
$x_2>x_1$, then
\[
F_i(x_2,y)>F_i(x_1,y)
\]
for every $y$, since $H_i$ is increasing and $f$ is strictly increasing. As $F_i(x_2,\cdot)$ is
strictly decreasing, its zero satisfies
\[
v_i(x_2)>v_i(x_1).
\]
Since $|\tanh|\le1$, the defining equation also shows that $v_i$ has at most linear growth. Thus
$v_i$ is continuous, increasing and of at most linear growth, closing the backward induction.

We compute a high-accuracy numerical reference for this Euler scheme with the settings in
Table~\ref{B:tab:snlref}, where $\Delta x$ is the spatial grid spacing, $[-L,L]$ the
interpolation domain and $Q$ the number of Gauss--Hermite nodes.

\begin{table}[tb]\centering\small
\caption{Numerical-reference settings for the state-nonlinear benchmark.}
\label{B:tab:snlref}
\begin{tabular}{@{}lccccc@{}}
\toprule
 & $\Delta x$ & $L$ & $Q$ & Poisson tail & root residual\\
\midrule
reported & $\snlRefDxVal$ & $\snlRefLVal$ & $\snlRefQVal$ & $10^{-14}$ & $10^{-12}$\\
refined & $\snlFinerDxVal$ & $\snlFinerLVal$ & $\snlFinerQVal$ & $\snlFinerTailVal$ & $\snlFinerRootVal$\\
\bottomrule
\end{tabular}
\end{table}

We compute $H_i$ as in Appendix~\ref{B:app:snlref} and obtain $\ex Z_i$ and $\ex U_i$ by
Lemma~\ref{B:lem:kw}.

The reference tuple is propagated using its own forward recursion. Its martingale component is
reconstructed pathwise from
\[
\Delta\ex M_i
=
\ex Y_{i+1}-\ex Y_i+\tau f(\ex X_i)
-\ex Z_i\Delta W_i-\ex U_i\Delta\Nt_i,
\]
and is nonzero for this benchmark. At each time step, we check the backward equation and the
projection identities for $\ex Z_i$ and $\ex U_i$ at 33 test states using Gauss--Legendre
quadrature. The largest absolute residual is below $8\times10^{-10}$.

The construction reproduces the analytic linear case to numerical precision. As a refinement
check, we repeat the construction with the refined settings in Table~\ref{B:tab:snlref}. Across
all reported checkpoints, replacing the reported reference by the refined reference changes
the computed $\Dpi$ by less than $\snlRefDirectVal$ relatively.

We use the FBSNN of Section~\ref{B:ssec:neural} with $\snlSeedsVal$ seeds and checkpoint
budgets $K\in\{0,100,400,800,1600,3200\}$. Since the coefficients are law-independent, $\Epi^{\mathrm{plug}}=\Epi$.

\section{Numerical results}\label{B:sec:results}

We report effectivity as $\Dpi/\Epi$. Reported effectivity and ratio bands are minima and maxima
over the cells specified for each study: grids, perturbation directions and amplitudes for the
closed-form sweeps, and seeds and
training budgets for the solver-produced approximations. Solver curves show means across
training seeds with one-standard-error bands ${}\pm s/\sqrt{n_{\mathrm{seed}}}$; target-error
bars show seed minima and maxima. Detrended correlations are Pearson correlations of the log
quantities after centring across seeds at each training budget. Training-budget axes are linear
from $K=0$ to the first positive budget and logarithmic thereafter.

\subsection{Closed-form benchmarks}\label{B:sec:results-closed}

\subsubsection{Discrete two-sided equivalence}
\label{B:ssec:twosided}
\begin{figure}[tb]\centering
\incfig[0.58\linewidth]{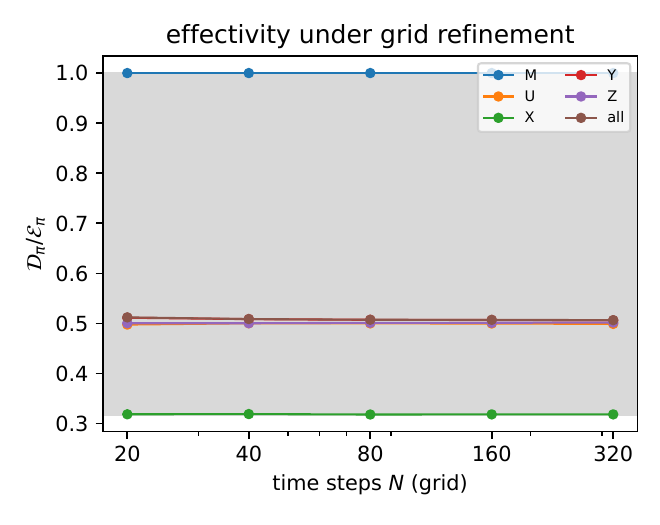}
\caption{Discrete two-sided equivalence (Theorem~\ref{B:thm:discrete}) on the linear
benchmark. Effectivity $\Dpi/\Epi$ across grids when $X,Y,Z,U,M$ are perturbed individually or
all together.}
\label{B:fig:twosided}
\end{figure}

For a fixed perturbation direction in the linear benchmark, both $\Dpi$ and $\Epi$ scale as
$\varepsilon^2$, so the effectivity is independent of amplitude. The measured effectivity depends on which components are perturbed
and, for a fixed perturbation, changes little under grid refinement (Figure~\ref{B:fig:twosided}).
Across the five grids from $N=\gridMinVal$ to $N=\gridMaxVal$, the minimum measured
effectivity stays between $\gridUniformLoVal$ and $\gridUniformHiVal$. The maximum rounds to
$\gridMaxEffVal$ on every grid.

Stronger feedback widens the measured band in these runs: for
$\kappa_{\mathrm F}=\kappa_{\mathrm B}=\rv{apocert.betaWeak},\rv{apocert.betaOne},\rv{apocert.betaStrong}$ the
bands are $\weakBandVal$, $\twosidedBandVal$ and $\strongBandVal$, with one seed each for the
smaller and larger values. Replacing the control-variate mean by the
empirical mean of the same paths introduces sampling error that persists as the perturbation
amplitude decreases. On grids $\lawCheckGridLoVal$--$\lawCheckGridHiVal$, the lower edge of the
measured effectivity band falls from $\oracleLawLoVal$ to $\particleLawLoVal$, with the minimum
attained at the smallest amplitude tested.
We repeat the base scalar sweep of Table~\ref{B:tab:studies} with two additional seeds, using
$\npathsVal$ paths per run. For each parameter combination, we express the range of the three
effectivity values as a percentage of their mean. The largest such range is below
$\seedSpreadPctVal\%$.

For the sine-terminal benchmark, the measured effectivity against the full exact tuple lies in
$\sineBandVal$. The exact discrete solution has a nonzero orthogonal martingale, whose terminal
energy $\E|\ex M_N|^2$ is estimated as $\sineMmassVal$ at $N=80$.

\subsubsection{Rank, law and jump structure}\label{B:ssec:rank}\label{B:ssec:fulllaw}%
\label{B:ssec:infmark-res}
\paragraph{Rank geometry.} Across all three rank regimes of Assumption~\ref{B:ass:mono}, the measured
effectivity ranges from roughly $0.3$ to $1$. Within each rectangular case, kernel perturbations
give effectivities closer to one than range perturbations at the same target squared error
(Table~\ref{B:tab:rankgeom}).

For $m<n$, the effectivity is about $\rankKerMltnVal$ for state perturbations in $\ker G$ and $\rankRanMltnVal$ for
those in $\operatorname{ran}(\Gstar)$. A perturbation in $\ker G$ leaves $b$, $f$ and $g$
unchanged, so only the initial and cumulative forward residuals contribute:
\[
\Epi
=\E|\delta X_0|^2+\max_{0\le i<N}\E|\delta X_{i+1}-\delta X_0|^2.
\]
For $m>n$, the same ordering holds for value perturbations in $\ker(\Gstar)$ and
$\operatorname{ran}(G)$.

\begin{table}[tb]\centering
\caption{Effectivity for target-error range and kernel perturbations in the two rectangular
cases, aggregated over the grid and perturbation-amplitude sweep.}
\label{B:tab:rankgeom}
\begin{tabular}{@{}lllcc@{}}
\toprule
case & component & direction & $\min\Dpi/\Epi$ & $\max\Dpi/\Epi$\\
\midrule
\rankGeometryTableBody
\bottomrule
\end{tabular}
\end{table}

\paragraph{Full-tuple law dependence.} Across all perturbations in Table~\ref{B:tab:fulllaw},
$\Dpi/\Epi$ remains in $\flBandVal$. Perturbing $Z$ or $U$ alone also changes the forward
drift through $\lambda_Z\E[\delta Z_i]$ or $\lambda_U\E[\int_E \delta U_i(e)\,\nu(\diff e)]$,
respectively. In both cases the measured effectivity stays close to $\tfrac12$ across the tested
grids.

\begin{table}[tb]
\centering\small
\caption{Full-tuple-law benchmark of Section~\ref{B:ssec:fulllaw-setup}: per-component
effectivity. Each row perturbs one component (or all jointly) and reports the effectivity range
$\Dpi/\Epi$ over the $N\times\varepsilon$ grid.}
\label{B:tab:fulllaw}
\begin{tabular}{lcc}
\toprule
Perturbed component & $\min\,\Dpi/\Epi$ & $\max\,\Dpi/\Epi$\\
\midrule
\flGridBody
\bottomrule
\end{tabular}
\end{table}

\paragraph{Jump structure.} Across the four jump structures (single mark, continuous marks, NIG
and the covariance-completed cutoff), the measured effectivity band is $\twosidedBandVal$, and
the unperturbed tuples score at round-off. Off-profile effectivities lie in
$[\contOffLoVal,\,\contOffHiVal]$ for continuous marks, $\nigOffBandVal$ for NIG and
$[\infOffLoVal,\,\infOffHiVal]$ for the covariance-completed cutoff.

At $\eta=\infTrapEpsVal$, the off-profile effectivity is $\infTrapMatchedVal$ when both the
driver and the error norm use the truncated measure. Evaluating $\Dpi$ with the full measure on
the same paths raises the effectivity to $\infTrapMismatchVal$, a factor of $\infTrapFacVal$ as
predicted in Section~\ref{B:ssec:infmark}. With covariance completion, the full-measure
effectivity is $\infTrapCompletedVal$. At the main cutoff $\eta=10^{-3}$, the effectivity for each parameter
combination differs from its base-seed value by less than $2\%$ in both repeat runs.

\begin{figure}[tb]\centering
\incfig[0.58\linewidth]{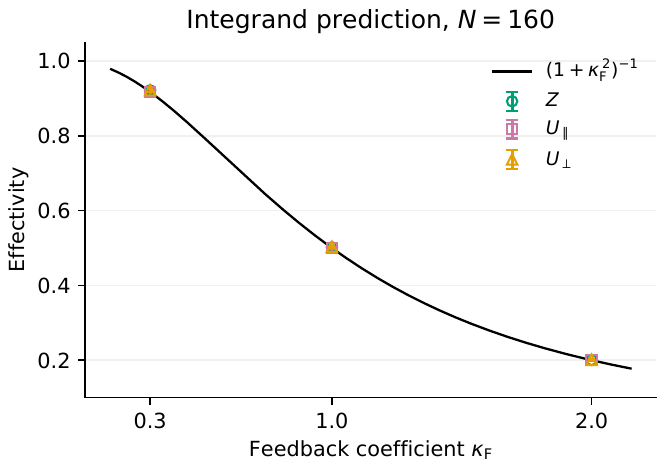}
\caption{Integrand effectivity on the continuous-mark benchmark at $N=160$, target squared error
$d=10^{-3}$ and $\kappa_{\mathrm B}=1$. The curve is the exact population prediction
$(1+\kappa_{\mathrm F}^2)^{-1}$. Points show means over three seeds; bars span their minimum and
maximum.}
\label{B:fig:matched}
\end{figure}

On the continuous-mark benchmark, the measured effectivities for perturbations in $Z$,
$U^\parallel$ and $U^\perp$ are close to the common prediction $(1+\kappa_{\mathrm F}^2)^{-1}$ at
all three feedback strengths (Figure~\ref{B:fig:matched}).

\subsubsection{Continuous-time estimate}
\label{B:ssec:tau}
\begin{figure}[tb]\centering
\incfig{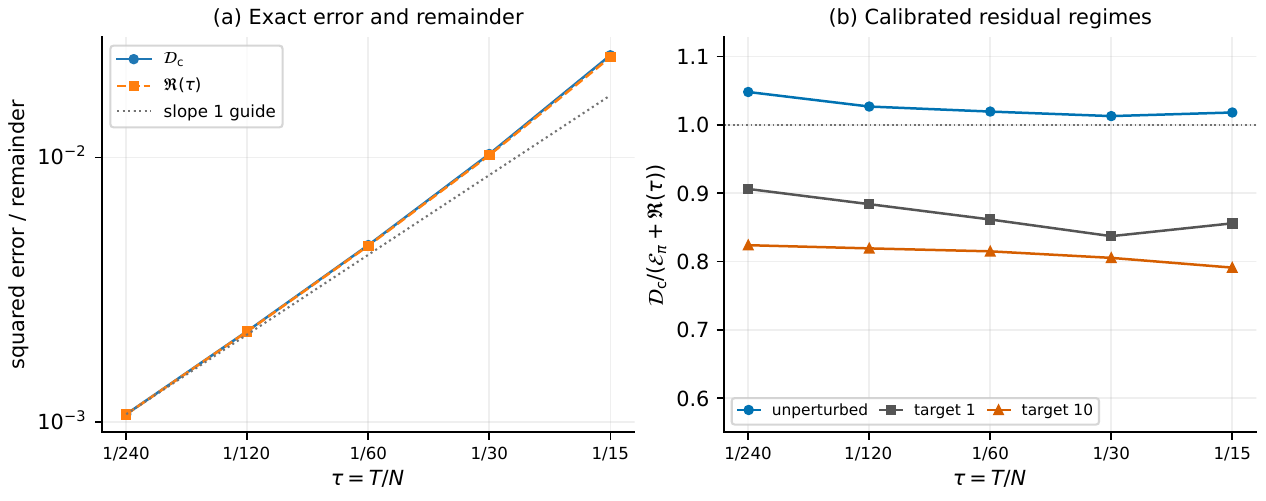}
\caption{Continuous-time estimate (Theorem~\ref{B:thm:cont}) on the coarse grids
$N\in\{\tauGridListVal\}$. (a) Exact population error and remainder for the exact discrete tuple.
(b) Normalised error proxy for the unperturbed tuple and independently calibrated targets $r=1$
and $r=\regimeEstTargetVal$, using a fine Euler reference with $N_{\mathrm{fine}}=\tauFineVal$.}
\label{B:fig:tau}
\end{figure}

Figure~\ref{B:fig:tau}(a) compares the continuous-time error $\Dc$ with the path-regularity
remainder on the exact discrete tuple ($\Epi=0$), both evaluated from exact second moments. Their
ratio lies in $\dcrBandVal$, and the two quantities \tauSlopeFitted, consistent with
Proposition~\ref{B:prop:benchreg}. Here $\omega\equiv0$. Since $Z_t=U_t=\tfrac14$, only the
increments of $X$ and $Y$ contribute to $\Reg(\tau)$.

We perturb all five components together as in Section~\ref{B:subsec:benchmark-protocol}. For
each grid and target $r\in\{1,\regimeEstTargetVal\}$, an independent calibration batch $A$ sets
\[
\varepsilon_A=\sqrt{\frac{r\,\Reg(\tau)}
{\widehat{\mathcal E}_{\pi,A}(1)}},
\]
where $\widehat{\mathcal E}_{\pi,A}(1)$ is the measured estimator at unit amplitude. We freeze
$\varepsilon_A$ and evaluate the approximation on an independent scoring batch $B$. Each batch
contains $\tauPathsVal$ paths.

The achieved estimator/remainder ratios lie in $\regimeBalAchievedBandVal$ and
$\regimeEstAchievedBandVal$ for targets $1$ and $\regimeEstTargetVal$, respectively. The
corresponding normalised errors lie in $[\regimeBalLoVal,\regimeBalHiVal]$ and
$[\regimeEstLoVal,\regimeEstHiVal]$ (Figure~\ref{B:fig:tau}(b)).

In a paired check over three seeds, refining $N_{\mathrm{fine}}$ from $960$ to $3840$ changes
each normalised-error band edge by at most $\refineWorstDevVal$.

\subsection{Solver-produced approximations}\label{B:sec:results-solver}

\subsubsection{Neural solver benchmarks}

\paragraph{Linear benchmark.} Across $\deepNseedsVal$ seeds, the seed-mean $\Epi^{\mathrm{plug}}$ falls from
$\deepEpiStartVal$ to $\deepEpiEndVal$ by budget $\deepBudgetMaxVal$, while $\Dpi$ reaches
$\deepDpiEndVal$. The realised effectivity band is $\deepBandVal$, and the detrended correlation is
$\stratValDetrendCorrVal$.

\begin{figure}[tb]\centering
\incfig{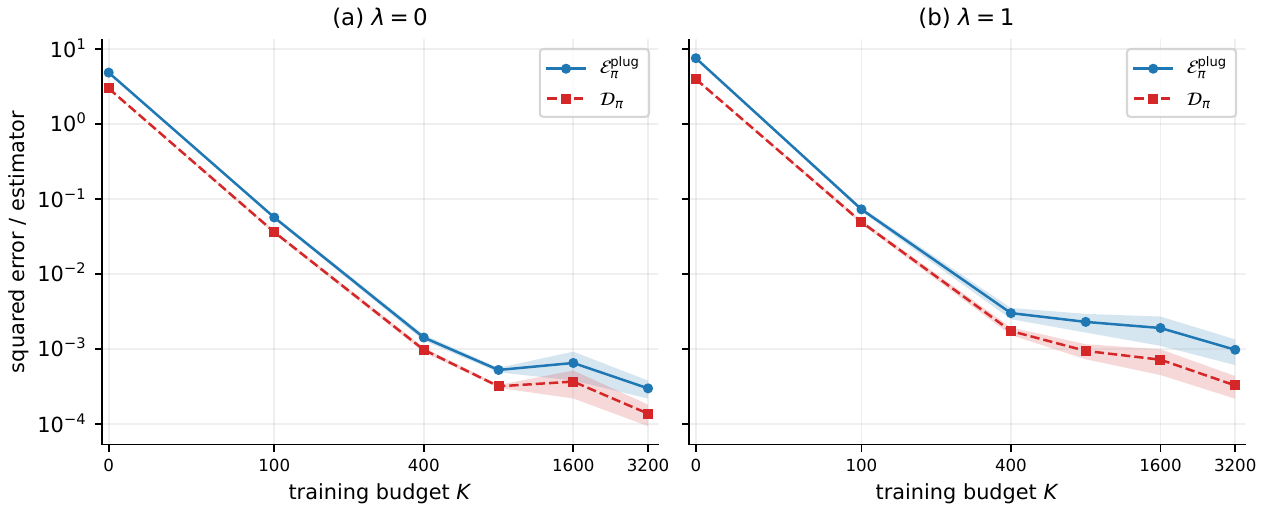}
\caption{Paired nonlinear-mean jump study: $\Epi^{\mathrm{plug}}$ and $\Dpi$
versus training budget, with $\lambda=0$ (a) and $\lambda=1$ (b), on common vertical and horizontal scales.
Curves show seed means with ${}\pm{}s/\sqrt{5}$ bands over five matched seeds.}
\label{B:fig:nonlinear-mean}
\end{figure}

\paragraph{Effect of the nonlinear mean term.}\label{B:ssec:nonlinear-mean-results}
The estimator follows the mean true error during training (Figure~\ref{B:fig:nonlinear-mean}).
With the nonlinear mean term present, the final mean error is higher. Across all saved
checkpoints, the ratio of true error to estimator lies in similar ranges with and without the
term, $\nonlinActiveBandVal$ and $\nonlinBaseBandVal$, respectively.

\begin{figure}[tb]\centering
\incfig[0.6\linewidth]{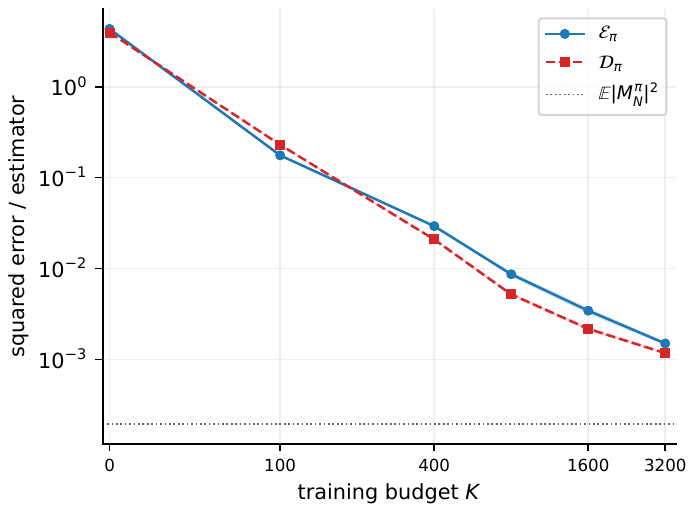}
\caption{State-nonlinear benchmark: $\Epi$ and $\Dpi$ versus training budget, seed means
with ${}\pm{}s/\sqrt{\snlSeedsVal}$ bands, and the seed-mean martingale contribution
$\E|\ex M_N|^2$ (dotted).}
\label{B:fig:snl}
\end{figure}

\paragraph{State-nonlinear benchmark.}\label{B:ssec:snl-res}
On the benchmark of Section~\ref{B:ssec:snl} (Figure~\ref{B:fig:snl}), across $\snlSeedsVal$
seeds the seed-mean $\Epi$ falls from $\snlEStartVal$ to $\snlEEndVal$ and $\Dpi$ from
$\snlDStartVal$ to $\snlDEndVal$ at budget $K=\snlBudgetMaxVal$, with $\Dpi$ still decreasing
from $\snlDPenultVal$ at $K=\snlBudgetPenultVal$. The measured effectivity $\Dpi/\Epi$ lies in
$\snlBandVal$ over the $\snlCellsVal$ checkpoints, with detrended correlation
$\snlDetrendVal$. Thus the observed estimator--error tracking is not confined to the affine
benchmarks. With $\ap M=0$, the reference
martingale contributes an irreducible error whose seed mean is $\snlFloorVal$, approximately
$\snlFloorShareVal$ of the final mean error.

\subsubsection{Picard/LSMC}\label{B:ssec:lsmc-res}
\begin{figure}[tb]\centering
\incfig[0.65\linewidth]{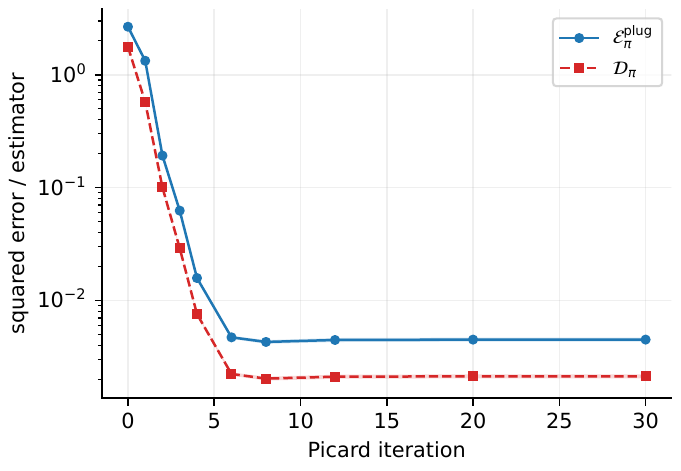}
\caption{Picard/LSMC estimator and error versus iteration. Curves show seed means with
${}\pm{}s/\sqrt{\lsmcNseedsVal}$ bands; each seed uses a fixed training cloud of
$\lsmcNtrainVal$ paths and common random numbers.}
\label{B:fig:lsmc}
\end{figure}

Across $\lsmcNseedsVal$ seeds and up to $\lsmcBudgetMaxVal$ Picard iterations,
$\Epi^{\mathrm{plug}}$ falls rapidly and then plateaus, while $\Dpi$ follows the same trend. The
measured effectivity band is $\lsmcBandVal$ (Figure~\ref{B:fig:lsmc}), with detrended correlation
$\stratLsmcDetrendCorrVal$, and the interior backward residual accounts for a median
$\lsmcCYshareVal$ of the scored estimator.

\subsubsection{Mean-field inventory control}\label{B:ssec:lqinv-res}

\begin{figure}[tb]\centering
\incfig{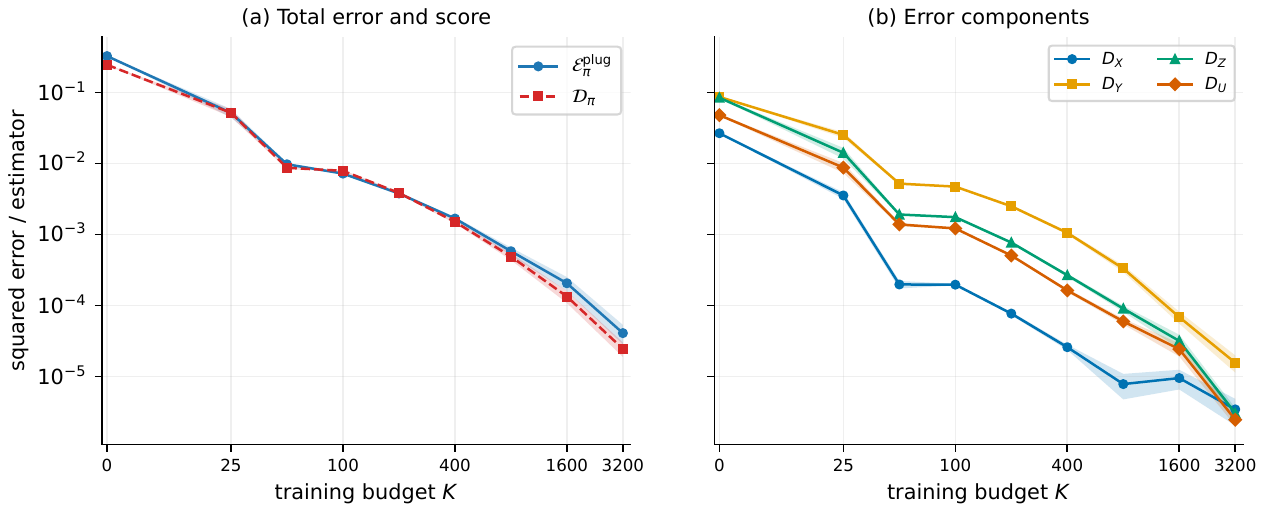}
\caption{Mean-field inventory benchmark. (a) Plug-in estimator and full discrete error
versus training budget. (b) Components $D_X,D_Y,D_Z,D_U$ of \eqref{B:eq:Dpi}, from the same
90 checkpoints; $D_M=0$ since $\ap M=\ex M=0$. The jump term is
$D_U=\sum_{i=0}^{N-1}\tau\nu(E)\,\E|\ap U_i-\ex U_i|^2$. Curves show seed means with
${}\pm{}s/\sqrt{\lqinvNseedsVal}$ bands.}
\label{B:fig:lqinv}
\end{figure}

When the same FBSNN solver is trained on this benchmark (Figure~\ref{B:fig:lqinv}),
$\Epi^{\mathrm{plug}}$ tracks the full discrete tuple error throughout training. Across
$\lqinvNseedsVal$ seeds, the measured effectivity $\Dpi/\Epi^{\mathrm{plug}}$ remains in
$\lqinvBandVal$, with detrended correlation $\stratLqinvDetrendCorrVal$, while $\Dpi$ falls from
$\lqinvDpiStartVal$ to $\lqinvDpiEndVal$. $D_Y$ is the largest seed-mean component, and the
jump-integrand error $D_U$ decreases throughout the reported budgets
(Figure~\ref{B:fig:lqinv}(b)).

\begin{figure}[tb]\centering
\incfig{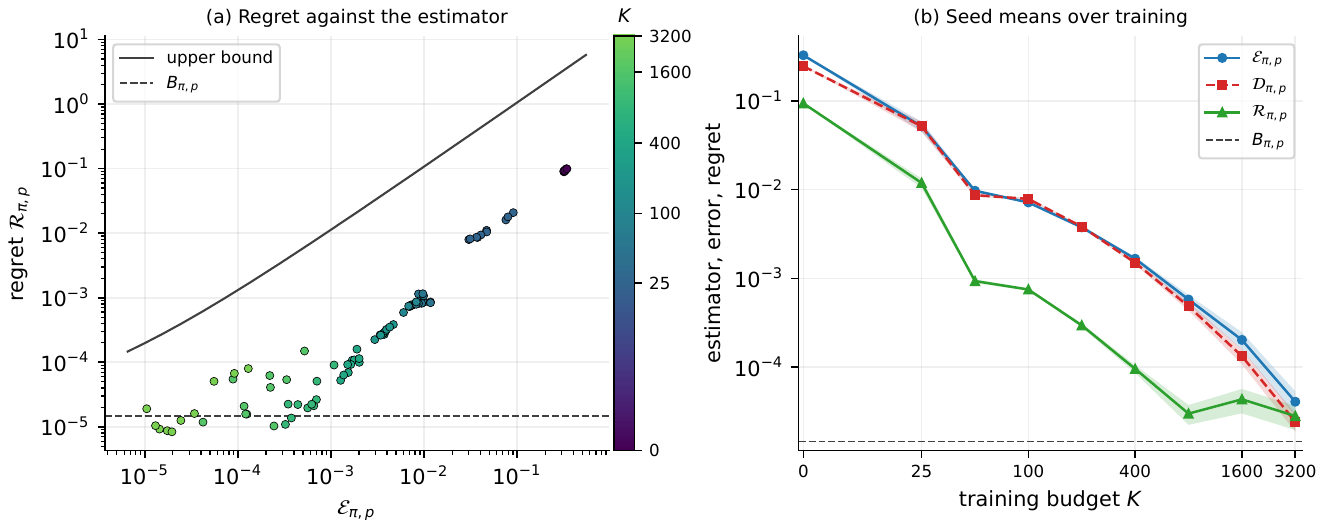}
\caption{Inventory control on the particle system, first cohort. (a) Regret against the
estimator for the $\ictlCellsVal$ checkpoints, coloured by training budget, with $B_{\pi,p}$
(dashed) and the upper bound $(\sqrt{C_{\mathrm{inv}}\Epip}+\sqrt{B_{\pi,p}})^2$ of
Proposition~\ref{B:prop:regret}; its lower bound is zero over the plotted range. (b) Seed means
of $\Epip$, $\Dpip$ and $\mathcal R_{\pi,p}$ against training budget, with
${}\pm{}s/\sqrt{\ictlSeedsVal}$ bands.}
\label{B:fig:lqinvctl}
\end{figure}

\paragraph{Control regret.} Figure~\ref{B:fig:lqinvctl} reports the interacting-particle
comparison of Section~\ref{B:ssec:lqinv}, with estimator $\Epip$ and error $\Dpip$.
On the reported grid the hypotheses of
Proposition~\ref{B:prop:regret} hold, with $C_{\mathrm{inv}}=\ictlCinvVal$ and
$B_{\pi,p}=\ictlBVal$. Over the $\ictlCellsVal$ checkpoints, $\Dpip/\Epip$ lies in
$\ictlBandVal$, with detrended correlation $\ictlDetrendEDVal$. The computed upper bound is
conservative: it is at least $\ictlSlackVal$ times the measured regret in every evaluation
(Figure~\ref{B:fig:lqinvctl}(a)).

Regret does not follow the estimator late in training. From $K=\ictlLateFromVal$ on, the
seed-mean $\Epip$ falls by a factor $\ictlLateEFacVal$ while the seed-mean regret stays between
$\ictlLateRegLoVal$ and $\ictlLateRegHiVal$ (Figure~\ref{B:fig:lqinvctl}(b)), and
$\ictlBelowBVal$ checkpoints have regret below $B_{\pi,p}$, closer to the discrete optimum than
the Euler-FBSDE feedback. Refining $\ictlRefCountVal$ checkpoints to $p=\ictlRefNpartVal$ changes their
regret by at most $\ictlRefShiftVal$.

\subsubsection{Pilot-statistic sensitivity}\label{B:ssec:candlaw}
When the same paths are rescored with $s^C$ in place of $s^A$, the median absolute percentage
change in plug-in effectivity is below $1\%$ in every study (Table~\ref{B:tab:candlaw}). The
largest change is about $10\%$, for the nonlinear-mean model with $\lambda=1$.

\begin{table}[tb]\centering
\caption{Absolute percentage change in effectivity relative to the frozen-pilot value under the
rescoring procedure of Section~\ref{B:sec:deepsolver}.}
\label{B:tab:candlaw}
\begin{tabular}{@{}lrccc@{}}
\toprule
study & cells & median & $95$th pct.\ & max\\
\midrule
\candStatsTableBody
\bottomrule
\end{tabular}
\end{table}

\subsubsection{Checkpoint selection and ranking}\label{B:ssec:ranking}
\paragraph{Solution-error selection.} To test the estimator as a checkpoint-selection criterion,
we select checkpoints by $\Epi^{\mathrm{plug}}$ and, where available, by the held-out consistency
loss $\mathcal L$ on a separate $100{,}000$-path selection batch, then evaluate their solution
error on an independent $\selBTVal$-path test cohort. For each seed, the saved checkpoint with the
smallest measured test error $\Dpi$ is called the test-best checkpoint.

\begin{table}[tb]\centering\small
\caption{Checkpoint selection on independent data. The test-error ratio is the
checkpoint's $\Dpi$ divided by that of the test-best checkpoint for the same seed.}
\label{B:tab:selection}
\begin{tabular}{@{}llcc@{}}
\toprule
study & criterion & test-best selections & worst test-error ratio\\
\midrule
neural & final checkpoint & $\selValFinalHitsVal/\selValFinalSeedsVal$ & $\selValFinalWorstVal$\\
 & $\Epi^{\mathrm{plug}}$ & $\selValEHitsVal/\selValESeedsVal$ & $\selValEWorstVal$\\
 & $\mathcal L$ & $\selValLHitsVal/\selValLSeedsVal$ & $\selValLWorstVal$\\
\addlinespace
Picard/LSMC & final checkpoint & $\selLsmcFinalHitsVal/\selLsmcFinalSeedsVal$ & $\selLsmcFinalWorstVal$\\
 & $\Epi^{\mathrm{plug}}$ & $\selLsmcEHitsVal/\selLsmcESeedsVal$ & $\selLsmcEWorstVal$\\
 & $\mathcal L$ & -- & --\\
\addlinespace
inventory & final checkpoint & $\selLqinvFinalHitsVal/\selLqinvFinalSeedsVal$ & $\selLqinvFinalWorstVal$\\
 & $\Epi^{\mathrm{plug}}$ & $\selLqinvEHitsVal/\selLqinvESeedsVal$ & $\selLqinvEWorstVal$\\
 & $\mathcal L$ & $\selLqinvLHitsVal/\selLqinvLSeedsVal$ & $\selLqinvLWorstVal$\\
\bottomrule
\end{tabular}
\end{table}

Across all three studies, the squared error of the estimator-selected checkpoint was within
$\rv{select.eWorstPct}\%$ of the test-best error in every seed (Table~\ref{B:tab:selection}). In
the linear neural study, the estimator selected the test-best checkpoint in all ten seeds. The
held-out loss did so in nine, with a worst test-error ratio of $\selValLWorstVal$.
Estimator-based selection also improved on the final-checkpoint baseline for Picard/LSMC, while
the final checkpoint was test-best in every inventory seed.

For each training budget, we compute the Spearman correlation across seeds between
$\Epi^{\mathrm{plug}}$ and $\Dpi$. The median correlations across budgets are
$\stratValSpearMedVal$, $\stratLsmcSpearMedVal$ and $\stratLqinvSpearMedVal$ for the linear
neural, Picard/LSMC and inventory studies, respectively. In the inventory study at $K=50$,
smaller estimator values do not consistently identify runs with smaller errors (Spearman
correlation $\stratLqinvSpearLoVal$).

Within each seed, we also compare how the two criteria rank all saved checkpoints against their
ranking by $\Dpi$. The median Spearman correlation \hoEstCorrPhrase. For the held-out loss, the lowest
correlation \hoObjLoPhrase. Thus
the loss can rank the checkpoints closely by error and still select a poor checkpoint, as in the
linear study (Table~\ref{B:tab:selection}).

\paragraph{Control-regret selection.} For control-regret selection, the estimator and consistency
loss select the same checkpoint in $\ictlSelAgreeVal$ of $\ictlSeedsVal$ seeds, but neither
criterion reliably identifies low-regret checkpoints: both agree with direct regret minimisation in
$\ictlSelEHitsVal$ of $\ictlSeedsVal$ seeds. On the independent test cohort, the largest regret was
$\ictlSelMaxVal$ times that of the checkpoint selected by direct regret minimisation. The mean paired
test-regret difference between estimator and loss selection is $\ictlSelDiffVal$ (estimator minus
loss), with standard error $\ictlSelDiffSeVal$, indicating no consistent advantage for either
criterion. Checkpoint selection for solution error and control regret are therefore distinct
objectives.

\FloatBarrier

\appendix

\section{Projection, discrete stability, and well-posedness}

\subsection{One-step Kunita--Watanabe projection}\label{B:app:kw}
\begin{proof}[Proof of Lemma~\ref{B:lem:kw}]
Fix $i$ (on the uniform grid, $\tau=t_{i+1}-t_i$), and let
$\mathcal L_i:=\{R\in L^2(\Omega,\F_{t_{i+1}};\R^m):\Ei[R]=0\}$, a closed subspace of
$L^2(\Omega,\F_{t_{i+1}};\R^m)$ with inner product $\inner{A}{B}_{L^2}:=\E\inner{A}{B}_{\R^m}$. Put
\[
\mathcal W_i:=\{Z\Delta W_i:Z\in L^2(\Omega,\F_{t_i};\R^{m\times d})\},\quad
\mathcal J_i:=\Big\{\int_E U\,\Delta\Nt_i:U\in L^2(\Omega,\F_{t_i};L^2(\nu;\R^m))\Big\}.
\]
Both lie in $\mathcal L_i$ because the Brownian and compensated-Poisson increments are
conditionally centred, and their conditional isometries give
\[
\Ei|Z\Delta W_i|^2=\tau|Z|_F^2,
\qquad
\Ei\Bigl|\int_E U\,\Delta\Nt_i\Bigr|^2=\tau\|U\|_{\Lnu}^2,
\]
Taking expectations shows that both maps are scaled $L^2$-isometries with factor
$\sqrt\tau>0$; their complete domains therefore have closed ranges
$\mathcal W_i,\mathcal J_i$.

The conditional Brownian--Poisson orthogonality gives
\[
\Ei\Bigl[(Z\Delta W_i)\Bigl(\int_E U\,\Delta\Nt_i\Bigr)^{\!\top}\Bigr]=0,
\]
hence $\mathcal W_i\perp\mathcal J_i$, and $\mathcal W_i\oplus\mathcal J_i$ is closed.

\emph{Decomposition.} For $R\in\mathcal L_i$, orthogonal projection onto $\mathcal W_i\oplus\mathcal J_i$
gives $R=Z\Delta W_i+\int_E U\,\Delta\Nt_i+\Delta M$ with $\Delta M\perp\mathcal W_i\oplus\mathcal J_i$.
Since $\mathcal W_i\perp\mathcal J_i$, the two projected parts are separately determined, and the
conditional isometries make $Z,U$ unique in their $L^2$-spaces.

\emph{Conditional orthogonality.} Let $A\in L^\infty(\Omega,\F_{t_i};\R^{m\times d})$ and
$C\in L^\infty(\Omega,\F_{t_i};\R^{m\times m})$ be bounded random test matrices, and let
$a\in L^2(\Omega,\F_{t_i};L^2(\nu;\R^m))$. Then
$A\Delta W_i\in\mathcal W_i$, while $Ca$ is again a square-integrable $\F_{t_i}$-measurable
$\Lnu$-integrand with $\int_E Ca\,\Delta\Nt_i=C\int_E a\,\Delta\Nt_i\in\mathcal J_i$, so
\[
\begin{aligned}
0&=\E\inner{\Delta M}{A\Delta W_i}
=\E\,\tr\bigl(A^{\!\top}\Ei[\Delta M\Delta W_i^{\!\top}]\bigr),\\[2pt]
0&=\E\Bigl\langle\Delta M,\;C\!\int_E a\,\Delta\Nt_i\Bigr\rangle
=\E\,\tr\Bigl(C^{\!\top}\Ei\bigl[\Delta M\bigl(\textstyle\int_E a\,\Delta\Nt_i\bigr)^{\!\top}\bigr]\Bigr).
\end{aligned}
\]
As $A$ and $C$ are arbitrary,
\[
\Ei[\Delta M\,\Delta W_i^{\!\top}]=0,
\qquad
\Ei\Bigl[
\Delta M\Bigl(\int_E a\,\Delta\Nt_i\Bigr)^{\!\top}
\Bigr]=0
\quad\text{for every such }a.
\]
Since $\Delta M\in\mathcal L_i$, we also have $\Ei[\Delta M]=0$, so $\Delta M$ is admissible.

By \eqref{B:eq:kw-decomp},
\[
\begin{aligned}
\Ei[R\,\Delta W_i^{\!\top}]
&=Z\,\Ei[\Delta W_i\Delta W_i^{\!\top}]
+\Ei\Bigl[
\Bigl(\int_E U\,\Delta\Nt_i\Bigr)\Delta W_i^{\!\top}
\Bigr]
+\Ei[\Delta M\,\Delta W_i^{\!\top}]\\
&=Z\,\Ei[\Delta W_i\Delta W_i^{\!\top}]
=\tau Z.
\end{aligned}
\]
Thus $Z=\tau^{-1}\Ei[R\,\Delta W_i^{\!\top}]$.

Similarly, for every such $a$, orthogonality and the conditional isometry give
\[
\begin{aligned}
\Ei\Bigl[
R\Bigl(\int_E a\,\Delta\Nt_i\Bigr)^{\!\top}
\Bigr]
&=\Ei\Bigl[
\Bigl(\int_E U\,\Delta\Nt_i\Bigr)
\Bigl(\int_E a\,\Delta\Nt_i\Bigr)^{\!\top}
\Bigr]\\
&=\tau\int_E U(e)a(e)^{\!\top}\,\nu(\diff e).
\end{aligned}
\]
Taking the trace gives
\[
\Ei\Bigl[
\Bigl\langle R,\int_E a\,\Delta\Nt_i\Bigr\rangle
\Bigr]
=\tau\int_E\inner{U(e)}{a(e)}\,\nu(\diff e).
\qedhere
\]
\end{proof}

\subsection{Proof of Proposition~\ref{B:prop:pert}}
\label{B:app:stab}
With the exact-minus-perturbed differences $\delta(\cdot)_i:=(\cdot)^\pi_i-(\cdot)'_i$ of
Section~\ref{B:sec:stab}, set
\[
\Sd:=\sum_{i=0}^{N-1}\tau\,\E\big[|\delta X_i|^2+|\delta Y_i|^2+|\delta Z_i|_F^2+\|\delta U_i\|_{\Lnu}^2\big].
\]
Introduce the $G$-weighted energy
\begin{equation}
\label{B:eq:EGdef}
\EG:=\alpha\,\E|G\delta X_N|^2+\beta_2\sum_{i=0}^{N-1}\tau\,\E|G\delta X_i|^2
+\beta_1\sum_{i=0}^{N-1}\tau\,\E\big[|\Gstar\delta Y_i|^2+|\Gstar\delta Z_i|_F^2+\|\Gstar\delta U_i\|_{\Lnu}^2\big]\ge0.
\end{equation}
Subtracting \eqref{B:eq:pfwd}--\eqref{B:eq:pbwd} from \eqref{B:eq:dfwd}--\eqref{B:eq:dbwd}
and grouping the $\F_{t_i}$-measurable coefficients gives
\begin{equation}
\label{B:eq:grouped}
\begin{aligned}
\Delta\delta X_i&=(\delta b_i-\chi_i)\tau+(\delta\sigma_i-\gamma_i)\Delta W_i+\int_E(\delta h_i-\eta_i)\,\Delta\Nt_i,\\
\Delta\delta Y_i&=-(\delta f_i+\phi_i)\tau+\delta Z_i\Delta W_i+\int_E\delta U_i\,\Delta\Nt_i+\Delta\delta M_i,
\end{aligned}
\end{equation}
with $\delta Y_N=[g(\ex X_N,\Law(\ex X_N))-g(X'_N,\Law(X'_N))]-\zeta$.

\emph{One-step expansion.} By the discrete product identity, \eqref{B:eq:grouped}, the
conditional centring and isometries of $\Delta W_i$ and $\Delta\Nt_i$, and the orthogonality of
$\Delta\delta M_i$ to both (Definition~\ref{B:def:incr-adm}),
\[
\begin{aligned}
\Ei\Delta\inner{G\delta X_i}{\delta Y_i}
&=\Ei\inner{G\Delta\delta X_i}{\delta Y_i}+\Ei\inner{G\delta X_i}{\Delta\delta Y_i}
+\Ei\inner{G\Delta\delta X_i}{\Delta\delta Y_i}\\
&=\tau\inner{\delta b_i-\chi_i}{\Gstar\delta Y_i}-\tau\inner{\delta f_i+\phi_i}{G\delta X_i}
+\Ei\inner{G\Delta\delta X_i}{\Delta\delta Y_i}\\
&=\tau\inner{\delta b_i-\chi_i}{\Gstar\delta Y_i}-\tau\inner{\delta f_i+\phi_i}{G\delta X_i}
+\tau\inner{\delta\sigma_i-\gamma_i}{\Gstar\delta Z_i}_F\\
&\quad+\tau\!\int_E\inner{\delta h_i-\eta_i}{\Gstar\delta U_i}\,\nu(\diff e)
-\tau^2\inner{G(\delta b_i-\chi_i)}{\delta f_i+\phi_i} .
\end{aligned}
\]
Writing
\[
\begin{aligned}
\Bc_i&:=\inner{\delta b_i}{\Gstar\delta Y_i}+\inner{\delta\sigma_i}{\Gstar\delta Z_i}_F
-\inner{\delta f_i}{G\delta X_i}+\int_E\inner{\delta h_i}{\Gstar\delta U_i}\,\nu(\diff e),\\
\Pi_i&:=\inner{\chi_i}{\Gstar\delta Y_i}+\inner{\gamma_i}{\Gstar\delta Z_i}_F
+\inner{\phi_i}{G\delta X_i}+\int_E\inner{\eta_i}{\Gstar\delta U_i}\,\nu(\diff e),\\
R_i&:=-\tau^2\inner{G(\delta b_i-\chi_i)}{\delta f_i+\phi_i},
\end{aligned}
\]
\begin{equation}
\label{B:eq:onestep}
\Ei\big[\Delta\inner{G\delta X_i}{\delta Y_i}\big]=\tau\,\Bc_i-\tau\,\Pi_i+R_i.
\end{equation}

\emph{Summation and monotonicity.} Taking full expectations and summing \eqref{B:eq:onestep},
then applying Assumption~\ref{B:ass:mono} at the terminal node and on the interior,
\[
\begin{aligned}
\E\inner{G\delta X_N}{\delta Y_N}-\E\inner{G\delta X_0}{\delta Y_0}
&=\sum_{i=0}^{N-1}\tau\E\Bc_i-\sum_{i=0}^{N-1}\tau\E\Pi_i+\sum_{i=0}^{N-1}\E R_i,\\
\E\inner{G\delta X_N}{\delta Y_N}&\ge\alpha\E|G\delta X_N|^2-\E\inner{G\delta X_N}{\zeta},\\
\sum_{i=0}^{N-1}\tau\E\Bc_i&\le-\beta_2\sum_{i=0}^{N-1}\tau\E|G\delta X_i|^2\\
&\qquad-\beta_1\sum_{i=0}^{N-1}\tau\E\big[|\Gstar\delta Y_i|^2+|\Gstar\delta Z_i|_F^2
+\|\Gstar\delta U_i\|_{\Lnu}^2\big].
\end{aligned}
\]
Moving the negative interior sum to the left and bounding $-\sum_{i=0}^{N-1}\tau\E\Pi_i$ and the remainder
by absolute values,
\begin{equation}
\label{B:eq:raw}
\EG\le\E\inner{G\delta X_0}{\delta Y_0}+\E\inner{G\delta X_N}{\zeta}+\sum_{i=0}^{N-1}\tau\,|\E\Pi_i|+\sum_{i=0}^{N-1}\E|R_i|.
\end{equation}

\emph{Young estimates.} By Assumption~\ref{B:ass:lip}, the law coupling of
Lemma~\ref{B:lem:law}, and Young's inequality on each pairing of a perturbation with a component
of $\delta$, with $\Pp$ the interior perturbation sum in \eqref{B:eq:Pnorm},
\[
\begin{aligned}
\sum_{i=0}^{N-1}\tau\,|\E\Pi_i|&\le\varepsilon\,\Sd+K_\varepsilon\,\Pp,&\qquad
\sum_{i=0}^{N-1}\E|R_i|&\le C\tau(\Sd+\Pp),\\
|\E\inner{G\delta X_0}{\delta Y_0}|&\le\varepsilon\E|\delta Y_0|^2+K_\varepsilon\E|\delta X_0|^2,
&\qquad
|\E\inner{G\delta X_N}{\zeta}|&\le\varepsilon\E|\delta X_N|^2+K_\varepsilon\E|\zeta|^2,
\end{aligned}
\]
with $K_\varepsilon$ a Young constant. Inserting these into \eqref{B:eq:raw} gives
\begin{equation}
\label{B:eq:rawabs}
\EG\le(C\tau+C_0\varepsilon)\,\Sd+\varepsilon\big(\E|\delta Y_0|^2+\E|\delta X_N|^2\big)+K_\varepsilon\big(\E|\delta X_0|^2+\|s\|_{\mathrm{src}}^2\big).
\end{equation}

\emph{Forward estimate.} Two energy estimates bound $\Sd$, $\E|\delta Y_0|^2$ and
$\E|\delta X_N|^2$, still on the right of \eqref{B:eq:rawabs}. By \eqref{B:eq:grouped}, the
conditional isometries of the driving increments, Assumption~\ref{B:ass:lip} and
Lemma~\ref{B:lem:law},
\[
\begin{aligned}
\E|\delta X_{i+1}|^2
&=\E|\delta X_i|^2+2\E\inner{\delta X_i}{\Delta\delta X_i}+\E|\Delta\delta X_i|^2\\
&=\E|\delta X_i|^2+2\tau\E\inner{\delta X_i}{\delta b_i-\chi_i}+\tau^2\E|\delta b_i-\chi_i|^2\\
&\qquad+\tau\E|\delta\sigma_i-\gamma_i|_F^2+\tau\E\|\delta h_i-\eta_i\|_{\Lnu}^2\\
&\le(1+C\tau)\E|\delta X_i|^2
+C\tau\E\big[|\delta Y_i|^2+|\delta Z_i|_F^2+\|\delta U_i\|_{\Lnu}^2\big]\\
&\qquad+C\tau\E\big[|\chi_i|^2+|\gamma_i|_F^2+\|\eta_i\|_{\Lnu}^2\big].
\end{aligned}
\]
With accumulated factor $(1+C\tau)^N\le e^{CT}$, independent of the grid, discrete
Gr\"onwall gives
\begin{equation}
\label{B:eq:ffwd}
\begin{split}
\max_{0\le i\le N}\E|\delta X_i|^2\le C\Big(&\E|\delta X_0|^2
+\sum_{i=0}^{N-1}\tau\E\big[|\delta Y_i|^2+|\delta Z_i|_F^2+\|\delta U_i\|_{\Lnu}^2\big]\\
&+\sum_{i=0}^{N-1}\tau\E\big[|\chi_i|^2+|\gamma_i|_F^2+\|\eta_i\|_{\Lnu}^2\big]\Big).
\end{split}
\end{equation}

\emph{Backward estimate.} Squaring the second equation of \eqref{B:eq:grouped} and using the
orthogonality of $\delta Z_i\Delta W_i$, $\int_E\delta U_i\Delta\Nt_i$ and $\Delta\delta M_i$,
\[
\begin{aligned}
\E|\delta Y_i|^2&+\tau\E|\delta Z_i|_F^2+\tau\E\|\delta U_i\|_{\Lnu}^2+\E|\Delta\delta M_i|^2\\
&=\E|\delta Y_{i+1}|^2+2\tau\E\inner{\delta Y_i}{\delta f_i+\phi_i}-\tau^2\E|\delta f_i+\phi_i|^2\\
&\le\E|\delta Y_{i+1}|^2+2\tau\,\E\big[|\delta Y_i|\,|\delta f_i+\phi_i|\big].
\end{aligned}
\]
By Assumption~\ref{B:ass:lip} and the synchronous law coupling of
Lemma~\ref{B:lem:law}, put
\[
a_i:=\bigl(\E|\delta Y_i|^2\bigr)^{1/2},\qquad
b_i:=\bigl(\E|\delta X_i|^2\bigr)^{1/2},\qquad
q_i:=\Bigl(\E\bigl[|\delta Z_i|_F^2+\|\delta U_i\|_{\Lnu}^2\bigr]\Bigr)^{1/2}.
\]
Minkowski's inequality, Assumption~\ref{B:ass:lip}, and
$W_2(\mu_i,\mu'_i)\le(a_i^2+b_i^2+q_i^2)^{1/2}$ imply
\[
\|\delta f_i\|_{L^2(\Omega)}
\le L\bigl(a_i+b_i+\|\delta Z_i\|_{L^2(\Omega)}
 +\|\delta U_i\|_{L^2(\Omega;\Lnu)}+W_2(\mu_i,\mu'_i)\bigr)
\le L(2a_i+2b_i+3q_i).
\]
Fix $\theta:=\tfrac12$ and set $C_\theta:=1+6L+18L^2$.
By Cauchy--Schwarz and Young's inequality \emph{after taking expectations},
\[
\begin{aligned}
2\E\bigl[|\delta Y_i|\,|\delta f_i+\phi_i|\bigr]
&\le 2a_i\bigl(L(2a_i+2b_i+3q_i)
          +(\E|\phi_i|^2)^{1/2}\bigr)\\
&\le C_\theta(a_i^2+b_i^2)+\theta q_i^2+2\E|\phi_i|^2.
\end{aligned}
\]
Indeed, the mixed products are bounded by
$4La_ib_i\le2La_i^2+2Lb_i^2$,
$6La_iq_i\le18L^2a_i^2+\tfrac12q_i^2$, and
$2a_i(\E|\phi_i|^2)^{1/2}\le a_i^2+\E|\phi_i|^2$.
Substituting into the preceding backward energy identity gives
\[
\begin{aligned}
&(1-C_\theta\tau)\E|\delta Y_i|^2
+(1-\theta)\tau\bigl(\E|\delta Z_i|_F^2
                 +\E\|\delta U_i\|_{\Lnu}^2\bigr)
+\E|\Delta\delta M_i|^2\\
&\qquad\le\E|\delta Y_{i+1}|^2+C_\theta\tau\E|\delta X_i|^2
+2\tau\E|\phi_i|^2.
\end{aligned}
\]
For $C_\theta\tau\le\tfrac12$, divide by $1-C_\theta\tau$ and use
$(1-C_\theta\tau)^{-1}\le1+2C_\theta\tau$. Dropping the resulting factors
larger than one on the nonnegative left-hand terms,
\begin{equation}
\label{B:eq:bwdrec}
\begin{aligned}
&\E|\delta Y_i|^2+(1-\theta)\tau
 \bigl(\E|\delta Z_i|_F^2+\E\|\delta U_i\|_{\Lnu}^2\bigr)
 +\E|\Delta\delta M_i|^2\\
&\qquad\le(1+2C_\theta\tau)\E|\delta Y_{i+1}|^2
+2C_\theta\tau\E|\delta X_i|^2+4\tau\E|\phi_i|^2.
\end{aligned}
\end{equation}
Since the $Z$-, $U$- and $M$-terms on the left are nonnegative, backward Gr\"onwall, using
\[
(1+2C_\theta\tau)^N\le e^{2C_\theta T},
\]
gives
\begin{equation}
\label{B:eq:bwdY}
\max_{0\le i\le N}\E|\delta Y_i|^2
\le C\Big(
\E|\delta Y_N|^2
+\sum_{i=0}^{N-1}\tau\,\E|\delta X_i|^2
+\sum_{i=0}^{N-1}\tau\,\E|\phi_i|^2
\Big).
\end{equation}
Summing over $i$ gives
\[
\begin{aligned}
&(1-\theta)\sum_{i=0}^{N-1}\tau
\bigl(
\E|\delta Z_i|_F^2
+\E\|\delta U_i\|_{\Lnu}^2
\bigr)
+\sum_{i=0}^{N-1}\E|\Delta\delta M_i|^2\\
&\qquad\le
\E|\delta Y_N|^2
+2C_\theta\sum_{i=0}^{N-1}\tau\,\E|\delta Y_{i+1}|^2
+2C_\theta\sum_{i=0}^{N-1}\tau\,\E|\delta X_i|^2
+4\sum_{i=0}^{N-1}\tau\,\E|\phi_i|^2.
\end{aligned}
\]
By orthogonality of the martingale increments and $\delta M_0=0$,
\[
\sum_{i=0}^{N-1}\E|\Delta\delta M_i|^2
=\E|\delta M_N|^2.
\]
Using \eqref{B:eq:bwdY} and $1-\theta=\tfrac12$ gives
\begin{equation}
\label{B:eq:bwdZU}
\begin{aligned}
&\sum_{i=0}^{N-1}\tau\,\E|\delta Z_i|_F^2
+\sum_{i=0}^{N-1}\tau\,\E\|\delta U_i\|_{\Lnu}^2
+\E|\delta M_N|^2\\
&\qquad\le
C\Big(
\E|\delta Y_N|^2
+\sum_{i=0}^{N-1}\tau\,\E|\delta X_i|^2
+\sum_{i=0}^{N-1}\tau\,\E|\phi_i|^2
\Big).
\end{aligned}
\end{equation}
Combining \eqref{B:eq:bwdY} and \eqref{B:eq:bwdZU} yields
\begin{equation}
\label{B:eq:fbwd}
\begin{aligned}
\max_{0\le i\le N}\E|\delta Y_i|^2
&+\sum_{i=0}^{N-1}\tau\,\E|\delta Z_i|_F^2
+\sum_{i=0}^{N-1}\tau\,\E\|\delta U_i\|_{\Lnu}^2
+\E|\delta M_N|^2\\
&\le C\Big(
\E|\delta Y_N|^2
+\sum_{i=0}^{N-1}\tau\,\E|\delta X_i|^2
+\sum_{i=0}^{N-1}\tau\,\E|\phi_i|^2
\Big).
\end{aligned}
\end{equation}
By the terminal Lipschitz bound,
\[
\begin{aligned}
\max_{0\le i\le N}\E|\delta Y_i|^2
&+\sum_{i=0}^{N-1}\tau\,\E|\delta Z_i|_F^2
+\sum_{i=0}^{N-1}\tau\,\E\|\delta U_i\|_{\Lnu}^2
+\E|\delta M_N|^2\\
&\le C\Big(
\E|\delta X_N|^2
+\E|\zeta|^2
+\sum_{i=0}^{N-1}\tau\,\E|\delta X_i|^2
+\sum_{i=0}^{N-1}\tau\,\E|\phi_i|^2
\Big).
\end{aligned}
\]

\emph{Rank closure.} By Assumption~\ref{B:ass:G}, $\Gstar$ is injective when $n>m$ and $G$ is
injective when $m>n$, and Assumption~\ref{B:ass:mono} then supplies $\beta_1>0$ when $n>m$ and
$\alpha,\beta_2>0$ when $m>n$. We write Case~A for $n\ge m$ with $\beta_1>0$ ($\Gstar$
injective) and Case~B for $m\ge n$ with $\alpha,\beta_2>0$ ($G$ injective).
When $n=m$, either $\beta_1>0$, or the conditions $\alpha+\beta_1>0$ and
$\beta_1+\beta_2>0$ imply $\alpha,\beta_2>0$; thus at least one case applies.

We set
\[
\mu_\star:=\begin{cases}
\beta_1,&n>m,\\[2pt]
\min\{\alpha,\beta_2\},&m>n,\\[2pt]
\max\{\beta_1,\min\{\alpha,\beta_2\}\},&n=m.
\end{cases}
\]

\emph{Case A} ($n\ge m$, $\beta_1>0$, $\Gstar$ injective). Then $|v|^2\le c_G|\Gstar v|^2$, and
by \eqref{B:eq:EGdef}, \eqref{B:eq:ffwd} and \eqref{B:eq:fbwd} with the terminal bound,
\[
\begin{aligned}
\sum_{i=0}^{N-1}\tau\E\big[|\delta Y_i|^2+|\delta Z_i|_F^2+\|\delta U_i\|_{\Lnu}^2\big]
&\le c_G\beta_1^{-1}\EG,\\
\max_{0\le i\le N}\E|\delta X_i|^2&\le C\big(\EG+\E|\delta X_0|^2+\Pp\big),\\
\max_{0\le i\le N}\E|\delta Y_i|^2+\sum_{i=0}^{N-1}\tau\E|\delta Z_i|_F^2+\sum_{i=0}^{N-1}\tau\E\|\delta U_i\|_{\Lnu}^2
+\E|\delta M_N|^2&\le C\big(\EG+\E|\delta X_0|^2+\|s\|_{\mathrm{src}}^2\big).
\end{aligned}
\]

\emph{Case B} ($m\ge n$, $\alpha,\beta_2>0$, $G$ injective). Then $|x|^2\le c_G|Gx|^2$, and by
the same three inputs,
\[
\begin{aligned}
\E|\delta X_N|^2\le c_G\alpha^{-1}\EG,&\qquad
\sum_{i=0}^{N-1}\tau\E|\delta X_i|^2\le c_G\beta_2^{-1}\EG,\\
\max_{0\le i\le N}\E|\delta Y_i|^2+\sum_{i=0}^{N-1}\tau\E|\delta Z_i|_F^2+\sum_{i=0}^{N-1}\tau\E\|\delta U_i\|_{\Lnu}^2
+\E|\delta M_N|^2&\le C\big(\EG+\|s\|_{\mathrm{src}}^2\big),\\
\max_{0\le i\le N}\E|\delta X_i|^2&\le C\big(\EG+\E|\delta X_0|^2+\|s\|_{\mathrm{src}}^2\big).
\end{aligned}
\]
In either case
\begin{equation}
\label{B:eq:caseclose}
\Sd+\E|\delta Y_0|^2+\E|\delta X_N|^2\le C\big(\EG+\E|\delta X_0|^2+\|s\|_{\mathrm{src}}^2\big).
\end{equation}

\emph{Absorption.} Inserting \eqref{B:eq:caseclose} into \eqref{B:eq:rawabs} and absorbing
fixed constants into $C_{\mathrm{abs}},C_0$,
\[
\EG\le(C_{\mathrm{abs}}\tau+C_0\varepsilon)\,\EG
+K_\varepsilon\big(\E|\delta X_0|^2+\|s\|_{\mathrm{src}}^2\big),
\qquad
C_{\mathrm{abs}}=C_{\mathrm{abs}}\bigl(L,T,\|G\|,c_G,\mu_\star^{-1}\bigr).
\]
Choose $\varepsilon$ with $C_0\varepsilon\le\tfrac14$ and set
\[
\tau_{\mathrm{stab}}:=\min\Bigl\{\frac{1}{4C_{\mathrm{abs}}},\ \frac{1}{2C_\theta}\Bigr\},
\]
where the second entry enforces the condition already used for
\eqref{B:eq:bwdrec}. For $\tau<\tau_{\mathrm{stab}}$ the bracket is at most $\tfrac12$, so
\[
\EG\le C_{\mathrm{stab}}\big(\E|\delta X_0|^2+\|s\|_{\mathrm{src}}^2\big),
\qquad
C_{\mathrm{stab}}=C_{\mathrm{stab}}\bigl(L,T,\|G\|,c_G,\mu_\star^{-1}\bigr),
\]
and \eqref{B:eq:caseclose}, \eqref{B:eq:ffwd} and \eqref{B:eq:fbwd} then bound the max-norms,
the $Z,U$-sums and $\E|\delta M_N|^2$ by the same right-hand side, which is
\eqref{B:eq:pert}. \qed

\subsection{Existence by continuation}
\label{B:app:wp}
We prove the existence assertion of Theorem~\ref{B:thm:wp} by the continuation
method~\cite{PengWu,RSZ}.
Throughout, $\mathcal S_\pi$ is the discrete tuple space of Section~\ref{B:sec:scheme}, equipped with the stability norm
\eqref{B:eq:Spi-norm}, exactly the norm controlled by Proposition~\ref{B:prop:pert}. For $\lambda\in[0,1]$ and
$\theta=(x,y,z,u)\in\mathcal H$, define
\begin{equation}\label{B:eq:homotopy}
\begin{gathered}
\psi^\lambda:=(1-\lambda)\psi^0+\lambda\psi,
\qquad \psi\in\{b,\sigma,h,f,g\},\\[4pt]
(b^0,\sigma^0,h^0)(\theta)
:=-\beta_1\bigl(\Gstar y,\Gstar z,\Gstar u\bigr),\\[4pt]
f^0(\theta):=\beta_2Gx,
\qquad g^0(x):=Gx.
\end{gathered}
\end{equation}
Let $E_\lambda$ denote the scheme \eqref{B:eq:dfwd}--\eqref{B:eq:dbwd} with each coefficient
$\psi$ replaced by $\psi^\lambda$. Thus $E_0$ is the linear monotone base problem and $E_1$ is
the original Euler scheme.

\emph{(i) A priori stability, uniform in $\lambda$.} Write
$\mathcal G^\lambda:=(1-\lambda)\,\mathcal G_{\mathrm{lin}}+\lambda\,\mathcal G_{\mathrm{sch}}$
for the interpolated generator, with $\mathcal G_{\mathrm{lin}}$ the linear monotone generator
in \eqref{B:eq:homotopy} and $\mathcal G_{\mathrm{sch}}=(b,\sigma,h,f,g)$ the original generator. By the triangle inequality, the interpolated coefficients satisfy Assumption~\ref{B:ass:lip}
uniformly in $\lambda$ with
\begin{equation}\label{B:eq:hom-lip-constant}
L':=\max\{L,\beta_1\|G\|,\beta_2\|G\|,\|G\|\}.
\end{equation}
Let $\mathcal M(\mathcal G)$ denote the left-hand side of the interior monotonicity inequality in
Assumption~\ref{B:ass:mono}. By linearity in the coefficient differences and
Assumption~\ref{B:ass:mono} at both endpoints,
\[
\begin{aligned}
\mathcal M(\mathcal G^\lambda)
&=(1-\lambda)\mathcal M(\mathcal G_{\mathrm{lin}})
+\lambda\mathcal M(\mathcal G_{\mathrm{sch}})\\
&\le -\beta_1\E\bigl[
|\Gstar\delta Y|^2+|\Gstar\delta Z|_F^2
+\|\Gstar\delta U\|_{\Lnu}^2
\bigr]
-\beta_2\E|G\delta X|^2.
\end{aligned}
\]
The terminal datum is
coercive with constant $(1-\lambda)+\lambda\alpha\ge\min\{1,\alpha\}=:\alpha_*$
\cite[App.~A]{RSZ}; Case~B uses this terminal coercivity and has $\alpha>0$, while in the
square case $\alpha=0$, $\beta_1>0$, Case~A supplies the interior coercivity. Hence
Proposition~\ref{B:prop:pert} holds for every $E_\lambda$, uniformly in $\lambda$, with $L$
replaced by $L'$ and $\alpha$ by $\alpha_*$; write $\tau^\star_{\mathrm{stab}}$ for the
resulting threshold. In particular each $E_\lambda$ has at most one solution in
$\mathcal S_\pi$.

\emph{(ii) Existence for $E_0$ \textup{(}the linear base case, \cite[Lemma~2.6]{RSZ}\textup{)}.}
The linear system is solved by the affine ansatz $Y_i=P_iX_i+p_i$, which yields the
\emph{semi-implicit} backward discrete matrix Riccati recursion
\[
P_N=G,\qquad P_i=(I+\beta_1\tau P_{i+1}\Gstar)^{-1}(P_{i+1}+\beta_2\tau G) .
\]
The recursion diagonalises in the singular value decomposition $G=V\Sigma Q^{\!\top}$, with
$V\in\R^{m\times m}$ and $Q\in\R^{n\times n}$ orthogonal and $\Sigma\in\R^{m\times n}$
rectangular-diagonal. If $P_{i+1}=VD_{i+1}Q^{\!\top}$ with $D_{i+1}$ rectangular-diagonal --- true at
the terminal node, where $D_N=\Sigma$ --- then
\[
P_{i+1}\Gstar=VD_{i+1}\Sigma^{\!\top}V^{\!\top},\qquad
P_{i+1}+\beta_2\tau G=V(D_{i+1}+\beta_2\tau\Sigma)Q^{\!\top},
\]
so $P_i=VD_iQ^{\!\top}$ with one decoupled scalar recursion per singular value $\sigma_k>0$, $1\le k\le\min(m,n)$,
\[
d^{(k)}_i=\frac{d^{(k)}_{i+1}+\beta_2\tau\sigma_k}{1+\beta_1\tau\sigma_k\,d^{(k)}_{i+1}},
\qquad d^{(k)}_N=\sigma_k .
\]
When $m>n$, the remaining $m-n$ diagonal entries of $D_{i+1}\Sigma^{\!\top}$ are zero and give
unit eigenvalues of $I+\beta_1\tau P_{i+1}\Gstar$.

Backward induction preserves $d^{(k)}_i\ge0$, since a nonnegative $d^{(k)}_{i+1}$ makes the
numerator nonnegative and the denominator at least $1$. Hence
\[
I+\beta_1\tau P_{i+1}\Gstar=V(I+\beta_1\tau D_{i+1}\Sigma^{\!\top})V^{\!\top}
\quad\text{has eigenvalues}\quad 1+\beta_1\tau\sigma_kd^{(k)}_{i+1}\ge1,
\]
and the inverted block exists \emph{unconditionally}.

Since the denominator in the scalar recursion is at least one,
\[
d_i^{(k)}\le d_{i+1}^{(k)}+\beta_2\tau\sigma_k.
\]
Iterating from $d_N^{(k)}=\sigma_k$ gives
\[
d_i^{(k)}\le \bigl(1+\beta_2(T-t_i)\bigr)\sigma_k.
\]
Hence $\max_{0\le i\le N}\|P_i\|\le (1+\beta_2T)\|G\|$ on every grid.

For an arbitrary square-integrable adapted source $(\chi_i,\gamma_i,\eta_i,\phi_i,\zeta)$ in the
perturbed scheme \eqref{B:eq:pfwd}--\eqref{B:eq:pbwd}, as the continuation step
requires, $P_i$ is unchanged and the ansatz $Y_i=P_iX_i+p_i$ reduces the backward equation to
\[
p_i=(I+\beta_1\tau P_{i+1}\Gstar)^{-1}\bigl(\Ei[p_{i+1}]+\tau(P_{i+1}\chi_i-\phi_i)\bigr),
\qquad p_N=\zeta ;
\]
$\phi_i$ enters with a minus sign because \eqref{B:eq:pbwd} carries $+\phi_i\tau$ against the
$-f_i\tau$ drift. Let
\[
\widetilde{\Delta p}_i:=p_{i+1}-\Ei[p_{i+1}]
=Z^{p}_i\,\Delta W_i+\int_E U^{p}_i(e)\,\Delta\Nt_i(\diff e)+\Delta M_i
\]
be its decomposition from Lemma~\ref{B:lem:kw}. The martingale part of $X_{i+1}$ is
$(-\beta_1\Gstar Z_i+\gamma_i)\Delta W_i+\int_E(-\beta_1\Gstar U_i+\eta_i)\,\Delta\Nt_i$, so
matching martingale components in $Y_{i+1}=P_{i+1}X_{i+1}+p_{i+1}$ gives, by uniqueness in
Lemma~\ref{B:lem:kw},
\[
B_i\,Z_i=P_{i+1}\gamma_i+Z^{p}_i,
\qquad
B_i\,U_i=P_{i+1}\eta_i+U^{p}_i,
\qquad
B_i:=I+\beta_1P_{i+1}\Gstar,
\]
with $\Delta M_i$ the admissible remainder of Lemma~\ref{B:lem:kw}. By the orthogonal diagonalisation of
$P_{i+1}\Gstar$, $B_i$ is symmetric with eigenvalues at least $1$. Hence
\[
\|B_i^{-1}\|=\frac{1}{\lambda_{\min}(B_i)}\le1.
\]
Substituting into the forward recursion gives
$\Theta\in\mathcal S_\pi$.

Thus $E_0$ has a unique solution for every square-integrable adapted source. For
$\tau<\tau^\star_{\mathrm{stab}}$, Proposition~\ref{B:prop:pert} with the uniform constants of
step~\textup{(i)}, applied against the zero solution, gives
\[
\|\Theta\|_{\mathcal S_\pi}^2\le C\bigl(\E|\xi_0|^2+\|s\|_{\mathrm{src}}^2\bigr)
\]
with $C$ independent of $\tau,N$,
and the same difference estimate makes the
source-to-solution map $\mathcal I_0$ uniformly Lipschitz. The admissibility constraints on
$M$ are closed linear conditions in $L^2$, so $\mathcal S_\pi$ is complete.

For the continuation step, take $\tau<\tau_{\mathrm{wp}}:=\tau^\star_{\mathrm{stab}}$.

\emph{(iii) Continuation step.} Write $\mu_i=\Law(\Theta_i)$ for $0\le i<N$ and define
$\mathcal Q(\Theta)=(\chi,\gamma,\eta,\phi,\zeta)$ by
\[
\begin{gathered}
\begin{aligned}
\chi_i &= b(t_i,\Theta_i,\mu_i)-b^0(\Theta_i),
&\qquad \gamma_i &= \sigma(t_i,\Theta_i,\mu_i)-\sigma^0(\Theta_i),\\
\eta_i &= h(t_i,\Theta_i,\cdot,\mu_i)-h^0(\Theta_i),
&\qquad \phi_i &= f^0(\Theta_i)-f(t_i,\Theta_i,\mu_i),
\end{aligned}\\[3pt]
\zeta=g\bigl(X_N,\Law(X_N)\bigr)-g^0(X_N).
\end{gathered}
\]
Applying Lemma~\ref{B:lem:law} and Assumption~\ref{B:ass:lip} to the terms in
\eqref{B:eq:srcnorm}, with $L'$ from \eqref{B:eq:hom-lip-constant}, gives
\[
\begin{aligned}
&\|\mathcal Q(\Theta)-\mathcal Q(\Theta')\|_{\mathrm{src}}^2\\
&\quad\le C(L')^2\Bigg(
\sum_{i=0}^{N-1}\tau\,\E\!\left[
|\delta X_i|^2+|\delta Y_i|^2+|\delta Z_i|_F^2
+\|\delta U_i\|_{\Lnu}^2
\right]
+\E|\delta X_N|^2
\Bigg).
\end{aligned}
\]
The stability norm in \eqref{B:eq:Spi-norm} gives
\begin{equation}\label{B:eq:hom-source-lip}
\|\mathcal Q(\Theta)-\mathcal Q(\Theta')\|_{\mathrm{src}}
\le L_{\mathrm{hom}}\|\Theta-\Theta'\|_{\mathcal S_\pi},
\end{equation}
where $L_{\mathrm{hom}}>0$ depends only on $L'$ and $T$.

Taking $C_{\mathrm{stab}}$ in Proposition~\ref{B:prop:pert} with the uniform bounds $L'$ and
$\alpha_*$ from step~\textup{(i)}, set
\[
K^\star_{\mathrm{stab}}:=\sqrt{C_{\mathrm{stab}}},
\qquad
c_\star:=\frac{1}{2K^\star_{\mathrm{stab}}L_{\mathrm{hom}}}.
\]
\begin{lemma}[Continuation step]\label{B:lem:contstep}
If $E_{\lambda_0}$ is well-posed in $\mathcal S_\pi$ \emph{for every
square-integrable adapted source}, then so is $E_{\lambda_0+c}$ (again for every
such source) for every $c\in(0,\min(c_\star,1-\lambda_0)]$; $c_\star$ does not depend on
$\lambda_0$.
\end{lemma}
\begin{proof}
Write $\mathcal I_{\lambda_0}$ for the solution operator of $E_{\lambda_0}$: given a
square-integrable adapted source it returns the unique tuple of $E_{\lambda_0}$ driven
by that source. By step~(i) and Proposition~\ref{B:prop:pert},
\[
\|\mathcal I_{\lambda_0}(s^1)-\mathcal I_{\lambda_0}(s^2)\|_{\mathcal S_\pi}
\le K^\star_{\mathrm{stab}}\|s^1-s^2\|_{\mathrm{src}},
\]
where $K^\star_{\mathrm{stab}}$ is independent of $\lambda_0$. Since
\[
\mathcal G^{\lambda_0+c}=\mathcal G^{\lambda_0}+c\,(\mathcal G_{\mathrm{sch}}-\mathcal G_{\mathrm{lin}}),
\]
the perturbation convention of \eqref{B:eq:pfwd}--\eqref{B:eq:pbwd}, which adds $+\phi_i\tau$
against $-f_i\tau$, turns this generator increment at a tuple $\Theta$ into exactly the source
$c\,\mathcal Q(\Theta)$:
\[
  (\chi_i,\gamma_i,\eta_i,\phi_i)=c\,(b-b^0,\,\sigma-\sigma^0,\,h-h^0,\,f^0-f)(t_i,\Theta_i,\mu_i),
  \qquad
  \zeta=c\,(g-g^0)(X_N,\Law(X_N)).
\]
For $c\in(0,1-\lambda_0]$ define
\[
  \Phi_c:\mathcal S_\pi\to\mathcal S_\pi,\qquad
  \Phi_c(\Theta'):=\mathcal I_{\lambda_0}\!\big(c\,\mathcal Q(\Theta')\big).
\]
Then $\Theta=\Phi_c(\Theta)$ if and only if $\Theta$ solves $E_{\lambda_0+c}$, and by
\eqref{B:eq:hom-source-lip}
\[
\begin{aligned}
  \|\Phi_c(\Theta'_1)-\Phi_c(\Theta'_2)\|_{\mathcal S_\pi}
  &\le K^\star_{\mathrm{stab}}\,c\,\|\mathcal Q(\Theta'_1)-\mathcal Q(\Theta'_2)\|_{\mathrm{src}}\\
  &\le c\,K^\star_{\mathrm{stab}}L_{\mathrm{hom}}\,\|\Theta'_1-\Theta'_2\|_{\mathcal S_\pi},
\end{aligned}
\]
so $\Phi_c$ is a contraction for $c\le c_\star$, and Banach's fixed-point theorem gives a
unique solution of $E_{\lambda_0+c}$. For a fixed extra source $s$, the map
$\Phi_c^{s}(\Theta'):=\mathcal I_{\lambda_0}\bigl(s+c\,\mathcal Q(\Theta')\bigr)$ contracts with
the same constant, since $s$ cancels in differences, so $E_{\lambda_0+c}$ is well-posed
\emph{with} sources, as the next continuation step requires.
\end{proof}
Starting from $\lambda_0=0$ and setting $\lambda_{k+1}:=\lambda_k+\min\{c_\star,1-\lambda_k\}$,
Lemma~\ref{B:lem:contstep} applies at every step and $\lambda_k=1$ after
$\lceil 1/c_\star\rceil$ steps, proving existence in Theorem~\ref{B:thm:wp} with
$\tau_0=\tau_{\mathrm{wp}}$. \qed

\section{Path-regularity and coefficient-freezing estimates}

\subsection{Proof of Lemma~\ref{B:lem:regvanish}}\label{B:app:regvanish}
\begin{proof}
Write $\psi_s:=\psi(s,\Theta_s,\mu_s)$ for $\psi\in\{b,\sigma,f\}$ and
$h_s(e):=h(s,\Theta_s^-,e,\mu_s^-)$. By
Assumptions~\ref{B:ass:integ}--\ref{B:ass:lip} and $\Theta\in\mathcal S$,
\[
b,f\in L^2(\diff t\otimes\diff\Prob),\qquad
\sigma\in L^2(\diff t\otimes\diff\Prob;\R^{n\times d}),\qquad
h\in L^2(\diff t\otimes\diff\Prob\otimes\nu;\R^n).
\]

\emph{$X$ and $Y$ increments.} For $t\in[t_i,t_{i+1}]$, \eqref{B:eq:mvfbsdej} gives
\[
\begin{aligned}
X_t-X_{t_i}&=\int_{t_i}^t b_s\diff s+\int_{t_i}^t\sigma_s\diff W_s
+\int_{t_i}^t\!\int_E h_s(e)\,\Nt(\diff s,\diff e),\\[2pt]
Y_t-Y_{t_i}&=-\int_{t_i}^t f_s\diff s+\int_{t_i}^t Z_s\diff W_s
+\int_{t_i}^t\!\int_E U_s(e)\,\Nt(\diff s,\diff e).
\end{aligned}
\]
Hence, by Cauchy--Schwarz on the drifts (with $t-t_i\le\tau$) and the Brownian and
compensated-Poisson isometries,
\begin{equation}\label{B:eq:XYincr-reg}
\begin{aligned}
&\sup_{t\in[t_i,t_{i+1}]}\big(\E|X_t-X_{t_i}|^2+\E|Y_t-Y_{t_i}|^2\big)
\le C\Big[\tau\!\int_{t_i}^{t_{i+1}}\!\big(\E|b_s|^2+\E|f_s|^2\big)\diff s\\
&\qquad\qquad+\int_{t_i}^{t_{i+1}}\!\big(\E|\sigma_s|_F^2+\E\|h_s\|_{\Lnu}^2
+\E|Z_s|_F^2+\E\|U_s\|_{\Lnu}^2\big)\diff s\Big].
\end{aligned}
\end{equation}
Each integrand $g$ in \eqref{B:eq:XYincr-reg} lies in $L^1(0,T)$,
so absolute continuity of the Lebesgue integral gives
\[
\max_{0\le i<N}\int_{t_i}^{t_{i+1}}g(r)\,\diff r\;\longrightarrow\;0
\qquad\text{as }\tau\to0 .
\]
Applied to each term on the right of \eqref{B:eq:XYincr-reg}, this shows that the
$(X,Y)$-increment part of $\Reg(\tau)$ tends to zero.

\emph{$(Z,U)$ projections.} Set $V=(Z,U)$ and $K=\R^{m\times d}\times L^2(\nu;\R^m)$, with the
product Hilbert norm. Write $\|\cdot\|_2$ for the norm in $L^2(\diff t\otimes\diff\Prob;K)$. On each
grid cell, the conditional average
\[
(\Pi_\tau V)_t
:=\frac1\tau\,\Ei\!\int_{t_i}^{t_{i+1}}V_s\,\diff s,
\qquad t\in(t_i,t_{i+1}],
\]
is the orthogonal projection onto the $\F_{t_i}$-measurable random variables that are constant
in time on that cell.

Given $\varepsilon>0$, density of simple predictable processes in the predictable $L^2$
space~\cite[\S4.1]{Applebaum} gives a partition $0=s_0<\cdots<s_J=T$ and
\[
S=\sum_{j=0}^{J-1}S_j\,\mathbf 1_{(s_j,s_{j+1}]},
\qquad S_j\in L^2(\Omega,\F_{s_j};K),
\qquad \|V-S\|_2<\varepsilon.
\]
Let $B_\tau$ be the union of grid cells containing a breakpoint $s_j$ strictly in their
interior. On every other cell, $S=S_j$ for some $s_j\le t_i$, so $\Pi_\tau S=S$ there. Since
$B_\tau$ contains at most $J-1$ cells, orthogonality gives
\[
\|S-\Pi_\tau S\|_2^2
\le \E\!\int_{B_\tau}\|S_t\|_K^2\,\diff t
\le (J-1)\tau\max_j\E\|S_j\|_K^2
\longrightarrow 0.
\]
As $\Pi_\tau$ is a contraction,
\[
\|V-\Pi_\tau V\|_2
\le 2\|V-S\|_2+\|S-\Pi_\tau S\|_2.
\]
Taking $\tau\to0$ and then $\varepsilon\to0$ gives
\[
\sum_{i=0}^{N-1}\E\!\int_{t_i}^{t_{i+1}}
\bigl(|Z_s-\sm Z_i|_F^2
+\|U_s-\sm U_i\|_{\Lnu}^2\bigr)\diff s
=\|V-\Pi_\tau V\|_2^2\longrightarrow0.
\]
Together with the $X,Y$ increment estimate \eqref{B:eq:XYincr-reg}, this proves
$\Reg(\tau)\to0$.
\end{proof}

\subsection{Freezing error for the sampled continuous tuple}
\label{B:app:disc}
\begin{lemma}\label{B:lem:freezing-error}
Let Assumptions~\ref{B:ass:integ}--\ref{B:ass:lip} hold, let $\Theta=(X,Y,Z,U)\in\mathcal S$ be a
solution of \eqref{B:eq:mvfbsdej}, and let the time-modulus \eqref{B:eq:omega} satisfy
$\omega(\tau)<\infty$. For $s\in[t_i,t_{i+1})$ write the coefficient-freezing differences
$\delta_i b(s):=b(s,\Theta_s,\mu_s)-b(t_i,\sm\Theta_i,\sm\mu_i)$, analogously
$\delta_i\sigma(s),\delta_i f(s)$, and
$\delta_i h(s,e):=h(s,\Theta_s^-,e,\mu_s^-)-h(t_i,\sm\Theta_i,e,\sm\mu_i)$. Then
\begin{equation}\label{B:eq:freezesum}
\sum_{i=0}^{N-1}\int_{t_i}^{t_{i+1}}\!\E\big[|\delta_i b(s)|^2+|\delta_i\sigma(s)|_F^2+|\delta_i f(s)|^2+\|\delta_i h(s,\cdot)\|_{\Lnu}^2\big]\diff s\le C\big(\omega(\tau)^2+\Reg(\tau)\big),
\end{equation}
and the auxiliary tuple $\au\Theta$ of \eqref{B:eq:auxcont} satisfies
\begin{equation}
\label{B:eq:discgap}
\max_{0\le i\le N}\E|\sm X_i-\au X_i|^2+\max_{0\le i\le N}\E|\sm Y_i-\au Y_i|^2\le C\big(\omega(\tau)^2+\Reg(\tau)\big).
\end{equation}
\end{lemma}
\begin{proof}
We estimate $\delta_i h$; the bounds for $b,\sigma,f$ are identical in their respective
norms. Split it as
\[
\delta_i h(s,\cdot)=\big[h(s,\Theta_s^-,\cdot,\mu_s^-)-h(s,\sm\Theta_i,\cdot,\sm\mu_i)\big]
+\big[h(s,\sm\Theta_i,\cdot,\sm\mu_i)-h(t_i,\sm\Theta_i,\cdot,\sm\mu_i)\big].
\]
Both input pairs lie in $\mathcal D_s$, since $\F_{t_i}\subset\F_s$, so
Assumption~\ref{B:ass:lip} bounds the first bracket for $\diff s$-a.e.\ $s$. The time-modulus
\eqref{B:eq:omega} bounds the second, whose input is frozen at $(\sm\Theta_i,\sm\mu_i)$. Hence
\[
\begin{aligned}
\|\delta_i h(s,\cdot)\|_{\Lnu}^2\le{}& C\,\omega(\tau)^2\big(1+\|\sm\Theta_i\|_{\mathcal H}^2+W_2(\sm\mu_i,\delta_0)^2\big)\\
&+C\big(|X_{s-}-X_{t_i}|^2+|Y_{s-}-Y_{t_i}|^2\big)\\
&+C\big(|Z_s-\sm Z_i|_F^2+\|U_s-\sm U_i\|_{\Lnu}^2+W_2(\mu_s^-,\sm\mu_i)^2\big).
\end{aligned}
\]
The synchronous coupling of Lemma~\ref{B:lem:law} gives
\[
W_2(\mu_s^-,\sm\mu_i)^2\le\E\|\Theta_s^--\sm\Theta_i\|_{\mathcal H}^2
=\E\big[|X_{s-}-X_{t_i}|^2+|Y_{s-}-Y_{t_i}|^2+|Z_s-\sm Z_i|_F^2+\|U_s-\sm U_i\|_{\Lnu}^2\big].
\]
For the time-freezing part, $W_2(\sm\mu_i,\delta_0)^2=\E\|\sm\Theta_i\|_{\mathcal H}^2$ and the
grid projection is a contraction, $\tau\,\E|\sm Z_i|_F^2\le\E\int_{t_i}^{t_{i+1}}|Z_s|_F^2\diff s$
and likewise for $\sm U_i$, so
\[
\begin{aligned}
&\sum_{i=0}^{N-1}\int_{t_i}^{t_{i+1}}\E\big[1+\|\sm\Theta_i\|_{\mathcal H}^2+W_2(\sm\mu_i,\delta_0)^2\big]\diff s\\
&\qquad\le C\Big(1+\sup_{t\le T}\E|X_t|^2+\sup_{t\le T}\E|Y_t|^2
+\E\!\int_0^T\!\big(|Z_s|_F^2+\|U_s\|_{\Lnu}^2\big)\diff s\Big)\le C,
\end{aligned}
\]
the last because $\Theta\in\mathcal S$. Since $X_{s-}=X_s$ and $Y_{s-}=Y_s$ for
$\diff s\otimes\diff\Prob$-a.e.\ $(s,\omega)$, the definition \eqref{B:eq:Reg} gives
\[
\sum_{i=0}^{N-1}\int_{t_i}^{t_{i+1}}\E\big[|X_{s-}-X_{t_i}|^2+|Y_{s-}-Y_{t_i}|^2
+|Z_s-\sm Z_i|_F^2+\|U_s-\sm U_i\|_{\Lnu}^2\big]\diff s\le(1+T)\,\Reg(\tau).
\]
Integrating the pointwise bound over each cell and summing over $i$ proves
\eqref{B:eq:freezesum}.

For \eqref{B:eq:discgap}, subtract \eqref{B:eq:auxcont} from \eqref{B:eq:mvfbsdej} on each
cell and sum, using $\sm X_0=\au X_0$ and $\sm Y_0=\au Y_0$: for $0\le k\le N$,
\[
\begin{aligned}
\sm X_k-\au X_k&=\sum_{i=0}^{k-1}\Big(\int_{t_i}^{t_{i+1}}\!\delta_i b(s)\diff s
+\int_{t_i}^{t_{i+1}}\!\delta_i\sigma(s)\diff W_s
+\int_{t_i}^{t_{i+1}}\!\!\int_E\delta_i h(s,e)\,\Nt(\diff s,\diff e)\Big),\\
\sm Y_k-\au Y_k&=\sum_{i=0}^{k-1}\Big(-\int_{t_i}^{t_{i+1}}\!\delta_i f(s)\diff s
+\int_{t_i}^{t_{i+1}}\!(Z_s-\sm Z_i)\diff W_s
+\int_{t_i}^{t_{i+1}}\!\!\int_E(U_s-\sm U_i)\,\Nt(\diff s,\diff e)\Big).
\end{aligned}
\]
Cauchy--Schwarz on the drift sums, and Doob's inequality with the It\^o and
compensated-Poisson isometries on the martingale sums, give
\[
\begin{aligned}
\E\max_{0\le k\le N}|\sm X_k-\au X_k|^2
&\le C\sum_{i=0}^{N-1}\int_{t_i}^{t_{i+1}}\E\big[|\delta_i b(s)|^2+|\delta_i\sigma(s)|_F^2
+\|\delta_i h(s,\cdot)\|_{\Lnu}^2\big]\diff s,\\
\E\max_{0\le k\le N}|\sm Y_k-\au Y_k|^2
&\le C\sum_{i=0}^{N-1}\int_{t_i}^{t_{i+1}}\E\big[|\delta_i f(s)|^2+|Z_s-\sm Z_i|_F^2
+\|U_s-\sm U_i\|_{\Lnu}^2\big]\diff s.
\end{aligned}
\]
By \eqref{B:eq:freezesum} and \eqref{B:eq:Reg}, both right-hand sides are bounded by
$C\bigl(\omega(\tau)^2+\Reg(\tau)\bigr)$.
\end{proof}

\section{Plug-in scoring}\label{B:app:lawbridge}
Theorem~\ref{B:thm:discrete} concerns the population estimator \eqref{B:eq:Epi}, whose interior
coefficients read the approximation's own laws $\ap\mu_i$. In solver-side evaluation we replace the
law inputs of the coefficients by estimates formed on an independent batch and score the resulting
\emph{plug-in} functional. This appendix bounds the gap between the population and plug-in
functionals and gives the additional error term required in the two-sided estimate.
Proposition~\ref{B:prop:lawbridge} treats a general plug-in law flow, and
Corollary~\ref{B:cor:statbridge} the finite-dimensional statistics used by the reported solver
studies.

\begin{assumption}[Off-diagonal law regularity; scoring only]\label{B:ass:offdiag}
There exists $L_\mu<\infty$ such that, for every $t\in[0,T]$, $\theta\in\mathcal H$ and
$\mu,\varrho\in\mathcal P_2(\mathcal H)$,
\[
\|\psi(t,\theta,\mu)-\psi(t,\theta,\varrho)\|
\le L_\mu W_2(\mu,\varrho),
\qquad \psi\in\{b,\sigma,h,f\},
\]
and, for every $x\in\R^n$ and $\rho,\rho'\in\mathcal P_2(\R^n)$,
\[
|g(x,\rho)-g(x,\rho')|
\le L_\mu W_2(\rho,\rho').
\]
\end{assumption}
All benchmark coefficients satisfy~\ref{B:ass:offdiag}.

\begin{proposition}[Population/plug-in bridge]\label{B:prop:lawbridge}
Let Assumption~\ref{B:ass:offdiag} hold. Fix $\ap\Theta\in\mathcal S_\pi$ with
$\Epi(\ap\Theta)<\infty$, and let $\varrho=(\varrho_i)_{i<N}\subset\mathcal P_2(\mathcal H)$,
$\varrho_N\in\mathcal P_2(\R^n)$ be any law flow. Write
$\Epi^{\mathrm{plug}}(\varrho)$ for \eqref{B:eq:Epi} with $\ap\mu_i$ replaced
by $\varrho_i$ in every interior coefficient evaluation and $\Law(\ap X_N)$ by
$\varrho_N$ in the terminal datum, and set
\[
\mathcal W(\varrho):=\max\Bigl\{\max_{i<N}W_2(\ap\mu_i,\varrho_i),\;
W_2\bigl(\Law(\ap X_N),\varrho_N\bigr)\Bigr\}.
\]
Then, with a universal $C$,
\[
\Bigl|\sqrt{\Epi^{\mathrm{plug}}(\varrho)}-\sqrt{\Epi}\Bigr|\;\le\;C\,L_\mu(1+T)\,\mathcal W(\varrho),
\]
and consequently
\[
\Epi\le2\,\Epi^{\mathrm{plug}}(\varrho)+CL_\mu^2(1+T)^2\,\mathcal W(\varrho)^2,
\qquad
\Epi^{\mathrm{plug}}(\varrho)\le2\,\Epi+CL_\mu^2(1+T)^2\,\mathcal W(\varrho)^2 .
\]
\end{proposition}
\begin{proof}
For $\psi\in\{b,\sigma,h,f\}$ write
\[
\delta_\varrho\psi_i:=\psi(t_i,\ap\Theta_i,\varrho_i)-\psi(t_i,\ap\Theta_i,\ap\mu_i),\qquad
\delta_\varrho g:=g(\ap X_N,\varrho_N)-g(\ap X_N,\Law(\ap X_N)),
\]
with $\ap\Theta_i:=(\ap X_i,\ap Y_i,\ap Z_i,\ap U_i)$. Assumption~\ref{B:ass:offdiag} gives,
for every $i<N$,
\[
\|\delta_\varrho b_i\|_{L^2},\ \|\delta_\varrho\sigma_i\|_{L^2},\
\|\delta_\varrho h_i\|_{L^2(\Omega;\Lnu)},\ \|\delta_\varrho f_i\|_{L^2},\
\|\delta_\varrho g\|_{L^2}\ \le\ L_\mu\mathcal W(\varrho).
\]
Only the coefficient evaluations change, so by \eqref{B:eq:CX}--\eqref{B:eq:CY}
\[
C^X_{i,\mathrm{plug}}-C^X_i=-\sum_{j=0}^{i}\Big(\tau\,\delta_\varrho b_j
+\delta_\varrho\sigma_j\Delta W_j+\int_E\delta_\varrho h_j\,\Delta\Nt_j\Big),\qquad
C^Y_{i,\mathrm{plug}}-C^Y_i=\sum_{j=0}^{i}\tau\,\delta_\varrho f_j,
\]
while $r^X_{0,\mathrm{plug}}=r^X_0$ and $r^Y_{T,\mathrm{plug}}-r^Y_T=-\delta_\varrho g$. The
increment isometries and the orthogonality of the Brownian and compensated-Poisson
increments give, uniformly in $i<N$,
\[
\|C^X_{i,\mathrm{plug}}-C^X_i\|_{L^2}\le\bigl(T+\sqrt{2T}\bigr)L_\mu\mathcal W(\varrho),
\qquad
\|C^Y_{i,\mathrm{plug}}-C^Y_i\|_{L^2}\le T\,L_\mu\mathcal W(\varrho).
\]
Collect the four estimator components in
\[
\mathbf R:=\Bigl(\|r_0^X\|_{L^2},\ \|r_T^Y\|_{L^2},\ \max_{0\le i\le N-1}\|C_i^X\|_{L^2},\
\max_{0\le i\le N-1}\|C_i^Y\|_{L^2}\Bigr),
\]
and $\mathbf R^{\mathrm{plug}}$ likewise, so that $\sqrt{\Epi}=|\mathbf R|$ and
$\sqrt{\Epi^{\mathrm{plug}}(\varrho)}=|\mathbf R^{\mathrm{plug}}|$. Since
$|\max_ia_i-\max_ib_i|\le\max_i|a_i-b_i|$ and $T+\sqrt{2T}\le1+2T$,
\[
\Bigl|\sqrt{\Epi^{\mathrm{plug}}(\varrho)}-\sqrt{\Epi}\Bigr|
\le|\mathbf R^{\mathrm{plug}}-\mathbf R|
\le C\,L_\mu(1+T)\,\mathcal W(\varrho).
\]
With $K:=CL_\mu(1+T)\mathcal W(\varrho)$, this gives
$\sqrt{\Epi}\le\sqrt{\Epi^{\mathrm{plug}}(\varrho)}+K$, hence
$\Epi\le2\,\Epi^{\mathrm{plug}}(\varrho)+2K^2$, and the reverse inequality in the same way.
\end{proof}

\begin{corollary}[Statistic plug-in]\label{B:cor:statbridge}
Suppose the coefficients depend on the law only through finite-dimensional statistics: there are
maps $S_i:\mathcal P_2(\mathcal H)\to\R^{q}$ for $i<N$ and $S_N:\mathcal P_2(\R^n)\to\R^{q_N}$, and
functions $\tilde b,\tilde f,\tilde g$, such that
\[
b(t_i,\theta,\mu)=\tilde b\bigl(t_i,\theta,S_i(\mu)\bigr),\qquad
f(t_i,\theta,\mu)=\tilde f\bigl(t_i,\theta,S_i(\mu)\bigr),\qquad
g(x,\rho)=\tilde g\bigl(x,S_N(\rho)\bigr),
\]
while $\sigma$ and $h$ do not depend on the law, and, in the Euclidean norm,
\[
\begin{gathered}
|\tilde b(t_i,\theta,s)-\tilde b(t_i,\theta,s')|\le K_b|s-s'|,\qquad
|\tilde f(t_i,\theta,s)-\tilde f(t_i,\theta,s')|\le K_f|s-s'|,\\
|\tilde g(x,s)-\tilde g(x,s')|\le K_g|s-s'|.
\end{gathered}
\]
Fix $\ap\Theta\in\mathcal S_\pi$ with $\Epi(\ap\Theta)<\infty$. Let $s=((s_i)_{i<N},s_N)$ be a
deterministic statistic path. Write $\Epi^{\mathrm{plug}}(s)$ for
\eqref{B:eq:Epi} with $\tilde b,\tilde f$ evaluated at $s_i$ and $\tilde g$ at $s_N$, and set
\[
\varepsilon_{\mathrm{stat}}(s):=\max\Bigl\{\max_{i<N}\bigl|s_i-S_i(\ap\mu_i)\bigr|,\;
\bigl|s_N-S_N\bigl(\Law(\ap X_N)\bigr)\bigr|\Bigr\}.
\]
Then $\Epi^{\mathrm{plug}}(s)=\Epi$ at the own statistic path
$s=\bigl((S_i(\ap\mu_i))_{i<N},S_N(\Law(\ap X_N))\bigr)$, and
\[
\Bigl|\sqrt{\Epi^{\mathrm{plug}}(s)}-\sqrt{\Epi}\Bigr|
\le\bigl(T^2K_b^2+T^2K_f^2+K_g^2\bigr)^{1/2}\,\varepsilon_{\mathrm{stat}}(s).
\]
\end{corollary}
\begin{proof}
At the candidate's own statistic path, the coefficient evaluations coincide with those in
\eqref{B:eq:Epi}, so $\Epi^{\mathrm{plug}}(s)=\Epi$.

For general $s$, only $b$, $f$ and $g$ change. The residual identities in the proof of
Proposition~\ref{B:prop:lawbridge} and the Lipschitz bounds in the statistic give, for $i<N$,
\[
\begin{aligned}
\|C^X_{i,\mathrm{plug}}-C^X_i\|_{L^2}
&\le K_b\sum_{j=0}^{i}\tau\,|s_j-S_j(\ap\mu_j)|
&&\le TK_b\,\varepsilon_{\mathrm{stat}}(s),\\
\|C^Y_{i,\mathrm{plug}}-C^Y_i\|_{L^2}
&\le K_f\sum_{j=0}^{i}\tau\,|s_j-S_j(\ap\mu_j)|
&&\le TK_f\,\varepsilon_{\mathrm{stat}}(s),\\
\|r^Y_{T,\mathrm{plug}}-r^Y_T\|_{L^2}
&\le K_g\,|s_N-S_N(\Law(\ap X_N))|
&&\le K_g\,\varepsilon_{\mathrm{stat}}(s).
\end{aligned}
\]
The initial residual $r_0^X$ is unchanged. By the reverse triangle inequality,
\[
\Bigl|\sqrt{\Epi^{\mathrm{plug}}(s)}-\sqrt{\Epi}\Bigr|
\le \bigl|\mathbf R^{\mathrm{plug}}-\mathbf R\bigr|
\le \sqrt{T^2K_b^2+T^2K_f^2+K_g^2}\,\varepsilon_{\mathrm{stat}}(s).
\]
\end{proof}

In the reported solver studies, the law dependence factors through the nodewise mean of $X$, and
additionally that of $Y$ in the inventory driver; the displayed coefficients are globally
Lipschitz in these statistics.

\begin{remark}[On the statistic plug-in]\label{B:rem:lawplugin}
Fix a realisation of the training and independent pilot data. The approximation rule and
pilot statistics $s^A$ are then fixed, while $\ap\Theta$ remains random on the evaluation
space of Section~\ref{B:sec:deepsolver}, on which all laws, statistics and expectations below
are taken. Empirical means over the scoring batch estimate $\Epi^{\mathrm{plug}}(s^A)$.
In particular, $\varepsilon_{\mathrm{stat}}:=\varepsilon_{\mathrm{stat}}(s^A)$ is measured against
the approximation's own statistics on this evaluation space.
Let $C_{\mathrm{eff}}$ be the efficiency constant in Theorem~\ref{B:thm:discrete}, and
write $C_{\mathrm{rel}}$ for its reliability constant. Set
$K_A:=\bigl(T^2K_b^2+T^2K_f^2+K_g^2\bigr)^{1/2}\varepsilon_{\mathrm{stat}}$.
Corollary~\ref{B:cor:statbridge} gives
\[
\sqrt{\Epi}\ge\Bigl(\sqrt{\Epi^{\mathrm{plug}}}-K_A\Bigr)_{\!+},\qquad
\sqrt{\Epi}\le\sqrt{\Epi^{\mathrm{plug}}}+K_A,
\]
and efficiency $\Epi\le C_{\mathrm{eff}}\Dpi$ and reliability $\Dpi\le C_{\mathrm{rel}}\Epi$
of Theorem~\ref{B:thm:discrete}, on grids satisfying the reliability restriction, then give
\[
C_{\mathrm{eff}}^{-1}\Bigl(\sqrt{\Epi^{\mathrm{plug}}}-K_A\Bigr)_{\!+}^{2}
\;\le\;\Dpi\;\le\;
C_{\mathrm{rel}}\Bigl(\sqrt{\Epi^{\mathrm{plug}}}+K_A\Bigr)^{2}.
\]
\end{remark}

\section{Benchmark calculations}

\subsection{Path regularity and exact effectivity}\label{B:app:deferred}

\medskip\noindent\emph{Path regularity of the affine benchmarks (Proposition~\ref{B:prop:benchreg}).}
\begin{proof}
Write $H_0(e)=h_0$ for the single-mark benchmark and
$H_0(e)=h_0\chi_\infty(e)$ for the weighted infinite-activity benchmark. Both profiles lie in
$\Lnu$: in the latter case $\chi_\infty(e)=e\,\chi(e)$ and $\nu(\diff e)=\Lambda e^{-2}\diff e$
on $E=(0,1]$, so
\[
\|\chi_\infty\|_{\Lnu}^2=\int_0^1|e\,\chi(e)|^2\,\Lambda e^{-2}\diff e
=\Lambda\int_0^1|\chi(e)|^2\diff e<\infty.
\]
Both benchmarks have centred initial data. Their continuous solutions satisfy
\[
\E X_t=\E Y_t=0,\qquad Y_t=P(t)X_t,
\]
where $P(T)=G$ and $\dot P=\kappa_{\mathrm F}P\Gstar P-\kappa_{\mathrm B}G$. With the singular value
decomposition $G=V\Sigma Q^{\!\top}$, $P(t)=VD(t)Q^{\!\top}$.
The diagonal entries satisfy
\[
\dot d_k(t)=\kappa_{\mathrm F}\sigma_k d_k(t)^2-\kappa_{\mathrm B}\sigma_k,
\qquad d_k(T)=\sigma_k\ge 0.
\]
Since $\dot d_k=-\kappa_{\mathrm B}\sigma_k\le 0$ whenever $d_k=0$, the solution remains nonnegative as $t$
decreases from $T$. Integrating from $t$ to $T$ gives
\[
0\le d_k(t)
=\sigma_k+\kappa_{\mathrm B}\sigma_k(T-t)
-\kappa_{\mathrm F}\sigma_k\int_t^T d_k(r)^2\,\diff r
\le \sigma_k\bigl(1+\kappa_{\mathrm B}(T-t)\bigr).
\]
Thus each $d_k$ is bounded and extends to $[0,T]$. Since $D(t)\Sigma^{\!\top}$ is diagonal with
nonnegative entries,
\[
I+\kappa_{\mathrm F}P(t)\Gstar
=V\bigl(I+\kappa_{\mathrm F}D(t)\Sigma^{\!\top}\bigr)V^{\!\top}
\]
is positive definite. Hence $\mathcal A(t):=(I+\kappa_{\mathrm F}P(t)\Gstar)^{-1}P(t)$ is $C^1$ on $[0,T]$.
\[
Z(t)=\mathcal A(t)\sigma_0,\qquad U(t,\cdot)=\mathcal A(t)H_0(\cdot),
\]
so $Z\in C^1([0,T];\R^{m\times d})$ and $U\in C^1([0,T];\Lnu)$, with
$\|U(t)-U(s)\|_{\Lnu}\le|\mathcal A(t)-\mathcal A(s)|\,\|H_0\|_{\Lnu}$. The jump coefficient
$h(t,\cdot)=H_0-\kappa_{\mathrm F}\Gstar U(t,\cdot)$ has the same regularity, so
\[
H_h^2:=\sup_{t\le T}\|h(t,\cdot)\|_{\Lnu}^2<\infty
\]
for both profiles.

\begin{itemize}
\item[\textup{(i)}] The $Z,U$-projection contribution is $O(\tau^2)$: $Z$ and $U$ are
Lipschitz in time, so there is $L_*$ with
\[
|Z_s-\sm Z_i|_F+\|U_s-\sm U_i\|_{\Lnu}\le L_*\tau
\]
on $[t_i,t_{i+1}]$, and
\[
\sum_{i=0}^{N-1}\E\!\int_{t_i}^{t_{i+1}}\bigl(|Z_s-\sm Z_i|_F^2+\|U_s-\sm U_i\|_{\Lnu}^2\bigr)\diff s
\le T L_*^2\tau^2.
\]

\item[\textup{(ii)}] The state increment is $O(\tau)$. The drift $b$ is affine in $(X,Y)$, so
$\sup_{s\le T}\E|b_s|^2<\infty$, and $\sigma$ is deterministic. Cauchy--Schwarz on the drift,
the Brownian isometry and the compensated-Poisson isometry
\[
\E\Big|\int_{t_i}^{t}\!\!\int_Eh(s,e)\,\Nt(\diff s,\diff e)\Big|^2
=\int_{t_i}^{t}\|h(s,\cdot)\|_{\Lnu}^2\,\diff s
\le H_h^2(t-t_i)
\]
give
\[
\E|X_t-X_{t_i}|^2
\le3\Bigl(\tau\sup_{s\le T}\E|b_s|^2+\sup_{s\le T}|\sigma(s)|_F^2+H_h^2\Bigr)(t-t_i)=O(\tau).
\]

\item[\textup{(iii)}] The value increment is $O(\tau)$: since $P\in C^1$ and
$\sup_t\E|X_t|^2<\infty$,
\[
Y_t-Y_{t_i}=P(t_i)(X_t-X_{t_i})+R_{i,t},\qquad \E|R_{i,t}|^2=O(\tau^2),
\]
so $\E|Y_t-Y_{t_i}|^2=O(\tau)$ by \textup{(ii)}.
\end{itemize}
Hence $\Reg(\tau)=O(\tau)$ for both benchmarks.
\end{proof}

\medskip\noindent\emph{Exact effectivity for integrand perturbations.}
In the scalar state-law affine benchmark with deterministic jump profile $h_0(e)$, let
$V_i:=\ap Z_i-\ex Z_i$ and $H_i:=\ap U_i-\ex U_i$ be arbitrary
square-integrable $\F_{t_i}$-measurable integrands, with $X,Y,M$
unchanged. Set
\[
S_i:=\sum_{j=0}^{i}\Bigl(V_j\Delta W_j+
\int_E H_j(e)\,\Delta\Nt_j(\diff e)\Bigr),\qquad
A:=\sum_{i=0}^{N-1}\tau\,\E\bigl[|V_i|^2+\|H_i\|_{\Lnu}^2\bigr].
\]
Since $X,Y,M$ are unchanged and the benchmark's law dependence is through $X$ alone, $b$, $f$,
$g$ and both boundary residuals are unchanged; only the $\sigma,h$ feedback and the explicit
$Z,U$ terms of the cumulative residuals change. From $\sigma=\sigma_0-\kappa_{\mathrm F}Z$ and
$h(e)=h_0(e)-\kappa_{\mathrm F}U(e)$,
\[
C_i^X=\kappa_{\mathrm F} S_i,\qquad C_i^Y=-S_i .
\]
By the conditional increment isometries and Brownian--Poisson orthogonality,
\[
\E|S_i|^2=\sum_{j=0}^{i}\tau\,\E\bigl[|V_j|^2+\|H_j\|_{\Lnu}^2\bigr],
\qquad\text{so}\qquad
\max_{0\le i\le N-1}\E|S_i|^2=A.
\]
Consequently
\[
\Dpi=A,\qquad \Epi=(1+\kappa_{\mathrm F}^{\,2})A,\qquad
\frac{\Dpi}{\Epi}=(1+\kappa_{\mathrm F}^{\,2})^{-1}\quad(A>0).
\]

\subsection{Monotonicity of the inventory benchmark}\label{B:app:inv-admiss}
The benchmark satisfies Assumptions~\ref{B:ass:integ}--\ref{B:ass:mono}.

Take $G=1$ ($n=m=d=1$; Assumption~\ref{B:ass:G} holds with $c_G=1$). For diagonal inputs
$\delta m=\E[\delta X]$; write $\bar x:=\E[\delta X]$ and $\bar y:=\E[\delta Y]$. With
constant $\sigma_0,h_0$ ($\delta\sigma=\delta h=0$),
\[
\delta b=-\kappa\,\delta X-\bar a\,\bar x-\tfrac1R\,\delta Y,\qquad
-\delta f=\kappa\,\delta Y+\bar a\,\bar y-Q(\delta X-\bar x)-\bar Q\,\bar x .
\]
In the pairing the $\kappa$- and $\bar a$-terms cancel,
\[
\E\inner{-\kappa\,\delta X}{\delta Y}+\E\inner{\kappa\,\delta Y}{\delta X}=0,
\qquad
-\bar a\,\bar x\,\bar y+\bar a\,\bar y\,\bar x=0,
\]
the second because the backward mean term is the adjoint of the forward crowding term;
without this cancellation an additional restriction on $\bar a$ would be required. Hence
\[
\E\big[\inner{\delta b}{\delta Y}+\inner{-\delta f}{\delta X}\big]
= -Q\,\E|\delta X-\bar x|^2 - \tfrac1R\,\E|\delta Y-\bar y|^2
- \bar Q\,\bar x^2 - \tfrac1R\,\bar y^2 .
\]

The additive noise contributes no $\delta Z$- or $\delta U$-coercive terms, so $\beta_1=0$ is
unavoidable. Since $n=m$, the reliability estimate uses Case~B with coercivity parameter
$\min\{\alpha,\beta_2\}$. Under
$Q,\bar Q,R,Q_T,\bar Q_T>0$, Assumption~\ref{B:ass:mono} holds with
\[
\beta_1=0,\qquad \beta_2=\min\{Q,\bar Q\}>0,\qquad \alpha=\min\{Q_T,\bar Q_T\}>0;
\]
the parameters of Section~\ref{B:ssec:lqinv} give $\beta_2=\min\{1,\,0.8\}=0.8$ and
$\alpha=\min\{0.9,\,0.7\}=0.7$.

\subsection{Control regret of the inventory feedback}\label{B:app:regret}
\begin{proof}[Proof of Proposition~\ref{B:prop:regret}]
The orthogonal split $v=\tilde v+\bar v\mathbf 1$ gives $|v|_p^2=|\tilde v|_p^2+\bar v^2$, and the
particle recursion splits into
\[
\tilde X_{i+1}=\ell\tilde X_i+\tau\tilde u_i+\tilde\xi_i,\qquad
\bar X_{i+1}=\bar\ell\,\bar X_i+\tau\bar u_i+\bar\xi_i,\qquad
\xi_i^a:=\sigma_0\Delta W_i^a+h_0\Delta\Nt_i^a,
\]
with $\xi_i$ centred and independent of $\F_{t_i}$. The cost splits in the same way.

\emph{(i)} For $\F_{t_i}$-measurable $x$ and $u$, completing the square with
$\ell-\tau k_i=\ell R/(R+\tau S_{i+1})$ gives
\[
\tau R|\tilde u|_p^2+S_{i+1}|\ell\tilde x+\tau\tilde u|_p^2
=(S_i-\tau Q)|\tilde x|_p^2+\tau(R+\tau S_{i+1})|\tilde u+k_i\tilde x|_p^2,
\]
and the analogue for the mean. With $V_i(x):=S_i|\tilde x|_p^2+\bar S_i\bar x^2$, and since the
cross terms with $\xi_i$ vanish under $\Ei$,
\[
\begin{aligned}
&\Ei\big[\tau\big(Q|\tilde X_i|_p^2+\bar Q\bar X_i^2+R|u_i|_p^2\big)+V_{i+1}(X_{i+1})\big]\\
&\qquad=V_i(X_i)+c_i+\tau\big[(R+\tau S_{i+1})|\tilde u_i+k_i\tilde X_i|_p^2
+(R+\tau\bar S_{i+1})|\bar u_i+\bar k_i\bar X_i|^2\big],
\end{aligned}
\]
where $c_i:=S_{i+1}\E|\tilde\xi_i|_p^2+\bar S_{i+1}\E\bar\xi_i^2$ does not depend on $u$. Summing
over $i$ and taking expectations, with $V_N$ the terminal cost, gives
$J_{\pi,p}(u)=\E V_0(X_0)+\sum_ic_i+\E\sum_i\tau[\,\cdot\,]$. The last term vanishes for
$u^\star$, which is therefore optimal, and (i) follows.

\emph{(ii)} Write $X,Y$ for the components of $\ap\Theta$ and
$e_i:=Y_i-P_i\tilde X_i-\Pi_i\bar X_i\mathbf 1$. The forward recursion gives
$\tilde X_{i+1}=\varphi_i\tilde X_i-\tau R^{-1}\tilde e_i+\tilde\xi_i$, and by~\eqref{B:eq:CY}
the backward residual $\rho_i:=C^Y_i-C^Y_{i-1}$, with $C^Y_{-1}:=0$, satisfies
\[
\tilde Y_{i+1}=(1+\kappa\tau)\tilde Y_i-\tau Q\tilde X_i+\tilde\rho_i+\tilde\eta_i,\qquad
\eta_i:=\ap Z_i\Delta W_i+\int_E\ap U_i\,\Delta\Nt_i+\Delta\ap M_i,\qquad\Ei[\eta_i]=0.
\]
Substituting $\tilde Y=\tilde e+P\tilde X$, the recursion for $P$ cancels the $\tilde X_i$-terms and
\[
\tilde e_{i+1}=D_i\tilde e_i+\tilde\rho_i+\tilde\eta_i-P_{i+1}\tilde\xi_i,\qquad
D_i:=1+\kappa\tau+\tau P_{i+1}/R=2-\varphi_{i+1}\ge1.
\]
Hence $\tilde e_i=D_i^{-1}\Ei[\tilde e_{i+1}-\tilde\rho_i]$ and, with
$\gamma_j:=\prod_{l=i}^jD_l^{-1}\in(0,1]$ decreasing in $j$,
\[
\tilde e_i=\Ei\Big[\gamma_{N-1}\tilde r_T^Y-\sum_{j=i}^{N-1}\gamma_j\tilde\rho_j\Big].
\]
Summation by parts gives
\[
\sum_{j=i}^{N-1}\gamma_j\tilde\rho_j
=\gamma_{N-1}\tilde C^Y_{N-1}
-\gamma_i\tilde C^Y_{i-1}
+\sum_{j=i}^{N-2}(\gamma_j-\gamma_{j+1})\tilde C^Y_j.
\]
Since $\gamma_j$ is nonincreasing,
\[
\gamma_{N-1}+\gamma_i
+\sum_{j=i}^{N-2}|\gamma_j-\gamma_{j+1}|
=2\gamma_i\le 2.
\]
The $L^2$ triangle inequality and conditional Jensen give
\[
\begin{aligned}
\E|\tilde e_i|_p^2
&\le \Bigl(
\bigl(\E|\tilde r_T^Y|_p^2\bigr)^{1/2}
+2\max_{0\le j<N}\bigl(\E|\tilde C_j^Y|_p^2\bigr)^{1/2}
\Bigr)^2\\
&\le 5\Bigl(
\E|\tilde r_T^Y|_p^2
+\max_{0\le j<N}\E|\tilde C_j^Y|_p^2
\Bigr).
\end{aligned}
\]
For $\bar e_i$, the same argument applies with
\[
\bar D_i:=1+(\kappa+\bar a)\tau+\tau\Pi_{i+1}/R
=2-\bar\varphi_{i+1}\ge1.
\]
Adding the bounds for $\tilde e_i$ and $\bar e_i$, and using $r_0^X=0$ and $C^X\equiv0$, gives
\[
\max_i\E|\tilde e_i|_p^2+\max_i\E|\bar e_i|^2
\le5\,\E|r_T^Y|_p^2+10\max_j\E|C_j^Y|_p^2
\le10\,\Epip(\ap\Theta).
\]
Let $\delta_i:=X_i-\ex X_i$. The candidate and the Euler solution driven by $u^{\mathrm{FB}}$ use
the same initial data and noise, so $\delta_0=0$ and
\[
\begin{aligned}
\tilde\delta_{i+1}
&=\varphi_i\tilde\delta_i-\tau R^{-1}\tilde e_i,\\
\bar\delta_{i+1}
&=\bar\varphi_i\bar\delta_i-\tau R^{-1}\bar e_i.
\end{aligned}
\]
Since $\varphi_i,\bar\varphi_i\in[0,1]$, iterating and applying Cauchy--Schwarz gives
\[
\begin{aligned}
\E|\tilde\delta_i|_p^2
&\le T^2R^{-2}\max_j\E|\tilde e_j|_p^2,\\
\E|\bar\delta_i|^2
&\le T^2R^{-2}\max_j\E|\bar e_j|^2.
\end{aligned}
\]

For a control $u$ with state $X$, set
\[
\begin{aligned}
d_i(u)&:=\tilde u_i+k_i\tilde X_i
       +(\bar u_i+\bar k_i\bar X_i)\mathbf 1,\\
\|d\|_W^2&:=\E\sum_{i=0}^{N-1}\tau
\Big[(R+\tau S_{i+1})|\tilde d_i|_p^2
 +(R+\tau\bar S_{i+1})\bar d_i^2\Big].
\end{aligned}
\]
By~(i), $\mathcal R_{\pi,p}(u)=\|d(u)\|_W^2$. For the two controls,
\[
d_i(\ap u)-d_i(u^{\mathrm{FB}})
=(k_i-P_i/R)\tilde\delta_i
+(\bar k_i-\Pi_i/R)\bar\delta_i\mathbf 1-R^{-1}e_i.
\]
The reverse triangle inequality and the bounds on $\delta$ and $e$ give
\[
\begin{aligned}
\left|\sqrt{\mathcal R_{\pi,p}(\ap u)}-\sqrt{B_{\pi,p}}\right|
&\le \Delta_k\|\delta\|_W+R^{-1}\|e\|_W\\
&\le \sqrt{w_{\max}T}\,R^{-1}(1+T\Delta_k)
       \sqrt{10\,\Epip(\ap\Theta)}.
\end{aligned}
\]

\emph{(iii)} By~(i), $B_{\pi,p}=\|d(u^{\mathrm{FB}})\|_W^2\le Tw_{\max}\Delta_k^2\max_i\E|\ex X_i|_p^2$,
and the maximum is bounded uniformly in $p$ and $\tau$ because $\varphi_i,\bar\varphi_i\in[0,1]$
and $\E|\xi_i|_p^2=(\sigma_0^2+h_0^2\nu(E))\tau$. Since $\ell=\varphi_i+\tau P_i/R\ge0$ and $D_i\ge1$, the Riccati recursions give
\[
\begin{aligned}
0\le P_i&\le(1+|\kappa|\tau)P_{i+1}+\tau Q,\\
0\le S_i&\le(1+|\kappa|\tau)^2S_{i+1}+\tau Q.
\end{aligned}
\]
Iterating these inequalities gives
\[
0\le S_i,P_i\le K:=e^{2|\kappa|T}(Q_T+TQ).
\]
Write $F_\tau(s):=\tau Q+\ell^2Rs/(R+\tau s)$, so $S_i=F_\tau(S_{i+1})$.
On $[0,K]$, both Riccati updates have the uniform expansion
\[
s+\tau\bigl(Q-2\kappa s-s^2/R\bigr)+O(\tau^2).
\]
Thus $P_i=F_\tau(P_{i+1})+O(\tau^2)$. Since
\[
0\le F_\tau'(s)=\frac{\ell^2R^2}{(R+\tau s)^2}
\le(1+|\kappa|\tau)^2\le1+C\tau,
\]
we obtain
\[
|S_i-P_i|\le(1+C\tau)|S_{i+1}-P_{i+1}|+C\tau^2.
\]
Discrete Gr\"onwall, with $S_N=P_N=Q_T$, gives
\[
\max_{0\le i\le N}|S_i-P_i|\le C\tau.
\]
The gain formula gives $k_i=S_{i+1}/R+O(\tau)$, and the Riccati expansion gives
$P_i=P_{i+1}+O(\tau)$. Hence
\[
k_i-\frac{P_i}{R}
=\frac{S_{i+1}-P_{i+1}}{R}+O(\tau)
=O(\tau)
\]
uniformly in $i$. Applying the same argument with $\kappa+\bar a,\bar Q,\bar Q_T$ gives
$\bar k_i-\Pi_i/R=O(\tau)$. Thus $\Delta_k=O(\tau)$, and the bound at the start of (iii) gives
$B_{\pi,p}=O(\tau^2)$ uniformly in $p$.
\end{proof}

\subsection{Coefficient bounds for the affine benchmarks}\label{B:app:benchmark-admiss}

\begin{remark}[Coefficient bounds for the state-law benchmark]\label{B:rem:direct-benchmark-mono}
For the benchmark of Section~\ref{B:subsec:benchmark-protocol}, set
$a:=|G\,\E\delta X|$ and $y:=|\Gstar\E\delta Y|$, and
\[
M:=\begin{pmatrix}\kappa_{\mathrm B}+\bar c&-|\bar a|\sqrt{c_G}/2\\[2pt]
-|\bar a|\sqrt{c_G}/2&\kappa_{\mathrm F}\end{pmatrix}.
\]
By Sylvester's criterion, $M$ is positive definite if and only if
\[
\kappa_{\mathrm B}+\bar c>0,\qquad (\kappa_{\mathrm B}+\bar c)\,\kappa_{\mathrm F}>\tfrac14\bar a^2c_G .
\]
Since $|G^+G\,\E\delta X|\le\sqrt{c_G}\,a$, the mean-field cross-term is bounded by
$|\bar a|\sqrt{c_G}\,ay$, and the interior monotonicity pairing is at most
\[
-\kappa_{\mathrm B}\bigl(\E|G\delta X|^2-a^2\bigr)
-\kappa_{\mathrm F}\bigl(\E|\Gstar\delta Y|^2-y^2
+\E|\Gstar\delta Z|_F^2+\E\|\Gstar\delta U\|_{\Lnu}^2\bigr)
-(a,y)M(a,y)^{\!\top}.
\]
When $M$ is positive definite, the last term is at most $-\lambda_{\min}(M)(a^2+y^2)$, so
Assumption~\ref{B:ass:mono} holds with
\[
\beta_1=\min\{\kappa_{\mathrm F},\lambda_{\min}(M)\},
\qquad
\beta_2=\min\{\kappa_{\mathrm B},\lambda_{\min}(M)\},
\qquad
\alpha=1,
\]
provided $\kappa_{\mathrm B}\ge0$, with $\kappa_{\mathrm B}>0$ when $m>n$, and $\bar g\ge0$: the terminal pairing is
$\E|G\delta X_T|^2+\bar g\,|G\,\E\delta X_T|^2$. Every reported setting has $\kappa_{\mathrm B}>0$,
$\bar g>0$ and satisfies the criterion. For the reported feedback choices
$\kappa_{\mathrm F}=\kappa_{\mathrm B}=\rv{apocert.betaOne},\rv{apocert.betaStrong},\rv{apocert.betaWeak}$
at $c_G=\nu(E)=1$, the corresponding two-decimal lower bounds for $\lambda_{\min}(M)$ are
$\monoBaseVal$, $\monoStrongVal$ and $\monoWeakVal$, respectively; each vector rank case is
checked with its own $c_G=s_{\min}(G)^{-2}$.

For the remaining hypotheses, the mean map satisfies $|m_X(\mu)-m_X(\varrho)|\le W_2(\mu,\varrho)$
and $\|G^+G\|\le1$, so Assumption~\ref{B:ass:lip} holds with
\[
L:=\max\bigl\{\kappa_{\mathrm F}\|G\|,\ \kappa_{\mathrm B}\|G\|,\ |\bar a|+|\bar c|\,\|G\|,\
\|G\|\max\{1,|\bar g|\}\bigr\}.
\]
With the tuple fixed, only $b,f,g$ change with the law, so Assumption~\ref{B:ass:offdiag}
holds with
\[
L_\mu:=\max\bigl\{|\bar a|,\ |\bar c|\,\|G\|,\ |\bar g|\,\|G\|\bigr\}.
\]
Time independence gives Convention~\ref{B:conv:nodewise} at every node and $\omega\equiv0$.
\end{remark}

\begin{remark}[Dissipation matrix for the full-tuple-law benchmark]\label{B:rem:fulllaw-dissipation}
The direct construction of Remark~\ref{B:rem:direct-benchmark-mono} extends to the
full-tuple benchmark of Section~\ref{B:ssec:fulllaw-setup}: with the component means ordered
$v=(\overline{\delta X},\overline{\delta Y},\overline{\delta Z},\overline{\delta U})$,
the $4\times4$ dissipation matrix is
\[
M=\begin{pmatrix}
\kappa_{\mathrm B} & -\lambda_X\sqrt{c_G}/2 & 0 & 0\\
-\lambda_X\sqrt{c_G}/2 & \kappa_{\mathrm F}-\lambda_Y & -\lambda_Z/2 & -\lambda_U\nu(E)/2\\
0 & -\lambda_Z/2 & \kappa_{\mathrm F} & 0\\
0 & -\lambda_U\nu(E)/2 & 0 & \kappa_{\mathrm F}\,\nu(E)
\end{pmatrix},
\]
the factors $\nu(E)$ coming from the $L^2(\nu)$ pairing in the $U$-component. For
$\kappa_{\mathrm F},\nu(E)>0$ the diagonal $(Z,U)$ block is positive definite, and eliminating it shows
that $M$ is positive definite if and only if the Schur complement
\[
\begin{pmatrix}
\kappa_{\mathrm B} & -\lambda_X\sqrt{c_G}/2\\
-\lambda_X\sqrt{c_G}/2 & \kappa_{\mathrm F}-\lambda_Y-\dfrac{\lambda_Z^2+\lambda_U^2\nu(E)}{4\kappa_{\mathrm F}}
\end{pmatrix}
\]
is, that is,
\[
\kappa_{\mathrm B}>0,\qquad
\kappa_{\mathrm B}\Bigl[\kappa_{\mathrm F}-\lambda_Y-\frac{\lambda_Z^2+\lambda_U^2\nu(E)}{4\kappa_{\mathrm F}}\Bigr]
>\frac{\lambda_X^2c_G}{4}.
\]
For the parameters of Section~\ref{B:ssec:fulllaw-setup},
$\kappa_{\mathrm F}=\kappa_{\mathrm B}=\nu(E)=c_G=1$ and
$(\lambda_X,\lambda_Y,\lambda_Z,\lambda_U)=\flLamVal$, the spectrum of $M$ is
$\{1,1,(38-\sqrt{17})/40,(38+\sqrt{17})/40\}$. Thus
\[
\lambda_{\min}(M)=\frac{38-\sqrt{17}}{40}\ge\monoFullLawVal.
\]
For any two diagonal inputs,
\[
\begin{aligned}
&\E\bigl[
\delta b\,\delta Y+\delta\sigma\,\delta Z
+\delta h\,\delta U-\delta f\,\delta X
\bigr]\\
&\qquad=-\E\|\delta\Theta-v\|_{\mathcal H}^{2}-v^{\!\top}Mv\\
&\qquad\le -\lambda_{\min}(M)\,\E\|\delta\Theta\|_{\mathcal H}^{2}.
\end{aligned}
\]
Since $g(x,\rho)=x$, Assumption~\ref{B:ass:mono} therefore holds with the choice
\[
\alpha=1,\qquad \beta_1=\beta_2=\monoFullLawVal.
\]
\end{remark}

\subsection{Numerical reference for the state-nonlinear benchmark}\label{B:app:snlref}
For the state-nonlinear benchmark, write $\ex Y_i=v_i(\ex X_i)$. With
$f(x)=g(x)=x+\tfrac12\tanh x$, the backward recursion is
\[
v_N=g,\qquad v_i(x)=H_i\bigl(x-\tau\tanh v_i(x)\bigr)+\tau f(x),
\]
where
\[
H_i(s)=\E\!\left[v_{i+1}\!\left(s+\sigma_0\sqrt{\tau}\,\xi+h_0(K-\tau)\right)\right],
\]
$\sigma_0=h_0=\tfrac12$, and $\xi\sim\mathcal N(0,1)$ and $K\sim\mathrm{Poisson}(\tau)$ are
independent. For $i<N$, we tabulate $v_i$ on a uniform grid on $[-L,L]$ with spacing $\Delta x$,
using the settings in Table~\ref{B:tab:snlref}.

Let $I_iv_i$ denote its monotone piecewise-cubic Hermite interpolant on this interval. We extend it
affinely by
\[
\widetilde v_i(x)=
\begin{cases}
v_i(-L)+A_i(x+L), & x<-L,\\
I_iv_i(x), & |x|\le L,\\
v_i(L)+A_i(x-L), & x>L,
\end{cases}
\]
where
\[
A_N=1,\qquad A_i=A_{i+1}+\tau.
\]
The implicit roots defining $v_i$ are solved to the root residual of
Table~\ref{B:tab:snlref}.

Let $(\xi_q,w_q)_{q=1}^{Q}$ denote the Gauss--Hermite rule normalised for a standard normal
variable, and let $p_k=\Prob(K=k)$. With $\widetilde v_N=g$, we evaluate
\[
H_i(s)\approx
\sum_{k=0}^{K_{\max}}p_k
\sum_{q=1}^{Q}w_q\,
\widetilde v_{i+1}
\!\left(
s+\sigma_0\sqrt{\tau}\,\xi_q
+h_0(k-\tau)
\right),
\]
where $K_{\max}$ is the smallest integer such that
$\E\bigl[(1+K)^2\mathbf 1_{\{K>K_{\max}\}}\bigr]$ is below the Poisson tail of
Table~\ref{B:tab:snlref}.

\section*{Declarations}
\noindent Siemens funded C. Garry's MSc studies at Durham University and employs
C. Garry. Siemens had no role in the conception, analysis, or preparation of this work.
The authors declare no other competing interests.

\medskip\noindent\textbf{Code and data availability.} The code, data and reproduction guide
are available at\\
\url{https://github.com/cgarryZA/mvfbsdej-aposteriori}.

\end{document}